\documentclass{amsart}
\usepackage{amssymb}

\newtheorem{thm}{Theorem}[section]
\newtheorem{prop}[thm]{Proposition}
\newtheorem{cor}[thm]{Corollary}

\newtheorem{rema}[thm]{Remark}
\newtheorem{defn}[thm]{Definition}
\newtheorem{conj}[thm]{Conjecture}
\newtheorem{ex}[thm]{Example}

\numberwithin{equation}{section}

\newcommand{\Lx}[0]{x^{j + 1} \frac{\partial}{\partial x} + \Bigl(\frac{j + 1}{2}\Bigr) x^j \Bigl( \varphi^+ \frac{\partial}{\partial \varphi^+} + \varphi^- \frac{\partial}{\partial \varphi^-}\Bigr)}

\newcommand{\Jx}[0]{x^j\Bigl(\varphi^+\frac{\partial}{\partial \varphi^+} - \varphi^- \frac{\partial}{\partial \varphi^-}\Bigr)}

\newcommand{\Gx}[0]{x^j \Bigl( \frac{\partial}{\partial \varphi^\pm} - \varphi^\mp \frac{\partial}{\partial x}\Bigr) \pm jx^{j-1} \varphi^+ \varphi^- \frac{\partial}{\partial \varphi^\pm}}

\newcommand{\dz}[0]{\frac{\partial}{\partial z}}
\newcommand{\nor}[0]{\Upsilon}
\newcommand{\sou}[0]{\Delta}

\newcommand{\Z}{\mathbb{Z}_+}

\begin{document}

\title[The $N=2$ moduli space]{The moduli space of $N=2$ super-Riemann spheres with tubes}

\author{Katrina Barron}
\address{Department of Mathematics, University of Notre Dame,
Notre Dame, IN 46556}
\email{kbarron@nd.edu}


\date{February 27, 2007}


\begin{abstract}
Within the framework of complex supergeometry and motivated by two-dimensional genus-zero holomorphic  $N = 2$  superconformal field theory, we define the moduli space of $N=2$  super-Riemann spheres with oriented and ordered half-infinite tubes (or equivalently, oriented and ordered punctures, and local superconformal coordinates vanishing at the punctures), modulo $N=2$ superconformal equivalence.  We develop a formal theory of infinitesimal $N = 2$ superconformal transformations based on a representation of the $N=2$ Neveu-Schwarz algebra in terms of superderivations.    In particular, via these infinitesimals we present the Lie supergroup of $N=2$ superprojective transformations of the $N=2$ super-Riemann sphere.  We give a reformulation of the moduli space in terms of these infinitesimals.  We introduce generalized $N=2$ super-Riemann spheres with tubes and discuss some group structures associated to certain moduli spaces of both generalized and non-generalized $N=2$ super-Riemann spheres.  We  define an action of the symmetric groups on the moduli space.  Lastly we discuss the nonhomogeneous (versus homogeneous) coordinate system associated to $N=2$ superconformal structures and the corresponding results in this coordinate system. 
\end{abstract}

\maketitle

\section{Introduction}

In this work, we give a detailed study of the differential supergeometry underlying two-dimensional genus-zero holomorphic $N=2$ superconformal field theory. 

Conformal field theory (or more specifically, string theory) and related theories, including superconformal field theories, (cf. \cite{BPZ}, \cite{Fd}, \cite{FS}, \cite{V}, and \cite{S}) are the most  promising attempts at developing a physical theory that combines all fundamental interactions of
particles, including gravity.  The geometry of (super)conformal field theory extends the use of Feynman diagrams, describing the interactions of point particles whose propagation in time sweeps out a line in space-time, to one-dimensional particles called ``strings" (or higher-dimensional ``superstrings") whose propagation in time sweeps out a two-dimensional surface (or higher-dimensional supersurface) called the ``worldsheet".  For two-dimensional genus-zero holomorphic (super)conformal field theory, algebraically, these interactions can be described by products of vertex operators or more precisely,  by vertex operator (super)algebras  (cf. \cite {Bo}, \cite{FLM}, \cite{KW},  \cite{B-vosas}, \cite{B-iso}).   

In \cite{H-thesis} and \cite{H-book}, motivated by the geometric notions arising in conformal field theory, Huang gives a precise geometric interpretation of the notion of vertex operator algebra by considering the geometric structure consisting of the moduli space of genus-zero Riemann surfaces with oriented and ordered half-infinite tubes (which are conformally equivalent to Riemann spheres with oriented and ordered punctures and local coordinates vanishing at the punctures), modulo conformal equivalence, together with the operation of sewing two such surfaces.  This work was then extended by the author in \cite{B-announce}--\cite{B-iso} to the notion of $N=1$ supergeometric vertex operator superalgebra motivated by the physical model of $N=1$ superconformal field theory and the notion of a superstring whose propagation in time sweeps out a supersurface \cite{Fd}, \cite{D}, \cite{Ro}, \cite{CR}, thereby giving a precise geometric interpretation of the notion of $N=1$ Neveu-Schwarz vertex operator superalgebra.

This rigorous foundation for the correspondence between the algebraic and geometric aspects of two-dimensional genus-zero holomorphic ($N=1$ super)conformal field theory has proved useful in furthering both the algebraic and geometric aspects of (super)conformal field theory.  On the one hand, it is much easier to rigorously construct many aspects of (super)conformal field theories and study many of their properties using the algebraic formulation of vertex operator (super)algebra (or equivalently, chiral algebra) (e.g., \cite{Wi}, \cite{BPZ}, \cite{Za}, \cite{FLM}, \cite{KT}, \cite{FZ}, \cite{FFR}, \cite{DL}, \cite{FF}, \cite{DMZ}, \cite{Wa}, \cite{Zh1}, \cite{Zh2}).  On the other hand,  the geometry of (super)conformal field theory can give insight and provide  tools useful for studying the algebraic aspects of the theory, for  example: giving rise to general results in Lie theory \cite{BHL1};  giving the necessary insight for developing a theory of tensor  products for vertex operator algebras \cite{HL1}, \cite{HL4}--\cite{HL7}, \cite{H-tensor}; giving rise to change of variables formulas for vertex operator algebras \cite{H-book} and $N=1$ Neveu-Schwarz vertex operator superalgebras  \cite{B-change}; and giving rise to constructions in orbifold conformal field theory \cite{BDM}, \cite{H-orbifold}, \cite{BHL2}.  But one of the most important applications arises {}from the fact that  this rigorous development of the differential geometric foundations of (super)conformal field theory (in particular an analytic development of the moduli space of (super-)Riemann surfaces and a sewing operation) is necessary for the construction of (super)conformal field theory in the sense of Segal \cite{S} and Kontsevich.  In fact the work of Huang in \cite{H-thesis}, \cite{H-book} along with \cite{H-tensor}, \cite{H1998}--\cite{H2005-2} solves the problem of constructing holomorphic genus-zero (weakly) conformal field theories {}from certain representations of vertex operator algebras.  

The purpose of this paper is to extend certain geometric aspects of Huang's work \cite{H-book} and previous work by the author \cite{B-thesis}, \cite{B-memoir} to the analogous $N=2$ supergeometry underlying two-dimensional genus-zero holomorphic $N = 2$ superconformal field theory.  That is to develop the differential supergeometric foundations underlying the ``worldsheet" approach to genus-zero holomorphic two-dimensional $N=2$ superconformal field theory necessary to give a rigorous correspondence between the geometric and algebraic aspects of the theory. 

In \cite{Fd}, Friedan describes the extension of the physical model of conformal field theory to that of $N=1$ superconformal field theory and the notion of a superstring whose propagation in time sweeps out a supersurface.  Whereas conformal field theory attempts to describe the interactions of bosons, superconformal field theory attempts to describe the interactions of bosons paired with fermions, adding additional boson-fermion symmetries.  $N=2$ superconformal field theory explores the theory of conformal fields if, rather than one boson-fermion symmetry is present as in $N=1$ theories, two symmetries are realized (cf. \cite{FMS}, \cite{LVW}, \cite{Ge}, \cite{Schwarz}, \cite{Wa}).  This theory has gained much interest due to the phenomenon of mirror symmetry present (cf., \cite{LVW}, \cite{GP}, \cite{CLS}, \cite{ALR}, \cite{CdGP}, \cite{FBC}, \cite{DHVW1}, \cite{DHVW2}, \cite{LY}, \cite{LLY1}--\cite{LLY4}, \cite{Gi}, \cite{Y}, \cite{GY}, \cite{PVY}).  However, a rigorous foundation for the correspondence between the underlying worldsheet differential geometry and the algebra of $N=2$ Neveu-Schwarz vertex operator superalgebras has not yet been realized.

The natural setting for $N=2$ superconformal field theory lies in $N = 2$ complex supergeometry \cite{FMS}, \cite{DRS}, \cite{D}.  This is the geometry of manifolds over a complex Grassmann algebra (i.e., complex supermanifolds) where there is one ``even'' dimension corresponding to the one bosonic component and two ``odd'' dimensions corresponding to the $N = 2$ fermionic components.  In this paper, within the framework of $N=2$ complex supergeometry and motivated by two-dimensional genus-zero holomorphic  $N = 2$  superconformal field theory, we define the moduli space of $N=2$  super-Riemann spheres with oriented and ordered half-infinite tubes (or equivalently, oriented and ordered punctures, and local $N=2$ superconformal coordinates vanishing at the punctures), modulo $N=2$ superconformal equivalence.  Physically, one can think of each $N=2$ super-Riemann sphere as representing some superstring interaction.   

We note that since there is yet no uniformization theorem for $N=2$ super-Riemann surfaces with genus-zero compact body, for much of this paper we are working with actual $N=2$ super-Riemann spheres with tubes, rather than more general genus-zero $N=2$ super-Riemann surfaces.  In Section 5, we present Conjecture 5.2 which states that any genus-zero $N=2$ super-Riemann surface is $N=2$ superconformally equivalent to the $N=2$ super-Riemann sphere.  If true, this conjecture implies that the results of this paper extend to the more general moduli space of $N=2$ genus-zero super-Riemann surfaces with tubes.   However, we do not need this conjecture for the results of this paper as stated. 

In Section 6, we develop a formal theory of infinitesimal $N=2$ superconformal transformations which are superderivations.  We show that any local $N=2$ superconformal coordinate can be expressed in terms of exponentials  of these superderivations.  This theory of infinitesimal $N=2$ superconformal transformations is the main machinery that will be used to reformulate the moduli space as we do in Section 8.  It is this formulation that will be needed in subsequent work to define a sewing operation on the moduli space and carry out the program of giving a rigorous geometric interpretation of $N=2$ Neveu-Schwarz vertex operator superalgebras.   

In Section 7, we discuss some of the group structures arising {}from this formulation of local coordinates in terms of infinitesimals.  We discuss the $N=2$ Neveu-Schwarz Lie superalgebra and point out that the infinitesimal $N=2$ superconformal transformations give a representation of the $N=2$ Neveu-Schwarz algebra  in terms of superderivations (cf. \cite{Ki}).  In particular, via these infinitesimals we present the Lie supergroup of $N=2$ superprojective transformations of the $N=2$ super-Riemann sphere.   In reviewing the literature, we found many claims of presentations of this group of automorphisms of the $N=2$ super-Riemann sphere, for example \cite{C}, \cite{Ki}, \cite{Melzer},  \cite{Sc}, \cite{Nogueira} and \cite{BL}.  However, none of these provide a correct presentation of the full group of $N=2$ superprojective transformations that are $N=2$ superconformal in the sense that we deal with in this paper.  We give a detailed analysis of each presentation and concrete counterexamples to show that the ``$N=2$ superprojective transformations" found in the literature cited above do not give the group of automorphisms of the $N=2$ super-Riemann sphere in the usual sense.   The calculus of $N=2$ superconformal functions we develop using exponentials of infinitesimals proves useful in unraveling the problem of determining the correct general form of the $N=2$ superprojective transformations.  

In Section 8, we give the main result of this paper -- a reformulation of the moduli space of $N=2$ super-Riemann spheres with tubes in terms of the infinitesimal $N=2$ superconformal transformations.  In Section 9, we introduce generalized $N=2$ super-Riemann spheres with tubes and discuss some group structures associated to the moduli space of both generalized and non-generalized $N=2$ super-Riemann spheres with one outgoing tube and one incoming tube.   There is no complex Lie group corresponding to the $N=2$ Neveu-Schwarz algebra (just as there is no complex Lie group corresponding to the Virasoro algebra \cite{L}).  However, the group structures we present here associated to the moduli space of $N=2$ super-Riemann spheres lay the groundwork for showing that there are group structures corresponding to certain infinite-dimensional sub-algebras of the $N=2$ Neveu-Schwarz algebra, and there is a partial monoid whose tangent space at the identity is the $N=2$ Neveu-Schwarz algebra.  This is analogous to the non-super and the $N=1$ superconformal case as studied in  \cite{H-book} and  \cite{B-memoir}, respectively.  

In Section 10, we  define an action of the symmetric groups on the moduli space.  In later work we will put a sewing operation on the moduli space of $N=2$ super-Riemann spheres with tubes, and this together with the action of the symmetric group presented here will give the moduli space the structure of a partial operad \cite{operads}, in analogy to the non-super case \cite{HL2}, \cite{HL3} and $N=1$ super case  \cite{B-memoir}.  Finally, in Section 11 we discuss the nonhomogeneous (versus homogeneous) coordinate systems associated to $N=2$ superconformal structures (cf. \cite{C}, \cite{DRS}, \cite{Nogueira}) and the corresponding results in this coordinate system.

Here we note some immediate applications of this work.  In analogy to the situation for geometric vertex operator algebras \cite{H-book} and $N=1$ supergeometric vertex operator superalgebras \cite{B-iso}, we expect to use the results of this paper to establish the notion of $N=2$ supergeometric vertex operator superalgebra and to establish the isomorphism between the category of such objects and the category of $N=2$ Neveu-Schwarz vertex operator superalgebras.  This rigorous correspondence between the geometric and algebraic aspects of $N=2$ Neveu-Schwarz vertex operator superalgebras is necessary for the construction of $N=2$ superconformal field theory in the sense of \cite{S} following the program of Huang in \cite{H-thesis}--\cite{H2005-2}.  

We take this opportunity to point out that in \cite{HeKac}, Heluani and Kac formulate an axiomatic notion of what they call ``strongly conformal $N_K = N$ SUSY vertex superalgebras" and in  \cite{He}, Heluani  gives certain change of variables formulas (near zero) for these algebras.  For $N=1$ these algebras correspond to ``$N=1$ Neveu-Schwarz vertex operator superalgebras" as formulated and studied previously by the author in  \cite{B-announce}--\cite{B-vosas}, and change of variable formulas for such algebras under any $N=1$ superconformal change of variables (not just near zero) had  been formulated and proved by the author in \cite{B-change}; the subtleties entering into the formulation and the proof of the change of variable formulas were explained in detail in \cite{B-change}.   It was actually the rigorous treatment of the correspondence between certain geometric and algebraic aspects of $N=1$ superconformal field theory as developed in \cite{B-thesis} and  \cite{B-memoir} that led the author to formulate and study the notion of $N=1$ Neveu-Schwarz vertex operator superalgebra with odd formal variables and to introduce (in \cite{B-announce} and \cite{B-thesis}; see also \cite{B-vosas}) such properties as the $G(-1/2)$-derivative property, the Jacobi identity and the commutator formula with odd variables, supercommutativity and associativity with odd variables, and skew-supersymmetry.   The author's $N=1$ axioms and properties mentioned here generalize routinely to $N \geq 1$, treated in \cite{HeKac} and \cite{He} (except for the Jacobi identity).

In fact, we can further motivate the present paper by recalling that actually, it was the substantial problem of formulating and studying not just the algebraic structures mentioned in the previous paragraph,
but rather the superconformal geometry underlying such structures, that was the main focus of \cite{B-announce}--\cite{B-change}.  Some of the main ingredients of this geometric study were the theory of the
sewing operation, in full generality, for $N=1$ genus zero super-Riemann surfaces with punctures and local coordinates; the intricate but essential proof of the convergence, necessary for handling the sewing with general local coordinates; the formulation of a {\it geometric} notion of ``$N=1$ Neveu-Schwarz vertex operator superalgebra''; the proof of the equivalence between this notion and the algebraic notion; as a consequence of this theory, a general change of variables formula including changes of coordinates not just near zero but also near infinity and more generally, on an annulus; and the solution of the problems of when and how one obtains isomorphic vertex operator superalgebras under general changes of variables.  The present paper is the first in a series of papers extending this entire algebraic and geometric program in the $N=1$ case to the $N=2$ case.

Acknowledgments:  The author thanks the referee for helpful suggestions.

Notational conventions: $\mathbb{F}$ denotes a field of characteristic zero, $\mathbb{N}$ denotes the nonnegative integers, $\mathbb{R}_+$ denotes the positive real numbers, $\mathbb{Z}_+$ denotes the positive integers, and $\mathbb{Z}_2$ denotes the integers modulo 2.

\section{Superanalytic functions}\label{preliminaries-section}

In this section, we recall the notion of superalgebra, Grassmann algebra and superanalytic function following \cite{B-memoir}, \cite{D}, \cite{Bc1}, \cite{Ro}.  

Let $\mathbb{F}$ be  a field of characteristic zero.  For a $\mathbb{Z}_2$-graded vector space $X = X^0 \oplus X^1$, define the {\it sign function} $\eta$ on the homogeneous subspaces of $X$ by $\eta(x) = j$, for $x \in X^j$ and $j \in \mathbb{Z}_2$.  If $\eta(x) = 0$, we say that $x$ is {\it even}, and if $\eta(x) = 1$, we say that $x$ is {\it odd}.    

A {\it superalgebra} is an (associative) algebra $A$ (with identity $1 \in A$), such that: (i) $A$ is a $\mathbb{Z}_2$-graded algebra; (ii) $ab = (-1)^{\eta(a)\eta(b)} ba$ for $a,b$ homogeneous in $A$.
Note that when working over a field of characteristic zero or of characteristic greater than two, property (ii), supercommutativity, implies that the square of any odd element is zero.

A $\mathbb{Z}_2$-graded vector space $\mathfrak{g} = \mathfrak{g}^0 \oplus \mathfrak{g}^1$ is said to be a {\it Lie superalgebra} if it has a bilinear operation $[\cdot,\cdot]$ on $\mathfrak{g}$ such that for $u,v$ homogeneous in $\mathfrak{g}$: (i) $\; [u,v] \in {\mathfrak g}^{(\eta(u) + \eta(v))\mathrm{mod} \; 2}$;  
(ii) skew symmetry holds $[u,v] = -(-1)^{\eta(u)\eta(v)}[v,u]$; (iii) the following Jacobi identity holds 
\[(-1)^{\eta(u)\eta(w)}[[u,v],w] + (-1)^{\eta(v)\eta(u)}[[v,w],u]+ \; (-1)^{\eta(w)\eta(v)}[[w,u],v] = 0. \]

\begin{rema}\label{envelope} 
{\em Given a Lie superalgebra $\mathfrak{g}$ and a superalgebra $A$, the space $(A^0 \otimes \mathfrak{g}^0) \oplus (A^1 \otimes \mathfrak{g}^1)$ is a Lie algebra with bracket given by
\begin{equation}\label{super-v.-non} 
[au , bv] = (-1)^{\eta(b) \eta(u)} ab [u,v]
\end{equation}
for $a,b \in A$ and $u,v \in \mathfrak{g}$ homogeneous (with obvious notation), where in (\ref{super-v.-non}) we have suppressed the tensor product symbol.  Note that the bracket on the left-hand side 
of (\ref{super-v.-non}) is a Lie algebra bracket, and the bracket on the right-hand side is a Lie superalgebra bracket.  The Lie algebra $(A^0 \otimes \mathfrak{g}^0) \oplus (A^1 \otimes \mathfrak{g}^1)$ is called the {\it $A$-envelope of $\mathfrak{g}$}.  }
\end{rema}

For any $\mathbb{Z}_2$-graded associative algebra $A$ and for $u,v \in A$ of homogeneous sign, we can define $[u,v] = u v - (-1)^{\eta(u)\eta(v)} v u$, making $A$ into a Lie superalgebra.  The algebra of endomorphisms of $A$, denoted $\mbox{End} \; A$, has a natural $\mathbb{Z}_2$-grading induced {}from that of $A$, and defining $[X,Y] = X Y - (-1)^{\eta(X) \eta(Y)} Y X$ for $X,Y$ homogeneous in $\mbox{End} \; A$, this gives $\mbox{End} \; A$ a Lie superalgebra structure.  An element $D \in 
(\mbox{End} \; A)^j$, for $j \in \mathbb{Z}_2$, is called a {\it superderivation of sign $j$} (denoted $\eta(D) = j$) if $D$ satisfies the {\it super-Leibniz rule}
\begin{equation}\label{leibniz} 
D(uv) = (Du)v + (-1)^{\eta(D) \eta(u)} uDv  
\end{equation} 
for $u,v \in A$ homogeneous.

Let $V$ be a vector space.  Then the exterior algebra generated by $V$, denoted  $\bigwedge (V)$, has the structure of a superalgebra.  Fix $V_L$ to be  an $L$-dimensional vector space 
over $\mathbb{C}$ with basis $\{\zeta_1,\zeta_2, \ldots, \zeta_L\}$ for $L \in \mathbb{N}$ such that $V_L \subset V_{L+1}$.  We denote $\bigwedge(V_L)$ by $\bigwedge_L$ and call this the {\it Grassmann 
algebra on $L$  generators}.  In other words, {}from now on we will consider the Grassmann algebras to have a fixed sequence of generators.  Note that $\bigwedge_L \subset \bigwedge_{L+1}$, and taking the direct limit as $L \rightarrow \infty$, we have the {\it infinite Grassmann algebra} denoted by $\bigwedge_\infty$.  Then $\bigwedge_L$ and $\bigwedge_\infty$ are the associative algebras over $\mathbb{C}$ with generators $\zeta_j$, for $j = 1, 2, \dots, L$ and $j = 1, 2, \dots$, respectively, and with relations
\[\zeta_j \zeta_k = - \zeta_k \zeta_j, \qquad \qquad \zeta_j^2 = 0 .\]    
Note that $\mathrm{dim}_{\mathbb{C}} \; \bigwedge_L = 2^L$, and if $L = 0$, then $\bigwedge_0 = \mathbb{C}$. We use the notation $\bigwedge_*$ to denote a Grassmann algebra, finite or infinite.  
The reason we take $\bigwedge_*$ to be over $\mathbb{C}$ is that we will mainly be interested in complex supergeometry.  However, formally, we could just as well have taken $\mathbb{C}$ to be any field of characteristic zero.  
 
Let 
\begin{eqnarray*}
I_L \! \!  &=&  \! \! \bigl\{ (i) = (i_1, i_2, \ldots, i_{2n}) \; | \; i_1 < i_2 < \cdots < i_{2n}, \; i_l \in \{1, 2, \dots, L\}, \; n \in \mathbb{N} \bigr\}, \\ 
J_L  \! \! &=&  \! \! \bigl\{(j) = (j_1, j_2, \ldots, j_{2n + 1}) \; | \; j_1 < j_2 < \cdots < j_{2n + 1}, \; j_l \in \{1, 2, \dots, L\}, \;  n \in \mathbb{N} \bigr\},
\end{eqnarray*}
and $K_L = I_L \cup J_L$.  Let
\begin{eqnarray*}
I_\infty \! \! &=& \! \! \bigl\{(i) = (i_1, i_2, \ldots, i_{2n})\; | \; i_1 < i_2 < \cdots < i_{2n}, \; i_l \in \mathbb{Z}_+, \; n \in \mathbb{N} \bigr\}, \\
J_\infty \! \! &=& \! \! \bigl\{(j) = (j_1, j_2, \ldots, j_{2n + 1})\; | \; j_1 < j_2 < \cdots < j_{2n + 1}, \; j_l \in \mathbb{Z}_+, \; n \in \mathbb{N} \bigr\},
\end{eqnarray*}
and $K_\infty = I_\infty \cup J_\infty$.  We use $I_*$, $J_*$, and $K_*$ to denote $I_L$ or $I_\infty$, $J_L$ or $J_\infty$, and $K_L$ or $K_\infty$, respectively.  Note that $(i) = (i_1,\dots,i_{2n})$ for $n = 0$ is in $I_*$, and we denote this element by $(\emptyset)$.  The $\mathbb{Z}_2$-grading of $\bigwedge_*$ is given explicitly by
\begin{eqnarray*}
\mbox{$\bigwedge_*^0$} \! &=& \! \Bigl\{a \in \mbox{$\bigwedge_*$} \; \big\vert \; a = \sum_{(i) \in I_*} a_{(i)}\zeta_{i_{1}}\zeta_{i_{2}} \cdots \zeta_{i_{2n}}, \; a_{(i)} \in \mathbb{C}, \; n \in \mathbb{N} \Bigr\}\\ 
\mbox{$\bigwedge_*^1$} \! &=& \! \Bigl\{a \in \mbox{$\bigwedge_*$} \; \big\vert \; a = \sum_{(j) \in J_*} a_{(j)}\zeta_{j_{1}}\zeta_{j_{2}} \cdots \zeta_{j_{2n + 1}}, \; a_{(j)} \in \mathbb{C}, \; n \in \mathbb{N}
\Bigr\} . 
\end{eqnarray*} 
Note that $a^2 = 0$ for all $a \in \bigwedge_*^1$.

We can also decompose $\bigwedge_*$ into {\it body}, $(\bigwedge_*)_B = \{ a_{(\emptyset)} \in \mathbb{C} \}$,  and {\it soul} 
\[(\mbox{$\bigwedge_*$})_S \; = \; \Bigl\{a \in \mbox{$\bigwedge_*$} \; \big\vert \;  a = \! \! \! \sum_{ \begin{tiny} \begin{array}{c}
(k) \in K_*\\
k \neq (\emptyset)
\end{array} \end{tiny}} \! \! \!
a_{(k)} \zeta_{k_1} \zeta_{k_2} \cdots \zeta_{k_n}, \; a_{(k)} \in \mathbb{C} 
\Bigr\}\] 
subspaces such that $\bigwedge_* = (\bigwedge_*)_B \oplus (\bigwedge_*)_S$.  For $a \in \bigwedge_*$, we write $a = a_B + a_S$ for its body and soul decomposition.  Note that for all $a \in \bigwedge_L$, we have $a_S^{L + 1} = 0$. 

For $n \in \mathbb{N}$, we introduce the notation $\bigwedge_{*>n}$ to denote a finite Grassmann algebra $\bigwedge_L$ with $L > n$ or an infinite Grassmann algebra.  We will use the corresponding index notations for the corresponding indexing sets $I_{*>n}, J_{*>n}$ and $K_{*>n}$.

Let $m,n \in \mathbb{N}$, and let $U$ be a subset of $(\bigwedge_*^0)^m \oplus (\bigwedge_*^1)^n$.  A $\bigwedge_*$-superfunction $H$ on $U$ in $(m,n)$-variables is given by
\begin{eqnarray*}
H: U &\longrightarrow& \mbox{$\bigwedge_*$} \\ 
(z_1,z_2,\dots,z_m ,\theta_1,\theta_2,\dots,\theta_n) &\mapsto& H(z_1,z_2,\dots,z_m ,\theta_1,\theta_2,\dots,\theta_n)
\end{eqnarray*}
where $z_k$, for $k = 1,\dots,m$, are even variables in $\bigwedge_*^0$ and $\theta_k$, for $k = 1,\dots,n$, are odd variables in $\bigwedge_*^1$.  Let $f((z_1)_B,(z_2)_B,\dots,(z_m)_B)$ be a complex analytic function in $(z_k)_B$, for $k = 1,\dots,m$.  For $z_k \in \bigwedge_*^0$, and $k = 1,..,m$, define
\begin{multline}\label{more-than-one-variable}
f(z_1,z_2,\dots,z_m) = \! \sum_{l_1,\dots,l_m \in \mathbb{N}} \! \! \frac{(z_1)_S^{l_1} (z_2)_S^{l_2} \cdots
(z_m)_S^{l_m}}{l_1 ! l_2 ! \cdots l_m !} \biggl(\frac{\partial \; \;}{\partial (z_1)_B}\biggr)^{l_1} \cdot \\
\biggl(\frac{\partial \; \;}{\partial (z_2)_B}\biggr)^{l_2} \! \! \cdots \biggl(\frac{\partial \; \;}{\partial (z_m)_B}\biggr)^{l_m} \! \! \cdot f((z_1)_B,(z_2)_B,\dots,(z_m)_B) .
\end{multline}

\begin{defn}\label{superanal-definition}
Let $m, n \in \mathbb{N}$.  Let $U \subset (\bigwedge_{*>n-1}^0)^m \oplus (\bigwedge_{*>n-1}^1)^n$, and let $H$ be a $\bigwedge_{*>n-1}$-superfunction in $(m,n)$-variables defined on $U$. Then $H$ is said to be {\em superanalytic} if $H$ is of the form
\begin{equation}
H(z_1,z_2,\dots,z_m ,\theta_1,\theta_2,\dots,\theta_n) = \sum_{ (l) \in K_n} \theta_{l_1} \cdots \theta_{l_{j}} f_{(l)}(z_1,z_2,\dots,z_m) , 
\end{equation} 
where each $f_{(l)}$ is of the form 
\[f_{(l)}(z_1,z_2,\dots,z_m) = \sum_{(k) \in K_{ * - n}} f_{(l),(k)}(z_1,z_2,\dots,z_m) \zeta_{k_1}\zeta_{k_2} \cdots \zeta_{k_{s}}, \]  
and each $f_{(l),(k)}((z_1)_B,(z_2)_B,\dots,(z_m)_B)$ is analytic  in $(z_i)_B$, for $i = 1,\dots,m$, and $((z_1)_B,(z_2)_B,\dots,(z_m)_B) \in  U_B \subset \mathbb{C}^m$.   We call $H$ an {\em even} superanalytic $(m,n)$-function if whenever $(l) \in I_n$ then $f_{(l),(k)} \equiv 0$ for $(k) \notin I_{*-n}$ and whenever $(l) \in J_n$ then $f_{(l),(k)} \equiv 0$ for $(k) \notin J_{*-n}$.  We call  $H$ an {\em odd} superanalytic $(m,n)$-function if whenever $(l) \in I_n$ then $f_{(l),(k)} \equiv 0$ for $(k) \notin J_{*-n}$ and whenever $(l) \in J_n$ then $f_{(l),(k)} \equiv 0$ for $(k) \notin I_{*-n}$. 
\end{defn}

We require the even and odd variables to be in $\bigwedge_{* > n-1}$, and we restrict the coefficients of the $f_{(l),(k)}$'s to be in $\bigwedge_{* - n} \subseteq \bigwedge_{*>n-1}$ in order for the partial derivatives with respect to each of the $n$ odd variables to be well defined and for multiple partials to be well defined (cf. \cite{D}, \cite{B-memoir}).  

We define the (left) partial derivatives $\frac{\partial}{\partial z_j}$ and $\frac{\partial}{\partial \theta_j}$ acting on some superanalytic superfunction $H(z_1,\dots,z_m ,\theta_1,\dots,\theta_n)$ defined on $U \subset (\bigwedge_{*>n-1}^0)^m \oplus (\bigwedge_{*>n-1}^1)^n$ by
\begin{multline*}
\Delta z_j \left(\frac{\partial}{\partial z_j} H(z_1,\dots,z_m ,\theta_1,\dots,\theta_n) \right) + O((\Delta z_j)^2) \\
= H(z_1,\dots,z_{j-1}, z_j + \Delta z_j,z_{j+1},\dots,z_m, \theta_1,\dots,\theta_n) - H(z_1,\dots,z_m ,\theta_1,\dots,\theta_n) 
\end{multline*}
for all $\Delta z_j \in \bigwedge_{*>0}^0$ such that $(z_1,\dots,z_{j-1},z_j + \Delta z_j,z_{j+1},\dots,z_m) \in U^0$, for $j=1,\dots,m$, and
\begin{multline*}
\Delta \theta_j \left(\frac{\partial}{\partial \theta_j} H(z_1,\dots,z_m ,\theta_1,\dots,\theta_n) \right) \\
= H(z_1,\dots,z_m ,\theta_1,\dots,\theta_{j-1},\theta_j +\Delta \theta_j, \theta_{j+1},\dots,\theta_n) - H(z_1,\dots,z_m ,\theta_1,\dots,\theta_n) 
\end{multline*}
for all $\Delta \theta_j \in \bigwedge_{*>0}^1$ such that $(\theta_1,\dots, \theta_{j-1},\theta_j + \Delta \theta_j,\theta_{j+1},\dots,\theta_n) \in U^1$.  Note that $\frac{\partial}{\partial z_i}$, for $j = 1,\dots,m$, and  $\frac{\partial}{\partial \theta_j}$, for $j = 1,\dots,n$, are endomorphisms of the superalgebra of superanalytic $(m,n)$-superfunctions, and in fact, are even and odd superderivations, respectively. 

Consider the projection
\begin{eqnarray}
\pi^{(m,n)}_B : (\mbox{$\bigwedge_{* > n-1}^0$})^m \oplus (\mbox{$\bigwedge_{* > n-1}^1$})^n & \longrightarrow & \mathbb{C}^m  \label{projection-onto-body}\\ 
(z_1,\dots,z_m, \theta_1,\dots,\theta_n) & \mapsto & ((z_1)_B,(z_2)_B,\dots,(z_m)_B) . \nonumber
\end{eqnarray}
We define the {\it DeWitt topology on $(\bigwedge_{* > n-1}^0)^m \oplus (\bigwedge_{* > n-1}^1)^n$} by letting 
\[U \subseteq ((\mbox{$\bigwedge_{*> n-1}^0$})^m \oplus (\mbox{$\bigwedge_{* > n-1}^1$})^n)\] 
be an open set in the DeWitt topology if and only if $U = (\pi^{(m,n)}_B)^{-1} (V)$ for some open set $V \subseteq \mathbb{C}^m$.  Note that the natural domain of a superanalytic $\bigwedge_{* > n-1}$-superfunction in $(m,n)$-variables is an open set in the DeWitt topology.

A ``superconformal" field theory based on ``superfields'' which are superanalytic superfunctions in $(1,n)$-variables satisfying certain symmetry conditions would be referred to as an ``$N = n$ superconformal field theory''.

Let $(\bigwedge_*)^\times$ denote the set of invertible elements in $\bigwedge_*$.  Then
\[(\mbox{$\bigwedge_*$})^\times = \{a \in \mbox{$\bigwedge_*$} \; | \; a_B \neq 0 \} \]
since 
\[\frac{1}{a} = \frac{1}{a_B + a_S} = \sum_{n \in \mathbb{N}} \frac{(-1)^n a_S^n}{a_B^{n + 1}} \]
is well defined if and only if $a_B \neq 0$.  In light of this fact, note that the DeWitt topology is non-Hausdorff.  For example, two points $a, b \in \bigwedge_*$ can be separated by disjoint open sets in 
the DeWitt topology if and only if $a_B \neq b_B$, i.e., if and only if their difference is an invertible element of $\bigwedge_*$.  In other words, the DeWitt topology fails to be Hausdorff exactly to the extent 
that the nonzero elements of $\bigwedge_*$ fail in general to be invertible.

\begin{rema}\label{H_L-remark}
{\em Recall that $\bigwedge_L \subset \bigwedge_{L+1}$ for $L \in \mathbb{N}$, and note that {}from (\ref{more-than-one-variable}), any superanalytic $\bigwedge_L$-superfunction, $H_L$, in $(m,n)$-variables for $L \geq n$ can naturally be extended to a superanalytic $\bigwedge_{L'}$-superfunction in $(m,n)$-variables for $L'>L$ and hence to a superanalytic $\bigwedge_\infty$-superfunction.  Conversely, if $H_{L'}$ is a superanalytic $\bigwedge_{L'}$-superfunction (or $\bigwedge_\infty$-superfunction) in $(m,n)$-variables for $L' >n$, then we can restrict $H_{L'}$ to a superanalytic $\bigwedge_L$-superfunction for $L'> L\geq n$ by restricting $(z_1,\dots,z_m,\theta_1,\dots,\theta_n) \in 
(\bigwedge_L^0)^m \oplus (\bigwedge_L^1)^n$ and setting $f_{(l),(k)} \equiv 0$ if $(k) \notin K_{L-n}$. }
\end{rema}

\section{Superconformal $(1,2)$-superfunctions and power series}\label{superconformal-section}

In this section we give the definition of $N=2$ superconformal superfunction and make note of the power series expansions of such functions vanishing at zero or infinity (cf. \cite{Ki}, \cite{Nogueira}). 

\begin{rema}{\em 
There are typically two different coordinate systems one can choose to work with in $N=2$ superconformal field theory (cf. \cite{C}, \cite{Nogueira}).  Throughout most of this work, we have chosen to work with what we call the ``homogeneous" coordinate system.  In Section \ref{nonhomo-section} we discuss the ``nonhomogeneous" coordinate system and show how to convert our results to that setting. 
}
\end{rema}

Let $z$ be an even variable in $\bigwedge_{*>1}^0$, and let $\theta^+$ and $\theta^-$ be odd variables in $\bigwedge_{*>1}^1$.   Define
\begin{equation}
D^+ = \frac{\partial}{\partial \theta^+} + \theta^- \frac{\partial}{\partial z}
\qquad \qquad 
D^- = \frac{\partial}{\partial \theta^-} + \theta^+ \frac{\partial}{\partial z}.
\end{equation}
Then $D^\pm$ are odd superderivations on $\bigwedge_{*>1}$-superfunctions in $(1,2)$-variables which are superanalytic in some DeWitt open subset $U \subseteq \bigwedge_{*>1}^0 \oplus (\bigwedge_{*>1}^1)^2$.   Note that
\begin{eqnarray}
[D^\pm, D^\pm] & = & 2(D^\pm)^2 \ = \ 0\\ 
\left[D^+,D^-\right] &=& D^+D^- + D^-D^+ \ = \ 2 \frac{\partial}{\partial z} .
\end{eqnarray}

Let
\begin{eqnarray}
\qquad H : U \subseteq \mbox{$\bigwedge_{*>1}^0$} \oplus (\mbox{$\bigwedge_{*>1}^1$})^2
\! \! &\rightarrow& \! \! \mbox{$\bigwedge_{*>1}^0$} \oplus (\mbox{$\bigwedge_{*>1}^1$})^2 \label{H-superanal}\\
(z,\theta^+,\theta^-) \! \! &\mapsto& \! \! (H^0(z, \theta^+, \theta^-), H^+(z, \theta^+, \theta^-), H^-(z, \theta^+, \theta^-) ) \nonumber\\
& & = (\tilde{z},\tilde{\theta}^+,\tilde{\theta}^-) \nonumber
\end{eqnarray} 
be superanalytic, i.e., $H^0(z, \theta^+, \theta^-) = \tilde{z}$ is an even superanalytic $(1,2)$-function in the sense of Definition \ref{superanal-definition}, and $H^\pm(z, \theta^+, \theta^-) = \tilde{\theta}^\pm$ are odd superanalytic $(1,2)$-functions.   Then $D^+$ and $D^-$ transform under $H(z,\theta^+,\theta^-)$ by
\begin{eqnarray}
D^+ = (D^+\tilde{\theta}^+)\tilde{D}^+ + (D^+\tilde{\theta}^-) \frac{\partial}{\partial \tilde{\theta}^-} + (D^+\tilde{z} - \tilde{\theta}^- D^+\tilde{\theta}^+) \frac{\partial}{\partial \tilde{z}} \label{transform-Dplus}\\
D^- = (D^-\tilde{\theta}^-)\tilde{D}^- + (D^-\tilde{\theta}^+) \frac{\partial}{\partial \tilde{\theta}^+} + (D^-\tilde{z} - \tilde{\theta}^+ D^-\tilde{\theta}^-) \frac{\partial}{\partial \tilde{z}} \label{transform-Dminus}
\end{eqnarray}
where $\tilde{D}^\pm = \frac{\partial}{\partial \tilde{\theta}^\pm} + \tilde{\theta}^\mp \frac{\partial}{\partial \tilde{z}}$ with $\frac{\partial}{\partial \tilde{z}}$ and $\frac{\partial}{\partial \tilde{\theta}^\pm}$ defined by the usual chain rule, i.e.,
\begin{eqnarray*} 
\dz &=& \frac{\partial \tilde{z}}{\partial z} \frac{\partial}{\partial \tilde{z}} + \frac{\partial \tilde{\theta}^+}{\partial z}  \frac{\partial}{\partial \tilde{\theta}^+} + \frac{\partial \tilde{\theta}^-}{\partial z} \frac{\partial}{\partial \tilde{\theta}^-}\\ 
\frac{\partial}{\partial \theta^+} &=&  \frac{\partial \tilde{z}}{\partial \theta^+} \frac{\partial}{\partial \tilde{z}} + \frac{\partial \tilde{\theta}^+}{\partial \theta^+} \frac{\partial}{\partial \tilde{\theta}^+} +  \frac{\partial \tilde{\theta}^-}{\partial \theta^+} \frac{\partial}{\partial \tilde{\theta}^-}\\
\frac{\partial}{\partial \theta^-} &=&  \frac{\partial \tilde{z}}{\partial \theta^-} \frac{\partial}{\partial \tilde{z}} + \frac{\partial \tilde{\theta}^+}{\partial \theta^-} \frac{\partial}{\partial \tilde{\theta}^+} +  \frac{\partial \tilde{\theta}^-}{\partial \theta^-} \frac{\partial}{\partial \tilde{\theta}^-}.
\end{eqnarray*}

Recall that a complex function $f$ defined on an open set $U_B$ in $\mathbb{C}$, of one complex variable $z_B$, is conformal in $U_B$ if and only if $\frac{d \;}{d z_B} f(z_B)$ exists for $z_B \in U_B$ and is not identically zero in $U_B$, i.e., if and only if $f(z_B) = \tilde{z}_B$ transforms $\frac{d \;}{d z_B}$ by $\frac{d \;}{d z_B} = f'(z_B) \frac{d \;}{d \tilde{z}_B}$ for $f'$ not identically zero.  Such a transformation of $\frac{d \;}{d z_B}$ is said to be {\it homogeneous of degree one}, i.e., $f$ transforms $\frac{d \;}{d z_B}$ by a non-zero analytic function times $\frac{d \;}{d \tilde{z}_B}$ to the first power with no higher order terms in $\frac{d\;}{d\tilde{z}_B}$.   Analogously we define an {\it $N=2$  superconformal function} $H$ on a DeWitt open subset $U$ of $\bigwedge_{*>1}^0 \oplus (\bigwedge_{*>1}^1)^2$ to $\bigwedge_{*>1}^0 \oplus (\bigwedge_{*>1}^1)^2$ to be a superanalytic function under which $D^+$ and $D^-$ transform homogeneously of degree one.   That is, $H$ transforms $D^\pm$ by non-zero superanalytic functions times $\tilde{D}^\pm$, respectively.   Since such a superanalytic function $H(z,\theta^+,\theta^-) = (\tilde{z}, \tilde{\theta}^+,\tilde{\theta}^-)$  transforms $D^+$ and $D^-$ according to  (\ref{transform-Dplus}) and (\ref{transform-Dminus}), respectively,  $H$ is superconformal if and only if, in addition to being superanalytic, $H$ satisfies
\begin{eqnarray}
D^+ \tilde{\theta}^- &=& 0, \label{basic-superconformal-condition1} \\
D^- \tilde{\theta}^+ &=& 0, \label{basic-superconformal-condition2}\\
D^+\tilde{z} - \tilde{\theta}^- D^+\tilde{\theta}^+ &=& 0, \label{basic-superconformal-condition3}\\
D^-\tilde{z} - \tilde{\theta}^+ D^-\tilde{\theta}^- &=& 0, \label{basic-superconformal-condition4}
\end{eqnarray}
for $D^\pm \tilde{\theta}^\pm$ not identically zero, thus transforming $D^\pm$ by $D^\pm = (D^\pm \tilde{\theta}^\pm)\tilde{D}^\pm$.

We can write $H(z,\theta^+,\theta^-) = (\tilde{z}, \tilde{\theta}^+, \tilde{\theta}^-)$ as
\begin{eqnarray*}
\tilde{z} &=& f(z) + \theta^+ \xi^+(z) + \theta^- \xi^- (z) + \theta^+ \theta^-f^{+,-}(z)\\
\tilde{\theta}^+ &=& \psi^+(z) + \theta^+ g^+(z) + \theta^- h^-(z) + \theta^+ \theta^-\phi^+(z)\\
\tilde{\theta}^- &=& \psi^-(z) + \theta^+ h^+(z) + \theta^- g^-(z) + \theta^+ \theta^- \phi^-(z)
\end{eqnarray*}
for $f,f^{+,-},g^\pm,h^\pm$ even and $\xi^\pm,\psi^\pm,\phi^\pm$ odd $(1,0)$-superfunctions in $z$.  Then the conditions (\ref{basic-superconformal-condition1}) -- (\ref{basic-superconformal-condition4}) are equivalent to the conditions
\begin{eqnarray}
\qquad \qquad \tilde{z} &=& f(z) + \theta^+ g^+(z) \psi^-(z) + \theta^- g^-(z) \psi^+(z) + \theta^+ \theta^- (\psi^+(z)\psi^-(z))' \label{superconformal-condition1} \\
\tilde{\theta}^+ &=& \psi^+(z) + \theta^+ g^+(z) + \theta^+ \theta^-(\psi^+)'(z) \label{superconformal-condition2} \\
\tilde{\theta}^- &=& \psi^-(z) + \theta^- g^-(z) - \theta^+ \theta^- (\psi^-)'(z) \label{superconformal-condition3}
\end{eqnarray}
and
\begin{equation}\label{superconformal-condition4}
f'(z) \; = \; (\psi^+)'(z)\psi^-(z) - \psi^+(z) (\psi^-)'(z) + g^+(z) g^-(z) ,
\end{equation}
and we also require that $D^+ \tilde{\theta}^+$ and $D^- \tilde{\theta}^-$ not be identically zero, which is equivalent to
\begin{eqnarray}\label{non-zero-condition}
g^+(z) + (\psi^+) ' (z)  \equiv \! \! \! \! \negthickspace /  \ \  0 \qquad \mathrm{and} \qquad 
g^-(z) + (\psi^-) ' (z)   \equiv \! \! \! \! \negthickspace / \ \  0 . 
\end{eqnarray}
Thus  an $N=2$ superconformal function $H$ is uniquely determined by $f(z)$, $\psi^\pm(z)$,  $g^\pm(z)$ satisfying the conditions (\ref{superconformal-condition4}) and (\ref{non-zero-condition}).

Note that the space of $N=2$ superconformal functions on $\bigwedge_{*>1}^0 \oplus (\bigwedge_{*>1}^1)^2$ is closed under composition when defined.  However, the sum of two $N=2$ superconformal functions is not in general superconformal.

In Section 5, we will study ``$N=2$ super-Riemann spheres with punctures and local superconformal coordinates vanishing at the punctures''.   These punctures can be thought of as being at $0 \in \bigwedge_{*>1}^0 \oplus (\bigwedge_{*>1}^1)^2$, a non-zero point in $\bigwedge_{*>1}^0 \oplus (\bigwedge_{*>1}^1)^2$, or at a distinguished point on the $N=2$ super-Riemann sphere we denote by ``$\infty$''.  As will be shown in Section 5, we can always shift a non-zero point in $\bigwedge_{*>1}^0 \oplus (\bigwedge_{*>1}^1)^2$ (or on the $N=2$ super-Riemann sphere) to zero via a global $N=2$ superconformal transformation.  Thus all local superconformal coordinates vanishing at the punctures can be expressed as power series vanishing at zero or vanishing as $(z,\theta^+,\theta^-) = (z_B + z_S, \theta^+,\theta^-) \longrightarrow (\infty  + 0, 0,0) = \infty$.  

If the puncture is at zero, we are interested in invertible $N=2$ superconformal functions $H(z,\theta^+,\theta^-) = (\tilde{z}, \tilde{\theta}^+, \tilde{\theta}^-)$ defined in a neighborhood  of zero vanishing at zero.  Such an $H$ satisfies  (\ref{superconformal-condition1}) -- (\ref{non-zero-condition}) 
where $f(z), g^\pm(z)$ are even superanalytic functions and $\psi^\pm(z)$ are odd superanalytic functions, satisfying $f(0) = \psi^\pm(0) = 0$ and $f'(0), g^\pm(0) \in (\bigwedge_{ * - 2}^0)^\times$.  Thus an $N=2$ superconformal power series invertible in a neighborhood of zero and vanishing at zero is uniquely determined by
\begin{eqnarray}
g^\pm(z) &= &\sum_{j \in \mathbb{N}} a^\pm_j z^j , \qquad \quad \mbox{for $a^\pm_j \in\bigwedge_{ * - 2}^0$ and $a^\pm_0 \in (\bigwedge_{ * - 2}^0)^\times$} \label{expansion-of-g}\\
\psi^\pm(z) &= &\sum_{j \in \mathbb{N}} m^\pm_{j+\frac{1}{2}} z^{j + 1} , \qquad \mbox{for $m^\pm_{j+\frac{1}{2}} \in \bigwedge_{ * - 2}^1$} , \label{expansion-of-psi}
\end{eqnarray}
and has the form (\ref{superconformal-condition1}) -- (\ref{superconformal-condition3}), where $f$ is uniquely determined by (\ref{superconformal-condition4}) since $f(0) = 0$.

\begin{rema}\label{determining-superconformal-remark}
{\em {}From the previous statement, we see that if $H(z, \theta^+, \theta^-) = (\tilde{z}, \tilde{\theta}^+, \tilde{\theta}^-)$ is $N=2$ superconformal and vanishing at zero, then it is completely determined by the components $\tilde{\theta}^\pm$.   In fact we can say more:  $H$ vanishing at zero is completely determined by $\theta^\mp \tilde{\theta}^\pm$, since then we are still able to pick out the $\psi^\pm$ and $g^\pm$ components.  Note that such an $H$ is not determined uniquely by $\tilde{z}$ alone.  This is in contrast to the $N=1$ case (cf. \cite{B-memoir}).
}
\end{rema}

Explicitly, we have that if $H(z,\theta^+,\theta^-) = (\tilde{z}, \tilde{\theta}^+, \tilde{\theta}^-)$ is superconformal and invertible in a neighborhood of zero and vanishing at zero, then 
\begin{eqnarray}
\tilde{z} &= & \! \! \sum_{j,k \in \mathbb{N}} \Bigl( \frac{1}{j+k+1} a^+_j a^-_k z^{j+k+1} +  \frac{j-k}{j+k+2} m^+_{j+\frac{1}{2}}  m^-_{k+\frac{1}{2}} z^{j+k + 2}   \Bigr. \label{powerseries1} \\
& & \quad + \;  \theta^+ a^+_j m^-_{k+\frac{1}{2}}  z^{j +k + 1} + \theta^-  a^-_j m^+_{k+\frac{1}{2}} z^{j +k+ 1} \nonumber \\
& & \Bigl. \quad +\;  \theta^+ \theta^-  (j+k +2) m^+_{j+\frac{1}{2}} m^-_{k+\frac{1}{2}} z^{j +k+ 1} \Bigr) \nonumber \\
\tilde{\theta}^\pm &=& \sum_{j \in \mathbb{N}} \Bigl( m^\pm_{j+\frac{1}{2}} z^{j + 1}  + \theta^\pm a^\pm_j z^j \pm \theta^+ \theta^-(j+1)m^\pm_{j+\frac{1}{2}} z^j  \Bigr) \label{powerseries2} 
\end{eqnarray}
for $a^\pm_0 \in (\bigwedge_{ * - 2}^0)^\times$, $a^\pm_j \in\bigwedge_{ * - 2}^0$, for $j \in \Z$, and 
$m^\pm_{j+\frac{1}{2}} \in \bigwedge_{ * - 2}^1$, for $j \in \mathbb{N}$.

Similarly, we would like to express an $N=2$ superconformal function vanishing as $(z,\theta^+,\theta^-) = (z_B + z_S, \theta^+,\theta^-) \longrightarrow (\infty  + 0, 0,0) = \infty$ as a power series in $z$ and $\theta^\pm$.  The superfunction 
\begin{eqnarray*}
I: \mbox{$(\bigwedge_{*>1}^0)^\times$} \oplus \mbox{$(\bigwedge_{*>1}^1)^2$} & \longrightarrow & \mbox{$(\bigwedge_{*>1}^0)^\times$} \oplus \mbox{$(\bigwedge_{*>1}^1)^2$}\\ 
(z,\theta^+,\theta^-) & \mapsto & \Bigl(\frac{1}{z}, \frac{i \theta^+}{z}, \frac{i\theta^-}{z} \Bigr)
\end{eqnarray*}
is $N=2$ superconformal, well defined and vanishing as $(z,\theta^+, \theta^-) \rightarrow \infty$.  In fact, $H$ is $N=2$ superconformal, well defined and invertible in a neighborhood of $(\infty,0,0)$ and vanishing at $(z,\theta^+, \theta^-) = \infty$ if  and only if $H(1/z, i\theta^+/z, i\theta^-/z)= (\tilde{z}, \tilde{\theta}^+, \tilde{\theta}^-)$ is of the form  (\ref{powerseries1}) -- (\ref{powerseries2}).

Explicitly, we have that if $H(z,\theta^+,\theta^-) = (\tilde{z}, \tilde{\theta}^+, \tilde{\theta}^-)$ is $N=2$ superconformal and invertible in a neighborhood of infinity and vanishing at infinity, then 
\begin{eqnarray}
\tilde{z} &= & \! \! \sum_{j,k \in \mathbb{N}} \Bigl( \frac{1}{j+k+1} a^+_j a^-_k z^{-j-k-1} +  \frac{j-k}{j+k+2} m^+_{j+\frac{1}{2}}  m^-_{k+\frac{1}{2}} z^{-j-k - 2}   \Bigr. \label{powerseries1-infinity} \\
& & \quad + \;  i\theta^+ a^+_j m^-_{k+\frac{1}{2}}  z^{-j -k -2} + i \theta^-  a^-_j m^+_{k+\frac{1}{2}} z^{-j -k- 2} \nonumber \\
& & \Bigl. \quad - \;  \theta^+ \theta^-  (j+k +2) m^+_{j+\frac{1}{2}} m^-_{k+\frac{1}{2}} z^{-j -k- 3} \Bigr) \nonumber \\
\tilde{\theta}^\pm &=& \sum_{j \in \mathbb{N}} \Bigl( m^\pm_{j+\frac{1}{2}} z^{-j - 1}  + i\theta^\pm a^\pm_j z^{-j-1} \mp \theta^+  \theta^-(j+1) m^\pm_{j+\frac{1}{2}} z^{-j-2}  \Bigr) \label{powerseries2-infinity} 
\end{eqnarray}
for $a^\pm_0 \in (\bigwedge_{ * - 2}^0)^\times$, $a^\pm_j \in\bigwedge_{ * - 2}^0$, for $j \in \Z$, and 
$m^\pm_{j+\frac{1}{2}} \in \bigwedge_{ * - 2}^1$, for $j \in \mathbb{N}$.

\section{Complex supermanifolds, $N=2$ super-Riemann surfaces and $N=2$ superspheres with tubes}  

In this section we give the definitions of supermanifold, $N=2$ super-Riemann surface and genus-zero $N=2$ super-Riemann surface which we call an ``$N=2$ supersphere".  We then study $N=2$ superspheres with ordered and oriented tubes and show that these are superconformally equivalent to $N=2$ superspheres with ordered and oriented punctures and local $N=2$ superconformal coordinates vanishing at the punctures.

A {\em DeWitt $(m,n)$-dimensional supermanifold over $\bigwedge_*$} is a topological space $X$ with a countable basis which is locally homeomorphic to an open subset of $(\bigwedge_*^0)^m \oplus (\bigwedge_*^1)^n$ in the DeWitt topology.  A {\em DeWitt $(m,n)$-chart on $X$ over $\bigwedge_*$} is 
a pair $(U, \Omega)$ such that $U$ is an open subset of $X$ and $\Omega$ is a homeomorphism of $U$ onto an open subset of $(\bigwedge_*^0)^m \oplus (\bigwedge_*^1)^n$ in the DeWitt topology.  A {\em superanalytic atlas of DeWitt $(m,n)$-charts on $X$ over $\bigwedge_{* > n-1}$} is a family of charts $\{(U_{\alpha}, \Omega_{\alpha})\}_{\alpha \in A}$ satisfying 
 
(i) Each $U_{\alpha}$ is open in $X$, and $\bigcup_{\alpha \in A} U_{\alpha} = X$. 

(ii) Each $\Omega_{\alpha}$ is a homeomorphism {}from $U_{\alpha}$ to a (DeWitt) open set in $(\bigwedge_{* > n-1}^0)^m \oplus (\bigwedge_{* > n-1}^1)^n$, such that $\Omega_{\alpha} \circ \Omega_{\beta}^{-1}: \Omega_{\beta}(U_\alpha \cap U_\beta) \longrightarrow \Omega_{\alpha}(U_\alpha \cap U_\beta)$ is superanalytic for all non-empty $U_{\alpha} \cap U_{\beta}$, i.e., $\Omega_{\alpha} \circ \Omega_{\beta}^{-1} = (\tilde{z}_1,\dots, \tilde{z}_m, \tilde{\theta}_1,\dots,\tilde{\theta}_n)$ where $\tilde{z}_i$ 
is an even superanalytic $\bigwedge_{* > n-1}$-superfunction in $(m,n)$-variables for $i = 1,\dots,m$, and $\tilde{\theta}_j$ is an odd superanalytic $\bigwedge_{* >n-1}$-superfunction in $(m,n)$-variables 
for $j = 1,\dots,n$.  

Such an atlas is called {\em maximal} if, given any chart $(U, \Omega)$ such that
\[\Omega \circ \Omega_{\beta}^{-1} : \Omega_{\beta} (U \cap U_\beta) \longrightarrow \Omega (U \cap U_\beta)\] 
is a superanalytic homeomorphism for all $\beta$, then $(U, \Omega) \in \{(U_{\alpha}, \Omega_{\alpha})\}_{\alpha \in A}$.
 
A {\em DeWitt $(m,n)$-superanalytic supermanifold over $\bigwedge_{* > n-1}$} is a DeWitt $(m,n)$-dimensional supermanifold $M$ together with a maximal superanalytic atlas of DeWitt $(m,n)$-charts over $\bigwedge_{* > n-1}$.  

Given a DeWitt $(m,n)$-superanalytic supermanifold $M$ over $\bigwedge_{* > n-1}$, define an equivalence relation $\sim$ on M by letting $p \sim q$ if and only if there exists $\alpha \in A$ such that $p,q \in U_\alpha$ and $\pi_B^{(m,n)} (\Omega_\alpha (p)) = \pi_B^{(m,n)} (\Omega_\alpha (q))$
where $\pi_B^{(m,n)}$ is the projection given by (\ref{projection-onto-body}).  Let $p_B$ denote the equivalence class of $p$ under this equivalence relation.  Define the {\it body} $M_B$ of $M$ to be the 
$m$-dimensional complex manifold with analytic structure given by the coordinate charts $\{((U_\alpha)_B, (\Omega_\alpha)_B) \}_{\alpha \in A}$ where $(U_\alpha)_B = \{ p_B \; | \; p \in U_\alpha \}$, and $(\Omega_\alpha)_B : (U_\alpha)_B \longrightarrow \mathbb{C}^m$ is given by $(\Omega_\alpha)_B (p_B) = \pi_B^{(m,n)} \circ \Omega_\alpha (p)$.     

We have that $M$ is a complex fiber bundle over the complex manifold $M_{B}$.  The fiber is $(\bigwedge_{* > n-1}^0)_S^m \oplus (\bigwedge_{* > n-1}^1)^n$, a possibly infinite-dimensional vector space over $\mathbb{C}$.  This bundle is not in general a vector bundle since the transition functions are not in general linear.  

For any DeWitt $(1,n)$-superanalytic supermanifold $M$, its body $M_{B}$ is a Riemann surface. An {\em $N=2$ super-Riemann surface over $\bigwedge_{*>1}$} is a DeWitt $(1,2)$-superanalytic supermanifold over $\bigwedge_{*>1}$ with coordinate atlas $\{(U_{\alpha}, \Omega_{\alpha})\}_{\alpha \in A}$ such that the coordinate transition functions $\Omega_{\alpha} \circ \Omega_{\beta}^{-1}$ in addition to being superanalytic are also $N=2$ superconformal for all non-empty $U_{\alpha} \cap U_{\beta}$.  Since the condition that the coordinate transition functions be superconformal instead of merely superanalytic is such a strong condition (unlike in the nonsuper case), we again stress the distinction between an $N=2$ supermanifold which has {\it superanalytic} transition functions versus an $N=2$ super-Riemann surface which has $N=2$ {\it superconformal} transition functions.  It would be perhaps more appropriate to refer to the later as an ``$N=2$ superconformal super-Riemann surface'' in order to avoid confusion.  In fact, in the literature one will find the term ``super-Riemann surface" or ``Riemannian supermanifold" used for both merely superanalytic structures (cf. \cite{D}) and for superconformal structures (cf. \cite{Fd}, \cite{CR}).  However, we will follow the terminology of \cite{Fd} for the $N=1$ super case and refer to an $N=2$ superconformal super-Riemann surface simply as an $N=2$ super-Riemann surface.

By {\em $N=2$ supersphere} we will mean a (superconformal) $N=2$ super-Riemann surface over $\bigwedge_{*>1}$ such that its body is a genus-zero one-dimensional connected compact complex manifold.  

An {\em $N=2$ supersphere with $1+n$ tubes} for $n \in \mathbb{N}$, is a supersphere $S$ with one negatively oriented point $p_0$ and $n$ positively oriented points $p_1,...,p_n$ (we call them {\em
punctures}) on $S$ which all have distinct bodies (i.e., $(p_j)_B \neq (p_k)_B$ if $j \neq k$, or equivalently $p_j$ is not equivalent to $p_k$ for $j \neq k$ under the equivalence relation $\sim$) and with local $N=2$ superconformal coordinates $(U_0,\Omega_0),\dots, (U_n, \Omega_n)$ vanishing at the  punctures $p_0, \dots ,p_n$, respectively.  We denote this structure by
\[(S;p_0,...,p_n;(U_0, \Omega_0),...,(U_n, \Omega_n)) . \]
We will always order the punctures so that the negatively oriented
puncture is $p_0$.  

The reason we call a puncture with local $N=2$ superconformal coordinate vanishing at the puncture a ``tube" is that such a structure is indeed superconformally equivalent to a half-infinite $N=2$ superconformal tube representing an incoming (resp., outgoing) $N=2$ ``superparticle" or ``superstring" propagating through space-time if the puncture is positively (resp., negatively) oriented.  For $r \in \mathbb{R}_+$, denote by  
\[ B_{z_B}^r = \{ w_B \in \mathbb{C} \; | \; |w_B - z_B| < r \} \qquad (\mbox{resp.,} \; \; \bar{B}_{z_B}^r = \{ w_B \in \mathbb{C} \; | \; |w_B - z_B| \leq r \}) \] 
an open (resp., closed) ball in the complex plane about the point $z_B$ with radius $r$.  Denote a DeWitt open (resp., closed) ball in $\bigwedge_{*>1}^0 \oplus (\bigwedge_{*>1}^1)^2$ about $(z, \theta^+,\theta^-)$ of radius $r$ by
\[ \mathcal{B}_z^r = B_{z_B}^r \times \bigl((\mbox{$\bigwedge_{*>1}^0$})_S \oplus \mbox{$(\bigwedge_{*>1}^1)^2$}\bigr) \qquad (\mbox{resp.,} \; \; \bar{\mathcal{B}}_z^r = \bar{B}_{z_B}^r \times \bigl((\mbox{$\bigwedge_{*>1}^0$})_S \oplus \mbox{$(\bigwedge_{*>1}^1)^2$}\bigr) .\]
(Note that $\mathcal{B}_z^r$ depends only on $z_B$ and $r$.)  Let $p$ be a positively oriented puncture on an $N=2$ supersphere with a local coordinate neighborhood $U$ and superconformal local coordinate map $\Omega: U \rightarrow \bigwedge_{*>1}^0 \oplus (\bigwedge_{*>1}^1)^2$ vanishing at the puncture.  Then for some $r \in \mathbb{R}_+$, we can find a DeWitt open disc 
$\mathcal{B}^r_0$ such that $\Omega^{-1}(\mathcal{B}^r_0) \subset U$.  Define the equivalence relation $\tilde{\sim}$ on $\bigwedge_{*>1}^0 \oplus (\bigwedge_{*>1}^1)^2$ by $(z_1,\theta_1^+,\theta_1^-) \tilde{\sim} (z_2,\theta_2^+,\theta^-_2)$ if  and only if $(z_1)_B = (z_2)_B + 2\pi i k$ for some integer $k$, $(z_1)_S = (z_2)_S$ and $\theta_1^\pm = \theta_2^\pm$.  Then the set $\tau_r$ of all equivalence classes of elements of $(z,\theta^+,\theta^-) \in \bigwedge_{*>1}^0 \oplus (\bigwedge_{*>1}^1)^2$ satisfying $\mathrm{Re}(z_B) < \log r$ (where $\mathrm{Re}(z_B)$ is the real part of the complex number $z_B$) together with the metric induced {}from the DeWitt metric on $\bigwedge_{*>1}^0 \oplus (\bigwedge_{*>1}^1)^2$ is a half-infinite tube in the body and is topologically trivial in the soul.  Letting $H(z,\theta^+,\theta^-) = (\log z, \theta^+ \sqrt{1/z}, \theta^- \sqrt{1/z})$, the map $q \mapsto H(\Omega(q))$ {}from $\Omega^{-1}(\mathcal{B}^r_0)$ to $\tau_r$ is a well-defined invertible $N=2$ superconformal map.  A closed curve on the $N=2$ supersphere shrinking to $p$ corresponds to a closed loop or ``superstring" around this half-infinite super-cylinder tending towards minus infinity in the body coordinate.  We can perform a similar superconformal transformation for the negatively oriented puncture with local coordinate in order to recover the half-infinite outgoing tube.

In $N=2$ superconformal field theory, one generally wants to consider $N=2$ superspheres and higher genus super-Riemann surfaces with $m \in \mathbb{Z}_+$ negatively oriented (i.e., outgoing) tubes and $n \in \mathbb{N}$ positively oriented (i.e., incoming) tubes.  However, for the purposes of this work, we restrict to genus zero and $m=1$. 

Let $(S_1;p_0,...,p_n;(U_0, \Omega_0),...,(U_n, \Omega_n))$ and $(S_2;q_0,...,q_n;(V_0, \Xi_0),...,(V_n, \Xi_n))$ be two $N=2$ superspheres with $1+n$ tubes.  A map $F : S_1 \rightarrow S_2$ will be said to be $N=2$ superconformal if $\Xi_\beta \circ F \circ \Omega_\alpha^{-1}$ is $N=2$ superconformal for all charts $(U_\alpha, \Omega_\alpha)$ of $S_1$, for all charts $(V_\beta, \Xi_\beta)$ of $S_2$, and for all $(z, \theta^+, \theta^-) \in \Omega_\alpha (U_\alpha)$ such that $F \circ \Omega_\alpha^{-1} 
(z, \theta^+, \theta^-) \in V_\beta$.  If there is an $N=2$ superconformal isomorphism $F : S_1 \rightarrow S_2$ such that for each $j = 0,..., n$, we have $F(p_j) = q_j$ and 
\[ \left. \Omega_j \right|_{W_j} = \left. \Xi_j \circ F \right|_{W_j} , \]
for $W_j$ some DeWitt neighborhood of $p_j$, then we say that these two $N=2$ superspheres with $1 + n$ tubes are {\em superconformally equivalent} and $F$ is a {\em superconformal equivalence} {}from   
\[ (S_1;p_0,...,p_n;(U_0, \Omega_0),...,(U_n, \Omega_n)) \]
to
\[ (S_2;q_0,...,q_n;(V_0, \Xi_0),...,(V_n, \Xi_n)) .\]
Thus the superconformal equivalence class of an $N=2$ supersphere with tubes depends only on the $N=2$ supersphere, the punctures, and the germs of the local coordinate maps vanishing at the punctures.

\section{The moduli space of $N=2$ super-Riemann spheres with tubes}\label{moduli-section}

In this section we define the $N=2$ super-Riemann sphere and  the moduli space of $N=2$ super-Riemann spheres with tubes.  We introduce {\it canonical $N=2$ superspheres with tubes}, and show that any $N=2$ super-Riemann sphere with tubes is $N=2$ superconformally equivalent to a canonical $N=2$ supersphere with tubes.  In addition, we show that two different canonical $N=2$ superspheres with tubes are not $N=2$ superconformally equivalent.  This shows that there is a bijection between the set of canonical $N=2$ superspheres with tubes and the moduli space of $N=2$ super-Riemann spheres with tubes.

\begin{rema}\label{infinite-variables1}
{\em In subsequent work, we will want to consider functions on the moduli space of $N=2$ super-Riemann spheres with tubes which are superanalytic or supermeromorphic.  These superfunctions will in general involve an infinite number of odd variables --- not only the odd part of the coordinate for the finite number of punctures but also the possibly infinite amount of odd data involved in describing the local coordinates about the punctures.  (See Remarks \ref{remark-giving-power-series} and \ref{infinite-variables2} below.)  In this case, we need to work over $\bigwedge_\infty$ if we want all multiple partial derivatives with respect to these odd variables to be well defined.  Since it is no harder to work over an infinite Grassmann algebra, {}from now on, this is what we will do.  One may always later restrict to some $\bigwedge_L$ for $L \in \mathbb{N}$ when substituting for these variables in the functional part of the theory or restrict to the supermanifold substructure defined over $\bigwedge_L$ (see \cite{B-memoir}) for geometric aspects of the theory. }
\end{rema}

Let $S^2\hat{\mathbb{C}}$ be the genus zero $N=2$ super-Riemann surface over $\bigwedge_\infty$ (meaning of course over $\bigwedge_\infty^0 \oplus (\bigwedge_\infty^1)^2$) with superconformal structure given by the covering of local coordinate neighborhoods $\{ U_\sou, U_\nor \}$ and the local coordinate maps
\begin{eqnarray}
\sou : U_\sou  & \longrightarrow & \mbox{$\bigwedge_\infty^0$} \oplus (\mbox{$\bigwedge_\infty^1$})^2 \\
\nor: U_\nor & \longrightarrow & \mbox{$\bigwedge_\infty^0$}  \oplus (\mbox{$\bigwedge_\infty^1$})^2,
\end{eqnarray}
which are homeomorphisms of $U_\sou$ and $U_\nor$ onto $\bigwedge_\infty^0 \oplus (\bigwedge_\infty^1)^2$, respectively, such that 
\begin{eqnarray}
\sou \circ \nor^{-1} : \mbox{$(\bigwedge_\infty^0)^\times$} \oplus \mbox{$(\bigwedge_\infty^1)^2$} &\longrightarrow& \mbox{$(\bigwedge_\infty^0)^\times$} \oplus \mbox{$(\bigwedge_\infty^1)^2$} \label{origin-of-I}\\ 
(w, \rho^+,\rho^-) & \mapsto & \Bigl(\frac{1}{w},\frac{i \rho^+}{w}, \frac{i \rho^-}{w} \Bigr) = I(w,\rho^+,\rho^-) . \nonumber
\end{eqnarray}
Thus the body of $S^2\hat{\mathbb{C}}$ is the Riemann sphere, $(S^2\hat{\mathbb{C}})_B = \hat{\mathbb{C}} = \mathbb{C} \cup \{\infty\}$, with coordinates $w_B$ near 0 and $1/w_B$ near $\infty$.   We will call $S^2\hat{\mathbb{C}}$ the {\em $N=2$ super-Riemann sphere} or just the {\em super-Riemann sphere} and will refer to $\nor^{-1} (0)$ as {\em the point at $(\infty, 0,0)$} or just {\em the point at infinity} and to $\sou^{-1}(0)$ as {\em the point at $(0, 0,0)$} or just {\em the point at zero}.  

An {\em $N=2$ super-Riemann sphere with $1+n$ tubes} for $n \in \mathbb{N}$, is an $N=2$ supersphere with $1+n$ tubes such that the underlying $N=2$ supersphere is the $N=2$ super-Riemann sphere.    The collection of all $N=2$ superconformal equivalence classes of $N=2$ super-Riemann spheres over $\bigwedge_\infty$ with $1+n$ tubes is called the {\it moduli space of $N=2$ super-Riemann spheres over $\bigwedge_\infty$ with $1+n$ tubes}.  The collection of all $N=2$ superconformal equivalence classes of  $N=2$ super-Riemann spheres over $\bigwedge_\infty$ with tubes is called the {\it moduli space of $N=2$ super-Riemann spheres over $\bigwedge_\infty$ with tubes}.

\begin{rema} 
{\em The results we obtain in this paper on the moduli space of $N=2$ super-Riemann spheres can be extended to the moduli space of all $N=2$ superspheres, i.e., all $N=2$ super-Riemann surfaces with genus-zero compact body, provided that the following conjecture is true.  }
\end{rema} 

\begin{conj} Any $N=2$ super-Riemann surface with genus-zero compact body is $N=2$ superconformally equivalent to the $N=2$ super-Riemann sphere $S^2 \hat{\mathbb{C}}$.
\end{conj}

However, for the results stated in this paper, we do not need or use this conjecture.  

The Lie supergroup of superconformal automorphisms of $S^2\hat{\mathbb{C}}$ is the group of $N=2$ superprojective transformations $(w,\rho^+,\rho^-) \mapsto (\tilde{w}, \tilde{\rho}^+,\tilde{\rho}^-)$ given by 
\begin{eqnarray}
\tilde{w} &=& \frac{aw + b}{cw + d} + \rho^+ \frac{e^+(\gamma^- w + \delta^-)}{(cw + d)^2} + \rho^+ \frac{(f^+ w + h^+) (\gamma^- w + \delta^-)}{(cw + d)^3}  \label{transform1}\\
& &  + \rho^- \frac{e^-(\gamma^+ w + \delta^+)}{(cw + d)^2} +  \rho^- \frac{(f^-w+ h^-)(\gamma^+ w + \delta^+)}{(cw + d)^3}  \nonumber\\
& & + \rho^+ \rho^- \frac{2\gamma^+ \gamma^-dw  - (\gamma^+ \delta^- + \delta^+ \gamma^-)(cw - d) - 2\delta^+ \delta^-c}{(cw + d)^3} \nonumber \\
\qquad \tilde{\rho}^+ &=& \frac{\gamma^+ w + \delta^+}{cw + d} + \rho^+ \frac{e^+}{cw+d} +  \rho^+ \frac{f^+ w+ h^+}{(cw+d)^2} + \rho^+ \rho^- \frac{\gamma^+ d - \delta^+ c}{(cw+d)^2}\\
\tilde{\rho}^- &=& \frac{\gamma^- w + \delta^-}{cw + d} + \rho^- \frac{e^-}{cw+d} +  \rho^- \frac{f^- w + h^-}{(cw+d)^2} - \rho^+ \rho^- \frac{\gamma^- d - \delta^-c}{(cw+d)^2} \label{transform3}
\end{eqnarray}
for $a,b,c,d,e^\pm, f^\pm, h^\pm \in \bigwedge_\infty^0$ and $\gamma^\pm, \delta^\pm \in \bigwedge_\infty^1$ satisfying
\begin{eqnarray}
ad - bc &=& 1 \label{transform-condition1}\\
e^+ e^-  &=& 1- \gamma^+ \delta^- + \delta^+ \gamma^- \label{transform-condition2}\\
f^\pm &=& \mp  e^\pm \gamma^+ \gamma^-d \label{transform-condition3} \\
h^\pm &=& \pm e^\pm ( \delta^+ \delta^-c - (\gamma^+ \delta^- + \delta^+ \gamma^-)d  \mp \delta^+ \delta^- \gamma^+ \gamma^- d)  \label{transform-condition4}. 
\end{eqnarray} 

Below we give an explicit description of how these $N=2$ superprojective transformations act on $S^2\hat{\mathbb{C}}$.   In Section \ref{NS-algebra-section}, we give a more detailed description of how this Lie supergroup and its action on $S^2\hat{\mathbb{C}}$ is derived and also point out errors in previous claims of descriptions of the group of $N=2$ superprojective transformations acting on $S^2\hat{\mathbb{C}}$, including counterexamples to several of the presentations given in the literature.  

We first note that using conditions (\ref{transform-condition3}) and (\ref{transform-condition4}), we have
\begin{equation} 
(f^\pm w + h^\pm) (\gamma^\mp w + \delta^\mp) = \pm  e^\pm \delta^+ \delta^-\gamma^\mp  (cw +  d )  ,
\end{equation}
and thus $\tilde{w}$ can be simplified {}from  (\ref{transform1}) and in general the $N=2$ superprojective transformations given by (\ref{transform1}) -- (\ref{transform3}) satisfying conditions 
(\ref{transform-condition1}) -- (\ref{transform-condition4}) are equivalent to
\begin{eqnarray}
\tilde{w}  &=& \frac{aw + b}{cw + d} + \rho^+ \frac{e^+(\gamma^- w + \delta^- + \delta^+ \delta^- \gamma^-) }{(cw + d)^2}   \label{short-transform1}\\
& &   + \rho^- \frac{e^-(\gamma^+ w + \delta^+ - \delta^+ \delta^- \gamma^+) }{(cw + d)^2}  \nonumber \\
& & + \rho^+ \rho^- \frac{2\gamma^+ \gamma^-dw  - (\gamma^+ \delta^- + \delta^+ \gamma^-)(cw - d) - 2\delta^+ \delta^-c}{(cw + d)^3} \nonumber \\
\qquad \tilde{\rho}^+ &=&  \frac{\gamma^+ w + \delta^+}{cw + d} + \rho^+ \frac{e^+}{cw+d} \\
& & +  \rho^+ \frac{e^+ (-  \gamma^+ \gamma^-d w + \delta^+ \delta^-c - (\gamma^+ \delta^- + \delta^+ \gamma^-)d - \delta^+ \delta^- \gamma^+ \gamma^- d)}{(cw+d)^2} \nonumber\\
& & + \rho^+ \rho^- \frac{\gamma^+ d - \delta^+ c}{(cw+d)^2} \nonumber 
\end{eqnarray}
\begin{eqnarray}
\qquad \tilde{\rho}^-  &=&  \frac{\gamma^- w + \delta^-}{cw + d} + \rho^- \frac{e^-}{cw+d}  \label{short-transform3} \\ 
& & +  \rho^- \frac{  e^-( \gamma^+ \gamma^-d w -   \delta^+ \delta^-c + (\gamma^+ \delta^- + \delta^+ \gamma^-)d - \delta^+ \delta^- \gamma^+ \gamma^- d)}{(cw+d)^2} \nonumber\\
& & - \rho^+ \rho^- \frac{\gamma^- d - \delta^-c}{(cw+d)^2} \nonumber
\end{eqnarray}
for $a,b,c,d,e^\pm \in \bigwedge_\infty^0$ and $\gamma^\pm, \delta^\pm \in \bigwedge_\infty^1$ satisfying
\begin{eqnarray}
ad - bc &=& 1 \label{short-transform-condition1}\\
e^+ e^-  &=& 1- \gamma^+ \delta^- + \delta^+ \gamma^- \label{short-transform-condition2}. 
\end{eqnarray} 

If $T : S^2\hat{\mathbb{C}} \longrightarrow  S^2\hat{\mathbb{C}}$ is an $N=2$ superprojective transformation, then $T$ can be uniquely expressed by $\sou \circ T \circ \sou^{-1}$ as follows.   We define $T_\sou = \sou \circ T \circ \sou^{-1}$ by
\begin{eqnarray*}
T_\sou : \bigl(\mbox{$\bigwedge_\infty^0$} \smallsetminus \bigl(\{- d_B/c_B \} \times (\mbox{$\bigwedge_\infty^0$})_S \bigr) \bigr) \oplus\mbox{$(\bigwedge_\infty^1)^2$}  &\longrightarrow& \\
& & \hspace{-1.1in} \bigl( \mbox{$\bigwedge_\infty^0$} \smallsetminus \bigl(\{a_B/c_B \} \times
(\mbox{$\bigwedge_\infty^0$})_S\bigr) \bigr) \oplus \mbox{$(\bigwedge_\infty^1)^2$} \\
(w,\rho^+,\rho^-) &\mapsto& (\tilde{w},\tilde{\rho}^+,\tilde{\rho}^-) 
\end{eqnarray*}
with $(\tilde{w},\tilde{\rho}^+,\tilde{\rho}^-)$ given by (\ref{transform1}) -- (\ref{transform3}) satisfying (\ref{transform-condition1}) -- (\ref{transform-condition4}), and we define $T_\nor = \nor \circ T \circ \nor^{-1}$ 
\begin{eqnarray*}
T_\nor: \bigl(\mbox{$\bigwedge_\infty^0$} \smallsetminus \bigl(\{-a_B/b_B \} \times
(\mbox{$\bigwedge_\infty^0$})_S \bigr)\bigr) \oplus \mbox{$(\bigwedge_\infty^1)^2$} &\longrightarrow&
\\
& & \hspace{-1.1in}  \bigl( \mbox{$\bigwedge_\infty^0$} \smallsetminus \bigr(\{d_B/b_B \} \times
(\mbox{$\bigwedge_\infty^0$})_S \bigr) \bigr) \oplus \mbox{$(\bigwedge_\infty^1)^2$}\\  
(w, \rho^+, \rho^-) &\mapsto& (\hat{w},\hat{\rho}^+,\hat{\rho}^-) ,
\end{eqnarray*}
with $(\hat{w},\hat{\rho}^+,\hat{\rho}^-)$ given by
\begin{eqnarray}
\hat{w} &=&  \frac{c + dw}{a + bw} -   i\rho^+ \frac{e^+(\gamma^- +( \delta^- + \delta^+ \delta^- \gamma^-) w)}{(a + bw)^2}  \\
& & - i \rho^- \frac{e^-(\gamma^+  + (\delta^+  - \delta^+ \delta^- \gamma^+)w)}{(a+ bw)^2} \nonumber  \\
& &    +  \rho^+ \rho^- \frac{2\gamma^+ \gamma^-b  - (\gamma^+ \delta^- + \delta^+ \gamma^-)(a-bw)- 2\delta^+ \delta^-aw  }{(a+bw)^3} \nonumber \\
\qquad \hat{\rho}^+ &=&  - i \frac{\gamma^+  + \delta^+ w}{a + bw}   + \rho^+ \frac{e^+  -  \gamma^+ \delta^- -  \delta^+  \gamma^- - \gamma^+ \delta^+ \delta^- \gamma^- }{a + bw} \\
& & +  \rho^+ \frac{e^+( \delta^+ \delta^-aw  -\gamma^+  \gamma^- b + ( \gamma^+ \delta^- +  \delta^+  \gamma^-)a  + \gamma^+ \delta^+ \delta^- \gamma^-a ) }{(a+bw)^2} \nonumber  \\
& &   + i \rho^+ \rho^-  \frac{  \gamma^+ b - \delta^+a}{(a + bw)^2} \nonumber \\
\hat{\rho}^- &=&  -i \frac{\gamma^-  + \delta^- w}{a + bw} + \rho^- \frac{e^- + \delta^+\gamma^- +  \gamma^+\delta^- - \delta^+ \delta^- \gamma^+ \gamma^-  }{a + bw}  \\
& & +  \rho^-\frac{ e^-(- \delta^+\delta^-aw +  \gamma^+\gamma^- b   - (\delta^+\gamma^- +  \gamma^+\delta^-)a + \delta^+ \delta^- \gamma^+ \gamma^-  a  )}{ (a+bw)^2} \nonumber \\
& & - i  \rho^+ \rho^-  \frac{ \gamma^- b-  \delta^-a}{(a + bw)^2}  . \nonumber
\end{eqnarray}
That is $T_\nor(w,\rho^+, \rho^-) =I^{-1} \circ T_\sou \circ I (w, \rho^+, \rho^-)$ for $(w,\rho^+, \rho^-) \in ((\bigwedge_\infty^0)^\times \smallsetminus (\{-a_B/b_B \} \times (\mbox{$\bigwedge_\infty^0$})_S) )\oplus  \mbox{$(\bigwedge_\infty^1)^2$} $.  

Note that $T_\nor(w, \rho, \rho) = (\hat{w},\hat{\rho}^+,\hat{\rho}^-)$ is of the form (\ref{short-transform1}) -- (\ref{short-transform3}) satisfying the conditions (\ref{short-transform-condition1}) and (\ref{short-transform-condition2}).  This can be seen by letting
\begin{equation}
\hat{a} =  d, \quad \hat{b} = c, \quad \hat{c} = b, \quad \hat{d} = a, \quad \hat{\gamma}^\pm = -i \delta^\pm, \quad  \hat{\delta}^\pm = - i \gamma^\pm 
\end{equation}
and
\begin{equation}
\hat{e}^\pm =  e^\pm \mp (\gamma^+ \delta^- + \delta^+ \gamma^-) - \delta^+ \delta^- \gamma^+ \gamma^-.
\end{equation}

Then we define $T: S^2\hat{\mathbb{C}} \longrightarrow S^2 \hat{\mathbb{C}}$ by 
\begin{equation}\label{T1}
T(p) = \left\{
  \begin{array}{ll} 
      \sou^{-1} \circ T_\sou \circ \sou (p) & \mbox{if $p \in
           U_\sou \smallsetminus \sou^{-1}( (\{- d_B/c_B \} \times
            (\bigwedge_\infty^0)_S ) \oplus  \mbox{$(\bigwedge_\infty^1)^2$} )$}, \\  
      \nor^{-1} \circ T_\nor \circ \nor (p) & \mbox{if $p \in
           U_\nor \smallsetminus \nor^{-1}(( \{-a_B/b_B \} \times
             (\bigwedge_\infty^0)_S ) \oplus  \mbox{$(\bigwedge_\infty^1)^2$} )$ }.  
\end{array} \right.
\end{equation}
This defines $T$ for all $p \in S^2\hat{\mathbb{C}}$ unless: \\
(i) $a_B = 0$ and  $p \in \nor^{-1}( (\{0 \} \times (\bigwedge_\infty^0)_S) \oplus  \mbox{$(\bigwedge_\infty^1)^2$}
)$; or \\
(ii) $d_B = 0$ and $p \in \sou^{-1}((\{0 \} \times (\bigwedge_\infty^0)_S) \oplus
\mbox{$(\bigwedge_\infty^1)^2$})$. \\
In case (i), we define
\begin{multline}\label{T2}
T(p) = \sou^{-1} \biggl(  \frac{a + bw}{c + dw} + i\rho^+ \frac{e^+(\gamma^- + \delta^-w)}{(c + dw)^2} + i\rho^+ \frac{(f^+  + h^+w) (\gamma^-  + \delta^-w)}{(c + dw)^3} \\
+ i \rho^- \frac{e^-(\gamma^+  + \delta^+w)}{(c + dw)^2} + i \rho^- \frac{(f^-+ h^-w)(\gamma^+  + \delta^+w)}{(c + dw)^3} \\
- \rho^+ \rho^- \frac{2\gamma^+ \gamma^-d  - (\gamma^+ \delta^- + \delta^+ \gamma^-)(c - dw) - 2\delta^+ \delta^-cw}{(c + dw)^3} ,\\
\frac{\gamma^+  + \delta^+w}{c + dw} + i\rho^+ \frac{e^+}{c+dw} + i \rho^+ \frac{f^+ + h^+w}{(c+dw)^2} - \rho^+ \rho^- \frac{\gamma^+ d - \delta^+ c}{(c+dw)^2},\\
\frac{\gamma^- + \delta^-w}{c + dw} + i\rho^- \frac{e^-}{c+dw} + i \rho^- \frac{f^-  + h^-w}{(c+dw)^2} + \rho^+ \rho^- \frac{\gamma^- d - \delta^-c}{(c+dw)^2} \biggr),
\end{multline} 
for $\nor (p) = (w,\rho^+, \rho^-) = (w_S, \rho^+, \rho^-)$. 

In case (ii), we define 
\begin{multline}\label{T3}
T(p) = \nor^{-1} \Biggl(   \frac{cw + d}{aw + b} -   \rho^+ \frac{e^+(\gamma^-w + \delta^- + \delta^+ \delta^- \gamma^- )}{(aw + b)^2}   -  \rho^- \frac{e^-(\gamma^+w  + \delta^+  - \delta^+ \delta^- \gamma^+)}{(aw+ b)^2} \\
-  \rho^+ \rho^- \frac{2\gamma^+ \gamma^-bw  - (\gamma^+ \delta^- + \delta^+ \gamma^-)(aw-b)- 2\delta^+ \delta^-a  }{(aw+b)^3}, \\
- i \frac{\gamma^+ w + \delta^+ }{aw + b}   -i  \rho^+ \frac{e^+  -  \gamma^+ \delta^- -  \delta^+  \gamma^- - \gamma^+ \delta^+ \delta^- \gamma^- }{aw + b} \\
-i  \rho^+ \frac{e^+( \delta^+ \delta^-a  -\gamma^+  \gamma^- bw + ( \gamma^+ \delta^- +  \delta^+  \gamma^-)aw  + \gamma^+ \delta^+ \delta^- \gamma^-aw ) }{(aw+b)^2}   \\
 - i \rho^+ \rho^-  \frac{  \gamma^+ b - \delta^+a}{(aw + b)^2} , 
 -i \frac{\gamma^-w  + \delta^- }{aw + b} -i \rho^- \frac{e^- + \delta^+\gamma^- +  \gamma^+\delta^- - \delta^+ \delta^- \gamma^+ \gamma^-  }{aw + b}  \\
-i  \rho^-\frac{ e^-(- \delta^+\delta^-a +  \gamma^+\gamma^- bw   - (\delta^+\gamma^- +  \gamma^+\delta^-)aw + \delta^+ \delta^- \gamma^+ \gamma^-  aw  )}{ (aw+b)^2}  \\
+ i  \rho^+ \rho^-  \frac{ \gamma^- b-  \delta^-a}{(aw + b)^2} 
 \Biggr) 
\end{multline} 
for $\sou (p) = (w,\rho^+, \rho^-) = (w_S, \rho+, \rho^-)$.

Note that with this definition, $T$ is uniquely determined by $T_\sou$, i.e., by its value on $\sou (U_\sou)$.  Or equivalently, $T$ is uniquely determined by $T_\nor$, i.e., by its value on $\nor (U_\nor)$.

Let 
\[\mathcal{B}^{r}_\infty = \Bigl\{ (w, \rho^+,\rho^-) \in \mbox{$\bigwedge_\infty^0$} \oplus (\mbox{$\bigwedge_\infty^1$})^2 \; \bigl| \; \Bigl|\frac{1}{w_B} \Bigr| < r  \bigr. \Bigr\} .\] 

\begin{prop}\label{canonicalcriteria}
Any $N=2$ super-Riemann sphere over $\bigwedge_\infty$ with $1+n$ tubes for $n \in \mathbb{Z}_+$ is $N=2$ superconformally equivalent to an $N=2$ super-Riemann sphere with $1+n$ tubes of the form
\begin{multline}\label{preliminarycanonical} 
\Bigl(S^2\hat{\mathbb{C}}; \nor^{-1}(0), \sou^{-1}(z_1,\theta^+_1,\theta^-_1),\dots.., \sou^{-1}(z_{n - 1}, \theta^+_{n - 1}, \theta^-_{n-1}), \sou^{-1}(0);  \Bigr. \\
\bigl(\nor^{-1}(\mathcal{B}^{1/r_0}_0), \Xi_0 \bigr), \bigl(\sou^{-1}( \mathcal{B}^{r_1}_{z_1}), H_1 \circ \sou \bigr),\dots,\\
\Bigl. \bigl(\sou^{-1}(\mathcal{B}^{r_{n - 1}}_{z_{n - 1}}),H_{n - 1} \circ \sou \bigr), \bigl(\sou^{-1}(\mathcal{B}^{r_n}_0),H_n \circ \sou \bigr) \Bigr) ,
\end{multline}
where 
\begin{equation}\label{at-infinity-for-southern-transform}
\Xi_0 |_{\sou^{-1}(\mathcal{B}^{r_0}_\infty)} = H_0 \circ \sou , 
\end{equation}
\begin{equation}\label{distinct-punctures}
(z_1,\theta^+_1,\theta^-_1),\dots.,(z_{n - 1}, \theta^+_{n - 1},\theta^-_{n-1}) \in \mbox{$(\bigwedge_\infty^0)^\times$} \oplus \mbox{$(\bigwedge_\infty^1$})^2, 
\end{equation} 
$(z_j)_B \neq  (z_k)_B$ for $j \neq k$,  $r_0,\dots, r_n \in \mathbb{R}_+$,  and $H_0, \dots,H_n$ are $N=2$ superconformal functions on $\mathcal{B}^{r_0}_\infty$, $\mathcal{B}^{r_1}_{z_1}$,\dots, $\mathcal{B}^{r_{n - 1}}_{z_{n - 1}}$, $\mathcal{B}^{r_n}_0$, respectively, 
such that if we let $H_0 (w, \rho^+,\rho^-) = (\tilde{w_0}, \tilde{\rho_0}^+, \tilde{\rho}^-_0)$, then
\begin{equation}\label{atinfinity}
\lim_{w \rightarrow \infty} H_0 (w, 0,0) = 0, \quad \mbox{and} \quad \lim_{w \rightarrow \infty} \frac{\partial}{\partial \rho^+} w \tilde{\rho}^+_0 = \lim_{w \rightarrow \infty} \frac{\partial}{\partial \rho^-} w
\tilde{\rho}^-_0 = i; 
\end{equation}
\begin{equation}\label{otherpunctures} 
H_j (z_j, \theta^+_j,\theta^-_j) = 0 , \quad \mbox{and} \quad  \lim_{ w \rightarrow z_j} \frac{H_j(w, \theta^+_j,\theta^-_j)}{w - z_j} \in \mbox{$(\bigwedge_\infty^0)^\times$} \oplus \mbox{$(\bigwedge_\infty^1)^2$},  
\end{equation}
for $j = 1,\dots, n - 1$; and 
\begin{equation}\label{atzero}
H_n (0) = 0 , \quad \mbox{and} \quad \lim_{ w \rightarrow 0} \frac{H_n(w,0,0)}{w}  \in \mbox{$(\bigwedge_\infty^0)^\times$} \oplus \mbox{$(\bigwedge_\infty^1)^2$}.
\end{equation}
\end{prop}

\begin{proof} 
Let
\begin{equation}\label{beginningsphere}
S = (S^2\hat{\mathbb{C}};p_0,\dots,p_n;(U_0, \Omega_0),\dots,(U_n, \Omega_n)) 
\end{equation} 
be an $N=2$ super-Riemann sphere with $1+n$ tubes.    We will build in stages an $N=2$ superconformal transformation that sends $S$ to an $N=2$ super-Riemann sphere with $1+n$ tubes of the form given by the proposition.

If $\sou(p_0)$ has even coordinate equal to zero, let $T_1  : S^2\hat{\mathbb{C}} \rightarrow  S^2\hat{\mathbb{C}}$ be given by $\sou \circ T_1 \circ \sou^{-1} (w, \rho^+,\rho^-) =  (T_1)_\sou (w, \rho^+, \rho^-)  = I(w, \rho^+, \rho^-) = (1/w, i\rho^+/w,i\rho^-/w)$.  Otherwise, $p_0 \in U_\nor$, and we take $T_1$ to be the identity.

Now $\nor \circ T_1(p_0) = (u_0,v_0^+, v_0^-)$ for some $(u_0,v_0^+, v_0^-) \in \bigwedge_\infty^0 \oplus (\bigwedge_\infty^1)^2$.  Let $T_2  : S^2\hat{\mathbb{C}} \rightarrow  S^2\hat{\mathbb{C}}$ be given by $\nor \circ T_2 \circ \nor^{-1} (w, \rho^+,\rho^-) =  (T_2)_\nor (w, \rho^+, \rho^-)  = (w - u_0 - \rho^+ v_0^- - \rho^- v_0^+, \rho^+ - v_0^+, \rho^- - v_0^-)$.

Note that $T_2 \circ T_1 (p_n) \in U_\sou$; let $\sou \circ T_2 \circ T_1 (p_n) = (u_n, v^+_n,v^-_n)$ and let $T_3 : S^2\hat{\mathbb{C}} \rightarrow  S^2\hat{\mathbb{C}} $ be given by $\sou \circ T_3 \circ \sou^{-1} (w, \rho^+,\rho^-) =  (T_3)_\sou (w, \rho^+, \rho^-)  =  ( w - u_n -  \rho^+ v_n^- -  \rho^-v_n^+  , \rho^+ - v_n^+, \rho^- - v_n^-)$. 

Let $T = T_3 \circ T_2 \circ T_1$.  Then $T(S)$ now has the outgoing puncture at $\nor^{-1} (0)$, i.e., at infinity, and the last incoming puncture at $\sou^{-1} (0)$, i.e., at zero.    Now we need to fix information about the local coordinate vanishing at infinity.    The local coordinate at infinity is now given by 
$\Omega_0 \circ T^{-1}$.
Let
\[\Omega_0  \circ T^{-1} \circ \sou^{-1} (w, \rho^+,\rho^-) = (\tilde{w}, \tilde{\rho}^+, \tilde{\rho}^-) .\]
Then 
\[  \lim_{w \rightarrow \infty} \frac{\partial}{\partial \rho^\pm } w \tilde{\rho}^\pm = i e^\pm \]
uniquely determines $e^\pm \in (\bigwedge_\infty^0)^\times$.  Let $T_4 :  S^2\hat{\mathbb{C}} \rightarrow  S^2\hat{\mathbb{C}} $ be given by $\sou \circ T_4 \circ \sou^{-1} (w, \rho^+,\rho^-) =  (T_4)_\sou (w, \rho^+, \rho^-)  =  ( (e^+ e^-)^{-1} w, (e^-)^{-1}\rho^+ , (e^+)^{-1} \rho^-)$.

Let $F = T_4 \circ T$.  Then $F$ is an $N=2$ superconformal equivalence {}from $S$ to 
\begin{multline}\label{nexttolastsphere}
( S^2\hat{\mathbb{C}} ; \nor^{-1} (0), \sou^{-1}(z_1,\theta^+_1, \theta^-_1),\dots.., \sou^{-1}(z_{n - 1}, \theta^+_{n - 1}, \theta^-_{n-1}), \sou^{-1}(0); \\ ( F(U_0), \Omega_0 \circ F^{-1} ),\dots,(F(U_n), \Omega_n \circ F^{-1}))  
\end{multline} 
where
\[(z_j, \theta^+_j, \theta^-_j) = \sou \circ F(p_j),\] for $j = 1,\dots, n - 1$.  Choose $r_0,\dots,r_n \in \mathbb{R}_+$ such that
\[\mathcal{B}^{r_0}_{\infty} \subset \sou \circ F(U_0) \]
\[\mathcal{B}^{r_j}_{z_j} \subset \sou \circ F(U_j) , \quad j = 1,\dots, n - 1,  \]
\[\mathcal{B}^{r_n}_0 \subset \sou  \circ F(U_n) . \]
Then the $N=2$ super-Riemann sphere with tubes (\ref{nexttolastsphere}) is $N=2$ superconformally equivalent to
\begin{multline*}
\Bigl( S^2 \hat{\mathbb{C}} ; \nor^{-1} (0), \sou^{-1}(z_1,\theta^+_1, \theta^-_1),\dots,  \sou^{-1} (z_{n - 1}, \theta^+_{n - 1}, \theta^-_{n-1}), \sou^{-1}(0); \\
\bigl(\sou^{-1} (\mathcal{B}^{r_0}_\infty) \cup \nor^{-1} (\{0 \} \times (\mbox{$\bigwedge_\infty$})_S) , \Omega_0 \circ F^{-1} \bigr), \bigl(\sou^{-1}(\mathcal{B}^{r_1}_{z_1}),\Omega_1 \circ F^{-1} \bigr), \dots,\\
\bigl(\sou^{-1}(\mathcal{B}^{r_{n - 1}}_{z_{n - 1}}), \Omega_{n - 2} \circ F^{-1} \bigr), \bigl(\sou^{-1}(\mathcal{B}^{r_n}_0), \Omega_n \circ F^{-1}\bigr)\Bigr) 
\end{multline*} 
where
\begin{eqnarray*}
H_0 &=& \left. \Omega_0 \circ F^{-1}  \circ \sou^{-1} \right|_{\mathcal{B}^{r_0}_\infty} , \\
H_j &=& \left. \Omega_j \circ F^{-1}  \circ \sou^{-1} \right|_{\mathcal{B}^{r_j}_{z_j}} , \quad j = 1,\dots, n-1 , \\ 
H_n &=& \left. \Omega_n \circ F^{-1} \circ \sou^{-1}  \right|_{\mathcal{B}^{r_n}_0}  
\end{eqnarray*}
satisfy (\ref{atinfinity}), (\ref{otherpunctures}), and (\ref{atzero}), respectively.
\end{proof}

An $N=2$  super-Riemann sphere with $1+n$ tubes, for $n \in \Z$,  of the form (\ref{preliminarycanonical}) is called a {\em canonical $N=2$ supersphere with $1+n$ tubes}.

\begin{rema}\label{remark-giving-power-series}
{\em A canonical $N=2$ supersphere with $1+n$ tubes, for $n \in \mathbb{Z}_+$, is determined by $(z_1,\theta^+_1,\theta^-_1),\dots, (z_{n - 1}, \theta^+_{n - 1}, \theta^-_{n-1}) \in (\bigwedge_\infty^0)^\times \oplus (\bigwedge_\infty^1)^2$, $r_0,\dots,r_n \in \mathbb{R}_+$ and $N=2$ superconformal functions $H_0,\dots,H_n$ satisfying (\ref{atinfinity}),  (\ref{otherpunctures}), and (\ref{atzero}), respectively.  Consider the $N=2$ superconformal power series obtained by expanding the superconformal functions $H_0,\dots,H_n$ around $\nor^{-1} (0) = \infty$, $(z_1,\theta^+_1,\theta^-_1),\dots,(z_{n-1},\theta^+_{n-1},\theta^-_{n-1})$, and $0$, respectively.  For $k =0, \dots, n$, we will denote by $H_k$ both the superconformal function $H_k$ and its power series expansion.  Moreover, if $H(w,\rho^+,\rho^-) = (\tilde{w}, \tilde{\rho}^+, \tilde{\rho}^-)$, then define $H^\pm (w,\rho^+, \rho^-) = \tilde{\rho}^\pm$.   By the fact that a power series expansion about zero of an $N=2$ superconformal superfunction vanishing at zero must be of the form (\ref{powerseries1}) -- (\ref{powerseries2}) and the power series expansion about zero of an $N=2$ superconformal superfunction vanishing at infinity composed with $I$ must also be of the form (\ref{powerseries1}) -- (\ref{powerseries2}), by the conditions (\ref{atinfinity}), (\ref{otherpunctures}), and (\ref{atzero}) and by the fact that the $H_k$'s are one-to-one as superanalytic functions on their domains, we have 
\begin{equation}\label{power-series-at-infinity2}
\rho^\mp H^\pm_0(w,\rho^+,\rho^-) =  \rho^\mp \biggl( \sum_{j \in \Z}m_{j-\frac{1}{2}}^{\pm,(0)} w^{-j} + \rho^\pm \Bigl( \frac{i}{w} + \sum_{j \in \mathbb{Z}_+}a_j^{\pm,(0)} w^{-j-1} \Bigr) \biggr), 
\end{equation}
\begin{equation}
\rho^\mp H^\pm_k(w,\rho^+,\rho^-) = \biggl. \rho^\mp \biggl( \sum_{j \in \Z} m_{j-\frac{1}{2}}^{\pm,(k)} w^j + \rho^\pm \sum_{j \in \mathbb{N}} a_j^{\pm,(k)} w^j \biggr) 
\biggr|_{ s_k(w,\rho^+,\rho^-)}   \label{power-series-at-ith-puncture} 
\end{equation}
for $k = 1,\dots,n$, where $(z_n, \theta^+_n,\theta^-_n) = 0$, $a_0^{\pm,(k)} \in (\bigwedge_\infty^0)^\times$, $a_j^{\pm,(k)} \in \bigwedge_\infty^0$, $m_{j-\frac{1}{2}}^{\pm,(k)} \in \bigwedge_\infty^1$, for $j\in \Z$, $s_k(w,\rho^+,\rho^-) =  (w - z_k -\rho^+ \theta^-_k - \rho^- \theta^+_k, \rho^+ - \theta^+_k, \rho^- - \theta^-_k)$ and $H_k$ are superconformal according to (\ref{superconformal-condition1}) -- (\ref{superconformal-condition4}).  Then given $H_k$ of the form (\ref{power-series-at-infinity2}), and (\ref{power-series-at-ith-puncture}), respectively, that are superconformal and vanishing at the punctures (at infinity for $k=0$, at $(z_k, \theta_k^+, \theta_k^-)$ for $k = 1,\dots n-1$, and at zero for $k = n$), these conditions uniquely determine the $H_k$'s by Remark \ref{determining-superconformal-remark}.

Thus a canonical $N=2$ supersphere with $1+n$ tubes, for $n \in \Z$, can be denoted by
\begin{equation}
((z_1,\theta^+_1,\theta^-_1),\dots, (z_{n - 1}, \theta^+_{n - 1},\theta^-_{n-1});r_0,\dots,r_n;H_0,\dots,H_n)  
\end{equation} 
where $H_0,\dots,H_n$ are power series of the form (\ref{power-series-at-infinity2}) and (\ref{power-series-at-ith-puncture}), respectively, which vanish at the punctures and satisfy the superconformal conditions (\ref{superconformal-condition1}) -- (\ref{superconformal-condition4}). }
\end{rema}

\begin{rema}\label{infinite-variables2}
{\em {}From Remark \ref{remark-giving-power-series} above, we can readily see that a point in the moduli space of $N=2$ super-Riemann spheres with $1+n$ tubes, for $n \in \mathbb{Z}_+$, will in general depend on an infinite number of odd variables -- the $\theta^\pm_1,\dots,\theta^\pm_{n-1}$ and the $m_{j-1/2}^{\pm,(k)} \in \bigwedge_\infty^1$,  for $k=0,\dots,n$, and $j\in \Z$. }
\end{rema}

\begin{prop}\label{equalspheres}
Two canonical $N=2$ superspheres with $1+n$ tubes, for $n \in \mathbb{Z}_+$,
\begin{equation}\label{stupidsphere1}
((z_1,\theta^+_1,\theta^-_1),\dots, (z_{n - 1}, \theta^+_{n - 1},\theta^-_{n-1}); r_0,\dots,r_n;H_0,\dots,H_n) 
\end{equation}
and 
\begin{equation}\label{stupidsphere2}
((\hat{z}_1,\hat{\theta}^+_1,\hat{\theta}^-_1),\dots, (\hat{z}_{n - 1}, \hat{\theta}^+_{n - 1},\hat{\theta}^-_{n-1});\hat{r}_0,\dots,\hat{r}_n; \hat{H}_0,\dots,\hat{H}_n) 
\end{equation}
are $N=2$ superconformally equivalent if and only if $(z_j,\theta^+_j,\theta^-_j) = (\hat{z}_j, \hat{\theta}^+_j, \hat{\theta}^-_j)$ for $j = 1,\dots, n-1$, and $H_j = \hat{H_j}$, for $j = 0,\dots,  n$, as superconformal power series.  
\end{prop}

\begin{proof} Let $F$ be an $N=2$ superconformal equivalence {}from (\ref{stupidsphere1}) to (\ref{stupidsphere2}).  The conclusion of the proposition is equivalent to the assertion that $F$ must be the identity map on $S^2\hat{\mathbb{C}}$.  By definition $F$ is a superconformal automorphism of $S^2\hat{\mathbb{C}}$, i.e., an $N=2$ superprojective transformation.  Also by definition we have
\begin{equation}\label{eq1}
F_\sou (0) = \sou \circ F \circ \sou^{-1} (0) = 0 , 
\end{equation}
\begin{equation}\label{eq2}
F_\nor (0) = \nor \circ F \circ \nor^{-1} (0) = 0 , 
\end{equation}
\begin{equation}\label{eq3}
\hat{H}_0 |_{\mathcal{B}^{\min (r_0, \hat{r}_0)}_\infty} = H_0 \circ F_\sou^{-1} |_{\mathcal{B}^{\min (r_0, \hat{r}_0)}_{\infty}} .  
\end{equation} 
{}From (\ref{eq1}) and (\ref{eq2}) and the fact that $F$ is a superprojective transformation, we obtain
\begin{equation}\label{eq4}
F_\sou (w, \rho^+, \rho^-) = (e^+ e^- w, e^+ \rho^+, e^- \rho^-)
\end{equation}
for some $e^\pm \in (\bigwedge_\infty^0)^\times$.  Let $H_0 (w,\rho^+, \rho^-) = (\tilde{w}_0, \tilde{\rho}^+_0, \tilde{\rho}^-_0)$ and $\hat{H}_0 (w,\rho^+, \rho^-) = (\hat{w}_0, \hat{\rho}^+_0, \hat{\rho}^-_0)$. {}From (\ref{atinfinity}), we know that
\[\lim_{w \rightarrow \infty} \frac{\partial}{\partial \rho^\pm} w 
\tilde{\rho}^\pm_0 (w, \rho^+, \rho^-) = \lim_{w \rightarrow \infty} \frac{\partial}{\partial \rho^\pm} w \hat{\rho}^\pm_0 (w,\rho^+,\rho^-) = i . \] 
Thus by (\ref{eq3}) and (\ref{eq4})
\begin{eqnarray*}
i &=& \lim_{w \rightarrow \infty} \frac{\partial}{\partial \rho^\pm} w \hat{\rho}_0^\pm (w, \rho^+, \rho^-) =\lim_{w \rightarrow \infty} \frac{\partial}{\partial \rho^\pm} w (H_0 \circ F_\sou^{-1})^\pm (w, \rho^+, \rho^-)\\
&=& \lim_{w \rightarrow \infty} \frac{\partial}{\partial \rho^\pm} w \tilde{\rho}^\pm_0 \Bigl(\frac{w}{e^+ e^-} , \frac{\rho^+}{e^+}, \frac{\rho^-}{e^-}\Bigr) \\ 
&=& \frac{i}{e^\pm},
\end{eqnarray*}
i.e., $e^\pm = 1$. Thus $F$ must be the identity map of $S^2\hat{\mathbb{C}}$.
 \end{proof}

For $N=2$ super-Riemann spheres with one tube, we have:

\begin{prop}\label{canonical-criteria-for-one-tube}
Any $N=2$ super-Riemann sphere with one tube is superconformally equivalent to an $N=2$ super-Riemann sphere with one tube of the form
\begin{equation}\label{onetube}
\bigl(S^2\hat{\mathbb{C}}; \nor^{-1} (0) ; (\nor^{-1} (\mathcal{B}^{1/r}_0), \Xi) \bigr) 
\end{equation}
where $\Xi |_{(\sou^{-1} (\mathcal{B}^{r}_{\infty})} = H \circ \Delta$, and $H(w,\rho^+, \rho^-) = (\tilde{w}, \tilde{\rho}^+, \tilde{\rho}^-)$ satisfies
\begin{eqnarray}
\lim_{w \rightarrow \infty} H(w,0, 0) &=&  0 \label{1-tube1}\\
\lim_{w \rightarrow \infty} \frac{\partial}{\partial \rho^+} w \tilde{\rho}^+ = \lim_{w \rightarrow \infty} \frac{\partial}{\partial \rho^-} w
\tilde{\rho}^- &=& i  \label{1-tube2} \\
\lim_{w \rightarrow \infty} w^2(\tilde{w}(w,0,0) - w^{-1}) =   \lim_{w \rightarrow \infty} w \tilde{\rho}^\pm(w,0,0)  &=& 0 . \label{1-tube3}
\end{eqnarray}
These conditions are equivalent to the condition that $H$ can be expanded in a power series about infinity of the form (\ref{power-series-at-infinity2}) with $a_1^{+,(0)} = - a_1^{-,(0)}$ and $m_{1/2}^{\pm,(0)} = 0$.
\end{prop}

\begin{proof} Given an $N=2$ super-Riemann sphere with one tube, 
\begin{equation}\label{onetube2}
S= (S^2\hat{\mathbb{C}}; p;(U, \Omega )), 
\end{equation}
if $\sou(p)$ has even coordinate equal to zero, let $T_1  : S^2\hat{\mathbb{C}} \rightarrow  S^2\hat{\mathbb{C}}$ be given by $\sou \circ T_1 \circ \sou^{-1} (w, \rho^+,\rho^-) =  (T_1)_\sou (w, \rho^+, \rho^-)  = I(w, \rho^+, \rho^-) = (1/w, i\rho^+/w,i\rho^-/w)$.  Otherwise, $p \in U_\nor$, and we take $T_1$ to be the identity.

Now $\nor \circ T_1(p) = (u,v^+, v^-)$ for some $(u,v^+, v^-) \in \bigwedge_\infty^0 \oplus (\bigwedge_\infty^1)^2$.  Let $T_2  : S^2\hat{\mathbb{C}} \rightarrow  S^2\hat{\mathbb{C}}$ be given by $\nor \circ T_2 \circ \nor^{-1} (w, \rho^+,\rho^-) =  (T_2)_\nor (w, \rho^+, \rho^-)  = (w - u - \rho^+ v^- - \rho^- v^+, \rho^+ - v^+, \rho^- - v^-)$.

Let $T =  T_2 \circ T_1$.  Then $T(S)$ now has the outgoing puncture at $\nor^{-1} (0)$, i.e., at infinity.    Now we need to fix certain properties of the local coordinate vanishing at infinity.    The local coordinate at infinity is now given by  $\Omega \circ T^{-1}$.
Let
\[\Omega  \circ T^{-1} \circ \sou^{-1} (w, \rho^+,\rho^-) = (\tilde{w}, \tilde{\rho}^+, \tilde{\rho}^-) .\]
Then 
\[  \lim_{w \rightarrow \infty} \frac{\partial}{\partial \rho^\pm } w \tilde{\rho}^\pm = i e^\pm \]
uniquely determines $e^\pm \in (\bigwedge_\infty^0)^\times$.  Let $T_3 :  S^2\hat{\mathbb{C}} \rightarrow  S^2\hat{\mathbb{C}} $ be given by $\sou \circ T_3 \circ \sou^{-1} (w, \rho^+,\rho^-) =  (T_3)_\sou (w, \rho^+, \rho^-)  =  ( (e^+ e^-)^{-1} w, (e^-)^{-1}\rho^+ , (e^+)^{-1} \rho^-)$.

Let $F = T_3 \circ T$.  Then $F$ is an $N=2$ superconformal equivalence {}from $S$ to 
\begin{equation}
( S^2\hat{\mathbb{C}} ; \nor^{-1} (0); ( F(U), \Omega \circ F^{-1} ))  
\end{equation} 
and $\Omega \circ F^{-1}  \circ \sou^{-1} = H_0$ where $H_0$ satisfies (\ref{atinfinity}), i.e., has a power series expansion of the form (\ref{power-series-at-infinity2}).  Let $H_0(w,\rho^+, \rho^-) = (\tilde{w}_0, \tilde{\rho}_0^+, \tilde{\rho}_0^-)$.  Then
\begin{eqnarray}
\lim_{w \rightarrow \infty} w^2(\tilde{w}_0(w,0,0) - w^{-1}) &=& a \\
\lim_{w \rightarrow \infty} w \tilde{\rho}_0^\pm(w,0,0) &=& m^\pm
\end{eqnarray}
uniquely determine $a \in \bigwedge_\infty^0$ and $m^\pm \in \bigwedge_\infty^1$.  Let $(T_4)_\sou (w,\rho^+, \rho^-) = (w - a - i \rho^+ m^- - i \rho^- m^+, -im^+ + \rho^+, -im^- + \rho^-)$.  Then $H = \Omega \circ F^{-1} \circ T_4^{-1} \circ \sou^{-1}$ satisfies conditions (\ref{1-tube1}) - (\ref{1-tube3}), and there exists some $r  \in \mathbb{R}_+$ such that $H$ is convergent in  $\mathcal{B}_\infty^{r}$.  Thus $T_4 \circ F$ is an $N=2$ superconformal equivalence {}from $S$ to an $N=2$ super-Riemann sphere with one tube of the required form.
 \end{proof}

An $N=2$ super-Riemann sphere with one tube of the form (\ref{onetube}) is called a {\em canonical $N=2$ supersphere with one tube}.  A canonical $N=2$ supersphere with one tube is determined by $r \in \mathbb{R}_+$ and an $N=2$ superconformal power series $H$ satisfying (\ref{1-tube1}) -- (\ref{1-tube3}) (or equivalently of the form (\ref{power-series-at-infinity2})  with $a_1^{+,(0)} = -a_1^{-,(0)}$ and $m_{1/2}^{\pm,(0)} = 0$), and can be denoted by $(r;H)$. The following proposition can be proved similarly to Proposition \ref{equalspheres}.

\begin{prop}\label{equal-spheres-with-one-tube}
Two canonical $N=2$ superspheres with one tube $(r;H)$ and $(\hat{r}; \hat{H})$ are $N=2$ superconformally equivalent if and only if $H = \hat{H}$.  
\end{prop} 

{}From Propositions \ref{canonicalcriteria}, \ref{equalspheres}, \ref{canonical-criteria-for-one-tube}, and \ref{equal-spheres-with-one-tube} we have the following theorem:

\begin{thm}\label{bijection} 
There is a bijection between the set of canonical $N=2$ superspheres with tubes and the moduli space of $N=2$ super-Riemann spheres with tubes.  In particular, the moduli space of $N=2$ super-Riemann spheres with $1+n$ tubes, for $n \in \Z$, can be identified with all $(3n+ n)$-tuples $((z_1,\theta^+_1,\theta^-_1),\dots,(z_{n-1}, \theta^+_{n-1}, \theta^-_{n-1});$ $H_0,\dots,H_n)$ satisfying $(z_j,\theta^+_j,\theta^-_j) \in (\bigwedge_\infty^0)^\times \oplus (\bigwedge_\infty^1)^2$, with $(z_j)_B \neq (z_k)_B$ if $j \neq k$,  for $j,k = 1,\dots,n-1$, and such that $H_0,\dots,H_n$ vanish at the corresponding punctures, are of the form (\ref{power-series-at-infinity2}) and (\ref{power-series-at-ith-puncture}), respectively, and are absolutely  convergent in neighborhoods of $\infty$, $(z_1,\theta^+_1, \theta^-_1),\dots,$ $(z_{n-1},\theta^+_{n-1},\theta^-_{n-1})$, and $0$, respectively.  The moduli space of $N=2$ super-Riemann spheres with one tube can be identified with the set of all power series $H_0$ of the form (\ref{power-series-at-infinity2}) such that $a_1^{+,(0)} = -a_1^{-,(0)}$ and $m_{1/2}^{\pm,(0)} = 0$ and such that $H_0$ vanishes at infinity and is absolutely convergent in a neighborhood of infinity. 
\end{thm}

\section{Infinitesimal $N=2$ superconformal transformations}\label{infinitesimals-section}

In this section, we develop a formal theory of infinitesimal $N = 2$ superconformal transformations based on a representation of the $N = 2$ Neveu-Schwarz algebra of superconformal symmetries in terms of superderivations.  The material in this section is algebraic and independent of the supergeometry studied in the previous sections.  However, the results of this section do of course have geometric motivation and meaning, and will be applied to the supergeometric setting in Section \ref{reformulate-moduli-section}.  

First we recall two generalizations of the ``automorphism property'' {}from \cite{FLM} and develop a formal supercalculus of $N=2$ superconformal power series.  Then using formal exponentiation, we characterize certain formal superconformal local coordinate maps in terms of exponentials of superderivations with infinitely many formal variable coefficients.  

An expression of the form $e^x$ denotes the formal exponential series in $x$. In the proof of Propositions \ref{superconformal}, \ref{auto2} and \ref{auto3}, we will need the following proposition which is a generalization of the ``automorphism property" of \cite{FLM} and was proved by the author in \cite{B-memoir}.

\begin{prop}\label{conformalproof} (\cite{B-memoir})
Let $A$ be a superalgebra; let $h \in A^0$; let $u,v \in A$; let $\mathcal{T} \in (\mbox{\em Der} \; A)^0$; and let $y$ be a formal variable commuting with $A$.  Then
\begin{eqnarray}
e^{y(h + \mathcal{T})} \cdot (u v) &=& (e^{y\mathcal{T}} \cdot u) (e^{y(h + \mathcal{T})} \cdot v) \\
&=&(e^{y(h + \mathcal{T})} \cdot u) (e^{y\mathcal{T}} \cdot v) .  \nonumber
\end{eqnarray}
\end{prop}

Let $R$ be a superalgebra over $\mathbb{C}$ (with identity).  Let $x$ be a formal variable which commutes with all elements of $R$, and let $\varphi^\pm$ be formal variables which commute with $x$ and elements of $R^0$ and anti-commute with elements of $R^1$, themselves and each other.  In general, we will use the term {\it even formal variable} to denote a formal variable which commutes with all formal variables and with all elements in any coefficient algebra.   We will use the term {\it odd formal variable} to denote a formal variable which anti-commutes with all odd elements and commutes with all even elements in any coefficient algebra, and in addition, odd formal variables will all
anti-commute with each other.  Consequently, an odd formal variable has the property that its square is zero.   

For a vector space $V$, and for even formal variables $x_1$, $x_2$,\dots, and odd formal variables $\varphi_1$, $\varphi_2$,\dots, consider the spaces 
\begin{multline*}
V[[x_1,\dots,x_n]][\varphi_1,\dots,\varphi_m]  \\
= \; \Bigl\{ \Bigl. \sum_{ \begin{tiny}
\begin{array}{c} k_1,\dots, k_n \in \mathbb{N}\\
l_1,\dots, l_m \in \mathbb{Z}_2
\end{array} \end{tiny}} \! a_{k_1,\dots,k_n, l_1,\dots, l_m} x_1^{k_1}\cdots x_n^{k_n} \varphi_1^{l_1} \cdots \varphi_m^{l_m} \; \Bigr| \; a_{k_1,\dots,k_n, l_1,\dots, l_m} \in V \Bigr\} 
\end{multline*}
and
\begin{multline*}
V((x))[\varphi^+,\varphi^-] = \Bigl\{ \Bigl. \sum_{n = N}^{\infty} a_n x^n + \varphi^+ \! \sum_{n = N}^{\infty} b_n x^n + \varphi^- \! \sum_{n = N}^{\infty} c_n x^n \; \Bigr| \; N \in \mathbb{Z} , \; \\
a_n, b_n, c_n \in V \Bigr\} \subset V[[x,x^{-1}]] [\varphi^+, \varphi^-] .
\end{multline*}  
Then $R((x))[\varphi^+,\varphi^-]$ is a superalgebra as is $R((x^{-1}))[\varphi^+,\varphi^-]$ with $\mathbb{Z}_2$-grading given by 
\begin{eqnarray*}
R((x))[\varphi^+,\varphi^-]^0 &=& R^0((x)) \oplus \varphi^+ R^1((x)) \oplus \varphi^- R^1((x)) \oplus \varphi^+ \varphi^- R^0((x)) \\
R((x))[\varphi^+,\varphi^-]^1  &=& R^1((x)) \oplus \varphi^+ R^0((x)) \oplus \varphi^- R^0((x)) \oplus  \varphi^+ \varphi^- R^1((x))  ,
\end{eqnarray*}  
and similarly for $R((x^{-1}))[\varphi^+,\varphi^-]$ and $R[[x^{-1}, x]][\varphi^+, \varphi^-]$.  

Define 
\begin{equation}\label{define-D}
D^\pm = \frac{\partial}{\partial \varphi^\pm} + \varphi^\mp \frac{\partial}{\partial x} . 
\end{equation}
Then $D^\pm$ are odd derivations in $\mbox{Der}(R((x))[\varphi^+, \varphi^-])$ and in $\mbox{Der}(R((x^{-1}))[\varphi^+, \varphi^-])$.  Furthermore, $D^+$ and $D^-$ satisfy the super-Leibniz rule (\ref{leibniz}) for the product of any two elements in $R[[x,x^{-1}]] [\varphi^+, \varphi^-]$ if that product is well defined.   Note that
\[ (D^\pm)^2 = 0 \qquad \mbox{and} \qquad [D^+,D^-] = D^+D^- + D^-D^+ = 
2 \frac{\partial}{\partial x} . \]

A superanalytic $(1,2)$-superfunction $H(z, \theta^+, \theta^-)$ {}from a DeWitt open set in $(\bigwedge_{*>0}^0 \oplus (\bigwedge_{*>0}^1)^2)$  to $(\bigwedge_{*>0}^0 \oplus (\bigwedge_{*>0}^1)^2)$ has a Laurent expansion about $z$ and $\theta^\pm$ which is an element of $(\bigwedge_{*>0}^0 \oplus (\bigwedge_{*>0}^1)^2) [[z,z^{-1}]] [\theta^+, \theta^-]$.  Taking a general coefficient superalgebra $R$, we can write a corresponding formal superfunction in one even formal variable and two odd formal variables over $R$ as $H(x,\varphi^+,\varphi^-) = (\tilde{x}, \tilde{\varphi}^+, \tilde{\varphi}^-)$ with
\begin{eqnarray*}
\tilde{x} &=& f(x) + \varphi^+ \xi^+(x) + \varphi^- \xi^-(x) + \varphi^+ \varphi^- f^{+,-}(x)\\
\tilde{\varphi}^+ &=& \psi^+(x) + \varphi^+ g^+(x) + \varphi^- h^-(x) + \varphi^+ \varphi^- \phi^+(x)\\
\tilde{\varphi}^- &=& \psi^-(x) + \varphi^+ h^+(x) + \varphi^- g^-(x) + \varphi^+\varphi^- \phi^-(x)
\end{eqnarray*}
where $\tilde{x} \in R[[x,x^{-1}]][\varphi^+,\varphi^-]^0$ and $\tilde{\varphi}^\pm \in  R[[x,x^{-1}]][\varphi^+,\varphi^-]^1$, i.e., where $f(x)$, $f^{+,-}(x)$, $g^\pm(x)$, $h^\pm(x) \in R^0[[x,x^{-1}]]$, and $\xi^\pm(x)$, $\psi^\pm(x)$, $\phi^\pm(x) \in  R^1[[x,x^{-1}]]$.  

In Section \ref{superconformal-section}, the operators $D^\pm = \frac{\partial}{\partial \theta^\pm} + \theta^\mp \frac{\partial}{\partial z}$ were used to define the notion of $N=2$ superconformal $(1,2)$-superfunction over $\bigwedge_{*>0}$ which is a $(\bigwedge_{*>0}^0 \oplus (\bigwedge_{*>0}^1)^2)$-valued superanalytic $(1,2)$-function with the condition that it transform $D^+$ and $D^-$ homogeneously of degree one.   This is equivalent to the conditions (\ref{superconformal-condition1}) -- (\ref{non-zero-condition}).  Thus formally, we  define a series 
\[H(x,\varphi^+,\varphi^-) = (\tilde{x}, \tilde{\varphi}^+,\tilde{\varphi}^-) \in R[[x, x^{-1}]][\varphi^+, \varphi^-]^0 \oplus (R[[x, x^{-1}]][\varphi^+, \varphi^-]^1)^2\] 
to be {\it formally $N=2$ superconformal} if  
\begin{eqnarray}
\qquad \quad \tilde{x} &=& f(x) + \varphi^+ g^+(x) \psi^-(x) + \varphi^- g^-(x) \psi^+(x) + \varphi^+ \varphi^- (\psi^+(x)\psi^-(x))' \label{formal-superconformal-condition1} \\
\tilde{\varphi}^+ &=& \psi^+(x) + \varphi^+ g^+(x) + \varphi^+ \varphi^-(\psi^+)'(x) \label{formal-superconformal-condition2} \\
\tilde{\varphi}^- &=& \psi^-(x) + \varphi^- g^-(x) - \varphi^+ \varphi^- (\psi^-)'(x) \label{formal-superconformal-condition3}
\end{eqnarray}
and
\begin{equation}\label{formal-superconformal-condition4}
f'(x) = (\psi^+)'(x)\psi^-(x) - \psi^+(x) (\psi^-)'(x) + g^+(x) g^-(x) ,
\end{equation}
with
\begin{eqnarray}\label{formal-non-zero-condition}
g^+(x) + (\psi^+) ' (x)  \equiv \! \! \! \! \negthickspace /  \ \  0 \qquad \mathrm{and} \qquad   
g^-(x) +  (\psi^-) ' (x)   \equiv \! \! \! \! \negthickspace /  \ \  0 . 
\end{eqnarray}

Therefore in general a formal $N=2$ superconformal series is uniquely determined by three even formal series, $f(x)$ and $g^\pm(x)$, and two odd formal series $\psi^\pm(x)$ satisfying the conditions (\ref{formal-superconformal-condition4}) and (\ref{formal-non-zero-condition}).   

We also note that for $H(x,\varphi^+,\varphi^-) =  (\tilde{x},\tilde{\varphi}^+, \tilde{\varphi}^-)$, the conditions (\ref{formal-superconformal-condition1}) -- (\ref{formal-superconformal-condition4}) are equivalent to the conditions
\begin{equation}\label{nice-superconformal-condition}
D^+ \tilde{\varphi}^- = D^- \tilde{\varphi}^+ = D^+\tilde{x} - \tilde{\varphi}^- D^+ \tilde{\varphi}^+ = D^-\tilde{x} - \tilde{\varphi}^+ D^-\tilde{\varphi}^- =0 
\end{equation}     
and the conditions (\ref{formal-non-zero-condition}) are equivalent to  the series $D^+ \tilde{\varphi}^+$ and $D^- \tilde{\varphi}^-$ not being identically zero. 

In Section \ref{moduli-section}, we began the study of the moduli space of $N=2$ super-Riemann spheres with punctures and local superconformal coordinates vanishing at the punctures.  The punctures on an $N=2$ super-Riemann sphere with tubes over $\bigwedge_\infty$ can be thought of as being at $0 \in \bigwedge_\infty^0 \oplus (\bigwedge_\infty^1)^2$,  a non-zero point in $\bigwedge_\infty^0 \oplus (\bigwedge_\infty^1)^2$, or at  $\infty$.  Since we can always shift a non-zero point in $\bigwedge_\infty^0 \oplus (\bigwedge_\infty^1)^2$ to zero, all local superconformal coordinates vanishing at the punctures can be expressed as power series vanishing at zero or infinity.  Thus we want to study in more detail certain formal $N=2$ superconformal power series $H(x, \varphi^+, \varphi^-) = (\tilde{x}, \tilde{\varphi}^+, \tilde{\varphi}^-)$ with $\tilde{x}, \tilde{\varphi}^\pm \in (xR[[x]]  \oplus  \varphi^+ R[[x]] [\varphi^-] \oplus \varphi^- R[[x]] ) \subset R[[x]] [\varphi^+, \varphi^-]$ or in $x^{-1} R[[x^{-1}]] [\varphi^+, \varphi^-]$ since these will include the formal $N=2$ superconformal coordinates over a superalgebra $R$ vanishing at $(x,\varphi^+,\varphi^-) = 0$ and $(x,\varphi^+,\varphi^-) = (\infty,0,0)$, respectively.    

First in order to characterize formal $N=2$ superconformal functions corresponding to $N=2$ superconformal local coordinates invertible in a neighborhood of zero, and vanishing at zero, we note that such an $H(x, \varphi^+, \varphi^-) = (\tilde{x}, \tilde{\varphi}^+, \tilde{\varphi}^-)$ with $\tilde{x}, \tilde{\varphi}^\pm \in (xR[[x]] \oplus \varphi^+ R[[x]][\varphi^-] \oplus \varphi^- R[[x]])$ is uniquely determined by even formal power series $g^\pm(x)$ and odd formal power series $\psi^\pm(x)$ satisfying $\psi^\pm(0) = 0$ and $g^\pm(0) \in (R^0)^\times$, where $f$ is uniquely determined by (\ref{formal-superconformal-condition4}) since $f(0) = 0$, and where $(R^0)^\times$ denotes the invertible elements of $R^0$.  Thus an $N=2$ superconformal formal power series invertible in a neighborhood of zero and vanishing at zero is uniquely determined by
\begin{eqnarray}
g^\pm(x) &= &\sum_{j \in \mathbb{N}} a^\pm_j x^j , \qquad \quad \mbox{for $a^\pm_j \in R^0$ and $a^\pm_0 \in (R^0)^\times$} \label{expansion of g}\\
\psi^\pm(x) &= &\sum_{j \in \Z} m^\pm_{j- \frac{1}{2}} x^j , \qquad \mbox{for $m^\pm_{j - \frac{1}{2}} \in R^1$} , \label{expansion of psi}
\end{eqnarray}
and has the form (\ref{formal-superconformal-condition1}) -- (\ref{formal-superconformal-condition3}).  Explicitly, then we have that if $H(x,\varphi^+,\varphi^-) = (\tilde{x}, \tilde{\varphi}^+, \tilde{\varphi}^-)$ is formal $N=2$ superconformal and invertible in a neighborhood of zero and vanishing at zero, then $\tilde{x}, \tilde{\varphi}^\pm$ are of the form (\ref{powerseries1}) and (\ref{powerseries2}), respectively, with $(z,\theta^+, \theta^-) = (x, \varphi^+ , \varphi^-)$.

Let $H(x,\varphi^+,\varphi^-) = (\tilde{x}, \tilde{\varphi}^+, \tilde{\varphi}^-)$ be an $N=2$ superconformal invertible formal power series  vanishing at zero with the coefficient of the $\varphi^\pm$ terms in $\tilde{\varphi}^\pm$, respectively, equal to one, i.e, such that $H$ satisfies
\begin{equation}\label{T-condition}
\frac{\partial}{\partial \varphi^\pm} \tilde{\varphi}^\pm \Bigl|_{(x,\varphi^+, \varphi^-) = (0,0,0)}  = 1.  
\end{equation}
Then $H$ is uniquely determined by 
\begin{eqnarray}
\varphi^\mp \tilde{\varphi}^\pm  &=&  \varphi^\mp \left( \varphi^\pm +  \sum_{j \in \Z} \Bigl( m^\pm_{j-\frac{1}{2}} x^j + \varphi^\pm a^\pm_j x^j  \Bigr) \right) \\
&=& \varphi^\mp (\psi^\pm(x) + \varphi^\pm g^\pm(x)) \nonumber .
\end{eqnarray}
See Remark \ref{determining-superconformal-remark}.

We wish to express any such $H$ in terms of a formal exponential of an infinite sum of certain superderivations.  For any two even formal series $g^\pm(x) \in R^0[[x]]$ with constant term equal to one, and any two odd formal series $\psi^\pm (x) \in xR^1[[x]]$, we first express $\varphi^\mp (\psi^\pm(x) + g^\pm(x)) \in \varphi^\mp R[[x]][\varphi^\pm]$ in terms of $\varphi^\mp$ times the exponential of an infinite sum of superderivations in $\mbox{Der} (R((x))[\varphi^+,\varphi^-])$ acting on $\varphi^\pm$, respectively.  Then we prove that this exponential of superderivations acting on $(x, \varphi^+, \varphi^-)$ is in fact $N=2$ superconformal with the coefficient of the $\varphi^\pm$ term in $\tilde{\varphi}^\pm$ equal to one, respectively, and that there is a one-to-one correspondence between such exponential expressions and formal $N=2$ superconformal power series in $R[[x]][\varphi^+, \varphi^-]$ vanishing at $(x,\varphi^+, \varphi^-) = 0$ with coefficient of $\varphi^\pm$ in $\tilde{\varphi}^\pm$  equal to one.  

Let $\mathcal{A}^\pm_j$, for $j \in \Z$, be even formal variables, and let  $\mathcal{M}^\pm_{j - 1/2}$, for $j \in \Z$, be odd formal variables.  Let  $\mathcal{A}^\pm  = \{\mathcal{A}^\pm_j\}_{j \in \Z}$ and $\mathcal{M}^\pm = \{\mathcal{M}^\pm_{j - 1/2}\}_{j \in \Z}$, and consider the $\mathbb{Q}$-superalgebra $\mathbb{Q}[\mathcal{A}^+, \mathcal{A}^-,\mathcal{M}^+, \mathcal{M}^-]$ of 
polynomials in the formal variables $\mathcal{A}^\pm_1, \mathcal{A}^\pm_2,\dots$ and $\mathcal{M}^\pm_{1/2}, \mathcal{M}^\pm_{3/2},\dots$.  Consider the even superderivations 
\begin{eqnarray}
L_j(x,\varphi^+,\varphi^-) &=& - \biggl( \Lx \biggr) \label{L-notation}\\
J_j(x,\varphi^+,\varphi^-) &=& - \Jx  \label{J-notation}
\end{eqnarray}
and the odd superderivations
\begin{eqnarray}
G^\pm_{j -\frac{1}{2}} (x,\varphi^+,\varphi^-) = - \biggl( \Gx \biggr) \label{G-notation}
\end{eqnarray}
in $\mbox{Der} (R[[x,x^{-1}]] [\varphi^+,\varphi^-])$, for $j \in \mathbb{Z}$. We define the sequences 
\begin{eqnarray}
E^{0,\pm}(\mathcal{A}^+, \mathcal{A}^-, \mathcal{M}^+, \mathcal{M}^-) &=& \left\{E_j^\pm (\mathcal{A}^+, \mathcal{A}^-, \mathcal{M}^+, \mathcal{M}^-)\right\}_{j \in \Z} \\
E^{1,\pm} (\mathcal{A}^+, \mathcal{A}^-, \mathcal{M}^+,  \mathcal{M}^-) &=& \bigl\{E^\pm_{j - \frac{1}{2}}(\mathcal{A}^+, \mathcal{A}^-, \mathcal{M}^+, \mathcal{M}^-)\bigr\}_{j \in \Z}
\end{eqnarray}
of even and odd elements, respectively, in $\mathbb{Q}[\mathcal{A}^+, \mathcal{A}^-, \mathcal{M}^+,  \mathcal{M}^-]$ by  
\begin{multline}\label{exp1}
\varphi^\mp \biggl( \varphi^\pm + \sum_{j \in \Z} \left( E^\pm_{j-\frac{1}{2}} (\mathcal{A}^+, \mathcal{A}^-, \mathcal{M}^+,   \mathcal{M}^-)x^j  + \varphi^\pm E^\pm_j (\mathcal{A}^+, \mathcal{A}^-, \mathcal{M}^+,   \mathcal{M}^-)x^j \right) \biggr) \\
=  \varphi^\mp \exp\Biggl( \! - \! \sum_{j \in \Z} \Bigl( \mathcal{A}^+_j L_j(x,\varphi^+, \varphi^-) + \mathcal{A}^-_j J_j(x,\varphi^+, \varphi^-) \\
+ \mathcal{M}^+_{j - \frac{1}{2}} G^+_{j -\frac{1}{2}} (x,\varphi^+, \varphi^-) + \mathcal{M}^-_{j - \frac{1}{2}} G^-_{j -\frac{1}{2}} (x,\varphi^+, \varphi^-) \Bigr) \! \Biggr) \cdot \varphi^\pm .
\end{multline}
As usual, ``exp'' denotes the formal exponential series, when it is defined, as it is in the case of the above exponential of superderivations in 
\[\mbox{Der} \; (\mathbb{Q}[\mathcal{A}^+, \mathcal{A}^-, \mathcal{M}^+, \mathcal{M}^-][[x]][\varphi^+, \varphi^-]).\]   
The reason for the $\varphi^\mp$ multiplier in (\ref{exp1}), is that we are in fact uniquely defining two series  
\begin{multline*}
\psi^\pm (x) + \varphi^\pm g^\pm (x) \\
=  \varphi^\pm + \sum_{j \in \Z} \left( E^\pm_{j-\frac{1}{2}} (\mathcal{A}^+, \mathcal{A}^-, \mathcal{M}^+,   \mathcal{M}^-)x^j  + \varphi^\pm E^\pm_j (\mathcal{A}^+, \mathcal{A}^-, \mathcal{M}^+,   \mathcal{M}^-)x^j \right)  
\end{multline*}  
in $\mathbb{Q}[\mathcal{A}^+, \mathcal{A}^-, \mathcal{M}^+, \mathcal{M}^-][[x]][\varphi^\pm]$, respectively, by means of $\varphi^\mp$ times certain series in $\mathbb{Q}[\mathcal{A}^+, \mathcal{A}^-, \mathcal{M}^+, \mathcal{M}^-][[x]][\varphi^+,\varphi^-]$ (cf. Remark \ref{determining-superconformal-remark}).

Let $(R^0)^\infty$ be the set of all sequences $\{A_j\}_{j \in \Z}$ of even elements in $R$,  let $(R^1)^\infty$ be the set of all sequences $\{M_{j - 1/2}\}_{j 
\in \Z}$ of odd elements in $R$, and let $R^\infty = (R^0)^\infty \oplus (R^1)^\infty$.  Given any 
\begin{eqnarray*}
(A^+, A^-,M^+,M^-) &=& \bigl(\{A^+_j \}_{ j \in \Z},\{A^-_j \}_{ j \in \Z}, \{M^+_{j - \frac{1}{2}} \}_{j \in \Z}, \{M^-_{j - \frac{1}{2}} \}_{j \in \Z}\bigr) \\
&=& \bigl\{(A^+_j, A^-_j, M^+_{j - \frac{1}{2}},  M^-_{j - \frac{1}{2}}) \bigr\}_{j \in \Z} \in (R^\infty)^2 , 
\end{eqnarray*}
since the $E^\pm_j(\mathcal{A}^+, \mathcal{A}^-, \mathcal{M}^+,\mathcal{M}^-)$ are in $\mathbb{Q}[\mathcal{A}^+, \mathcal{A}^-, \mathcal{M}^+, \mathcal{M}^-]$, for $j \in \frac{1}{2} \Z$, we have a well-defined sequence 
\begin{multline*}
E(A^+,A^-,M^+, M^-) = (E^{0,+}(A^+,A^-,M^+,M^-),E^{0,-}(A^+,A^-,M^+,M^-), \\
E^{1,+}(A^+,A^-,M^+,M^-), E^{1,-}(A^+,A^-,M^+,M^-))
\end{multline*} 
in $(R^\infty)^2$ by substituting $A^+$, $A^-$, $M^+$ and $M^-$ into $E^{0,\pm}(\mathcal{A}^+, \mathcal{A}^-, \mathcal{M}^+, \mathcal{M}^-)$ and $E^{1,\pm}(\mathcal{A}^+, \mathcal{A}^-, \mathcal{M}^+, \mathcal{M}^-)$, respectively.  This defines a map    
\begin{eqnarray*}
E : (R^\infty)^2 &\longrightarrow& (R^\infty)^2 \\
(A^+, A^-,M^+, M^-) &\mapsto& E(A^+,A^-,M^+,M^-) .
\end{eqnarray*}

\begin{prop}\label{Ebijection}
The map $E$ is a bijection.  In particular, $E$ has an inverse
$E^{-1}$.
\end{prop}

\begin{proof} 
Define the following grading by weight on $\mathbb{Q} [x, \varphi^+, \varphi^-,\mathcal{A}^+, \mathcal{A}^-, \mathcal{M}^+, \mathcal{M}^- ]$:  For $j \in \Z$,
\begin{equation*}
\begin{array}{llllll}
\mathrm{wt} \;  x^j &=& j & \mathrm{wt} \; \varphi^\pm &=& \frac{1}{2} \\
\mathrm{wt} \; \mathcal{A}^\pm_j &=& -j & \mathrm{wt} \; \mathcal{M}^\pm_{j-\frac{1}{2}} &=& -j + \frac{1}{2}\\
\mathrm{wt} \; c &=& 0 \quad \mbox{for $c \in \mathbb{Q}$.} & 
\end{array}
\end{equation*}
Note that extending this grading of $\mathbb{Q}[x, \varphi^+, \varphi^-, \mathcal{A}^+, \mathcal{A}^-, \mathcal{M}^+, \mathcal{M}^- ]$ to derivations on $\mathbb{Q}[x, \varphi^+, \varphi^-, \mathcal{A}^+, \mathcal{A}^-, \mathcal{M}^+, \mathcal{M}^- ]$, we have that $\mathrm{wt} \; \frac{\partial}{\partial x} = -1$ and $\mathrm{wt} \; \frac{\partial}{\partial \varphi^\pm} = -\frac{1}{2}$.  
Thus letting 
\begin{multline}
\mathcal{T}= - \! \sum_{j \in \Z} \Bigl( \mathcal{A}^+_j L_j(x,\varphi^+, \varphi^-) + \mathcal{A}^-_j J_j(x,\varphi^+, \varphi^-) \\
+ \mathcal{M}^+_{j - \frac{1}{2}} G^+_{j -\frac{1}{2}} (x,\varphi^+, \varphi^-) + \mathcal{M}^-_{j - \frac{1}{2}} G^-_{j -\frac{1}{2}} (x,\varphi^+, \varphi^-) \Bigr)
\end{multline}
we have $\mathrm{wt} \; \mathcal{T} = 0$.  Therefore, $e^\mathcal{T} \cdot \varphi^\pm$ are both homogeneous of weight $1/2$ in $\mathbb{Q} [\mathcal{A}^+, \mathcal{A}^-, \mathcal{M}^+, \mathcal{M}^- ] [[x]][\varphi^+, \varphi^-]$.  

Note that 
\begin{eqnarray*}
\lefteqn{\varphi^\mp e^\mathcal{T} \cdot \varphi^\pm} \\
&=&\varphi^\mp \left(  \varphi^\pm + \mathcal{T}\varphi^\pm + \frac{1}{2} \mathcal{T}^2 \varphi^\pm + \frac{1}{3!} \mathcal{T}^3 \varphi^\pm + \cdots \right)  \\
&=& \varphi^\mp \biggl( \varphi^\pm + \sum_{j \in \Z} \biggl(  \Bigl(\frac{j + 1}{2}\Bigr) \mathcal{A}^+_j   \varphi^\pm x^j \pm \mathcal{A}^-_j \varphi^\pm x^j + \mathcal{M}^\pm_{j - \frac{1}{2}} x^j \\
& & \quad \pm j \mathcal{M}^\mp_{j - \frac{1}{2}} \varphi^+ \varphi^- x^{j-1}  \biggr)  + \frac{1}{2} \mathcal{T}^2 \varphi^\pm + \frac{1}{3!} \mathcal{T}^3 \varphi^\pm + \cdots \biggr) \\
&=& \varphi^\mp \biggl( \varphi^\pm + \sum_{j \in \Z} \biggl(  \Bigl(\frac{j + 1}{2}\Bigr) \mathcal{A}^+_j   \varphi^\pm x^j \pm \mathcal{A}^-_j \varphi^\pm x^j + \mathcal{M}^\pm_{j - \frac{1}{2}} x^j \biggr) + \frac{1}{2} \mathcal{T}^2 \varphi^\pm \\
& & \quad  + \frac{1}{3!} \mathcal{T}^3 \varphi^\pm + \cdots \biggr) \\
\end{eqnarray*}
where any term in $\mathcal{T}^k \varphi^\pm$ has order $k$ in the $A^\pm_j$'s and $M^\pm_{j - \frac{1}{2}}$'s.

Thus by the definition of $E^\pm_k (\mathcal{A}^+, \mathcal{A}^-, \mathcal{M}^+, \mathcal{M}^-)$, for $k \in \frac{1}{2} \Z$ given by (\ref{exp1}), we see that, for $j \in \Z$,
\begin{equation}\label{E-even}
E^\pm_j(\mathcal{A}^+, \mathcal{A}^-, \mathcal{M}^+, \mathcal{M}^-) = \frac{j+1}{2} \mathcal{A}^+_j \pm \mathcal{A}_j^- + r_j^\pm (\mathcal{A}^+, \mathcal{A}^-, \mathcal{M}^+, \mathcal{M}^-) ,
\end{equation}
and 
\begin{equation}\label{E-odd}
E^\pm_{j - \frac{1}{2}}(\mathcal{A}^+, \mathcal{A}^-, \mathcal{M}^+, \mathcal{M}^-) = \mathcal{M}^\pm_{j - \frac{1}{2}} + r_{j-\frac{1}{2}} ^\pm(\mathcal{A}^+, \mathcal{A}^-, \mathcal{M}^+, \mathcal{M}^-)   
\end{equation}
where $r_k^\pm (\mathcal{A}^+, \mathcal{A}^-, \mathcal{M}^+, \mathcal{M}^-) \in \mathbb{Q} [\mathcal{A}^+, \mathcal{A}^-, \mathcal{M}^+, \mathcal{M}^- ]$ is homogeneous of weight $-k$ for $k \in \frac{1}{2} \Z$ and only involves terms of degree at least two in the $\mathcal{A}^\pm_l$'s and $\mathcal{M}^\pm_{l-\frac{1}{2}}$'s for $l \in \Z$.  Therefore $r^\pm_k$ must involve only $\mathcal{A}^\pm_l$'s and $\mathcal{M}^\pm_{l-\frac{1}{2}}$'s of weight strictly greater than $k$.  
That is 
\begin{multline}\label{r1}
r^\pm_j (\mathcal{A}^+, \mathcal{A}^-, \mathcal{M}^+, \mathcal{M}^-) =  r^\pm_j (\mathcal{A}^+_1, \dots,  \mathcal{A}^+_{j-1}, \mathcal{A}^-_1, \dots, \mathcal{A}^-_{j-1}, \mathcal{M}^+_{\frac{1}{2}}, \dots, \\
\mathcal{M}^+_{j- \frac{1}{2}} , \mathcal{M}^-_{\frac{1}{2}}, \dots, \mathcal{M}^-_{j- \frac{1}{2}} ) 
\end{multline}
is in 
\[\mathbb{Q} [ \mathcal{A}^+_1, \dots,  \mathcal{A}^+_{j-1}, \mathcal{A}^-_1, \dots, \mathcal{A}^-_{j-1}, \mathcal{M}^+_{\frac{1}{2}}, \dots, \mathcal{M}^+_{j- \frac{1}{2}} , \mathcal{M}^-_{\frac{1}{2}}, \dots, \mathcal{M}^-_{j- \frac{1}{2}} ], \]
and
\begin{multline}\label{r2}
r^\pm_{j-\frac{1}{2}} (\mathcal{A}^+, \mathcal{A}^-, \mathcal{M}^+, \mathcal{M}^-) =  r^\pm_{j-\frac{1}{2}} (\mathcal{A}^+_1, \dots,  \mathcal{A}^+_{j-1}, \mathcal{A}^-_1, \dots, \mathcal{A}^-_{j-1}, \mathcal{M}^+_{\frac{1}{2}}, \\
\dots, \mathcal{M}^+_{j- \frac{3}{2}} , \mathcal{M}^-_{\frac{1}{2}}, \dots, \mathcal{M}^-_{j- \frac{3}{2}} ) 
\end{multline}
is in 
\[\mathbb{Q} [ \mathcal{A}^+_1, \dots,  \mathcal{A}^+_{j-1}, \mathcal{A}^-_1, \dots, \mathcal{A}^-_{j-1}, \mathcal{M}^+_{\frac{1}{2}}, \dots, \mathcal{M}^+_{j- \frac{3}{2}} , \mathcal{M}^-_{\frac{1}{2}}, \dots, \mathcal{M}^-_{j- \frac{3}{2}} ]. \]

Given $(a^+,a^-,m^+,m^-) \in (R^\infty)^2$ consider the infinite system of equations  
\begin{eqnarray}
E^\pm_j(A^+,A^-, M^+, M^-) &=& a^\pm_j \label{system-of-equations1} \\
E^\pm_{j - \frac{1}{2}}(A^+,A^-, M^+, M^-) &=& m^\pm_{j-\frac{1}{2}} \label{system-of-equations2}
\end{eqnarray}
for the unknown sequence $(A^+,A^-, M^+, M^-) \in (R^\infty)^2$.   This system of equations is equivalent to
\begin{eqnarray*}
\frac{1}{j+1} \left( a^+_j+ a^-_j \right) &=&  \frac{1}{j+1} \left( E^+_j(A^+,A^-, M^+, M^-) + E^-_j(A^+,A^-, M^+, M^-) \right)  \\
&=& A_j^+ +  \frac{1}{j+1}  \left( r^+_j(A^+,A^-, M^+, M^-) + r^-_j(A^+,A^-, M^+, M^-) \right)  \\
 \frac{1}{2} \left( a^+_j - a^-_j \right) &=& \frac{1}{2} \left( E^+_j(A^+,A^-, M^+, M^-) - E^-_j(A^+,A^-, M^+, M^-) \right)  \\
&=& A_j^- + \frac{1}{2} \left( r^+_j(A^+,A^-, M^+, M^-) - r^-_j(A^+,A^-, M^+, M^-) \right)  \\
 m^\pm_{j-\frac{1}{2}} &=& E^\pm_{j - \frac{1}{2}}(A^+,A^-, M^+, M^-) \\
 &=& M^\pm_{j-\frac{1}{2}} + r^\pm_{j-\frac{1}{2}} (A^+, A^-, M^+, M^-)
\end{eqnarray*}
i.e.,
\begin{eqnarray}
\frac{1}{j+1} \left( a^+_j+ a^-_j \right) &=& A_j^+ +  p^+_j(A^+,A^-, M^+, M^-)  \label{system-of-equations3}\\
\frac{1}{2} \left( a^+_j - a^-_j \right) &=& A_j^- + p^-_j(A^+,A^-, M^+, M^-) \\
m^\pm_{j-\frac{1}{2}} &=& M^\pm_{j-\frac{1}{2}} + r^\pm_{j-\frac{1}{2}} (A^+, A^-, M^+, M^-), \label{system-of-equations5}
\end{eqnarray}
where $p^\pm_j$ and $r^\pm_{j-\frac{1}{2}}$ only depend on $A_k^\pm$'s for $k<j$ and on 
$M^\pm_{k-\frac{1}{2}}$'s for $k=j$ and $k<j$, respectively.   Thus we can determine $(A^+, A^-, M^+, M^-)$ in terms of $(a^+, a^-, m^+, m^-)$ by determining  $M^\pm_\frac{1}{2}$, then $A^\pm_1$, then $M^\pm_\frac{3}{2}$, then $A^\pm_2$, etc.   That is, we determine the $A_j^\pm$'s and $M^\pm_{j-\frac{1}{2}}$ in order of weight.  It follows that the system of equations (\ref{system-of-equations3}) -- (\ref{system-of-equations5}) (and thus the system of equations (\ref{system-of-equations1}) and (\ref{system-of-equations2})) has a unique solution 
\begin{eqnarray}
\lefteqn{ (A^+,A^-,M^+,M^-) } \\
&= & E^{-1} (a^+,a^-,m^+,m^-) \nonumber \\
&=& \Bigl( (E^{-1})^{0,+} (a^+,a^-,m^+,m^-) , (E^{-1})^{0,-} (a^+,a^-,m^+,m^-), \nonumber \\
& & \quad (E^{-1})^{1,+} (a^+,a^-,m^+,m^-) , (E^{-1})^{1,-} (a^+,a^-,m^+,m^-) \Bigr). \nonumber
\end{eqnarray}
The proposition follows immediately. 
\end{proof}

\begin{cor}
For any two formal power series in $R^1[[x]] \oplus \varphi^\pm R^0[[x]]$, respectively, of the form 
\begin{equation}\label{gpsi} 
\psi^\pm(x) + \varphi^\pm g^\pm(x)  = \varphi^\pm +  \sum_{j \in \Z} \Bigl( m^\pm_{j-\frac{1}{2}} x^j + \varphi^\pm a^\pm_j x^j  \Bigr)  
\end{equation}
we have 
\begin{eqnarray}
\lefteqn{\varphi^\mp  ( \psi^\pm(x) + \varphi^\pm g^\pm(x)  )}  \\
&=& \varphi^\mp \exp  \Biggl(\! - \! \sum_{j \in \Z} \biggl( (E^{-1})^+_j(a^+,a^-,m^+,m^-) L_j(x,\varphi^+,\varphi^-) \nonumber \\
& & \quad + (E^{-1})^-_j(a^+,a^-,m^+,m^-) J_j(x,\varphi^+,\varphi^-) \nonumber \\
& & \quad + (E^{-1})^+_{j - \frac{1}{2}}(a^+,a^-,m^+,m^-) G^+_{j - \frac{1}{2}}(x, \varphi^+, \varphi^-) \nonumber \\
& & \quad + (E^{-1})^-_{j - \frac{1}{2}}(a^+,a^-,m^+,m^-) G^-_{j - \frac{1}{2}}(x, \varphi^+, \varphi^-)  \biggr) \! \Biggr) \! \cdot \varphi^\pm .  \nonumber
\end{eqnarray}
\end{cor}

\begin{proof} Using equation (\ref{exp1}) and the fact that $E$ is a bijection, we have
\begin{eqnarray*}
\lefteqn{\varphi^\mp  ( \psi^\pm(x) + \varphi^\pm g^\pm(x)  )  } \\
&=&   \varphi^\mp \biggl( \varphi^\pm +  \sum_{j \in \Z} \Bigl( m^\pm_{j-\frac{1}{2}} x^j + \varphi^\pm a^\pm_j x^j  \Bigr) \biggr) \\
&=&  \varphi^\mp \biggl( \varphi^\pm + \sum_{j \in \Z} E^\pm_{j-\frac{1}{2} } (E^{-1}(a^+,a^-,m^+,m^-)) x^j \\
& & \quad + \varphi^\pm E^\pm_j (E^{-1}(a^+,a^-,m^+,m^-)) x^j \biggr) 
\end{eqnarray*}
\begin{eqnarray*}
&=& \varphi^\mp \exp  \Biggl(\! - \! \sum_{j \in \Z} \biggl( (E^{-1})^+_j(a^+,a^-,m^+,m^-) L_j(x,\varphi^+,\varphi^-) \\
& & \quad + (E^{-1})^-_j(a^+,a^-,m^+,m^-) J_j(x,\varphi^+,\varphi^-)  \\
& & \quad + (E^{-1})^+_{j - \frac{1}{2}}(a^+,a^-,m^+,m^-) G^+_{j - \frac{1}{2}}(x, \varphi^+, \varphi^-) \\
& & \quad + (E^{-1})^-_{j - \frac{1}{2}}(a^+,a^-,m^+,m^-) G^-_{j - \frac{1}{2}}(x, \varphi^+, \varphi^-)  \biggr) \! \Biggr) \! \cdot \varphi^\pm  .  
\end{eqnarray*}
\end{proof}

\begin{prop}\label{superconformal}
Let $R$ be a superalgebra and 
\[(A^+,A^-,M^+,M^-) = \{(A_j^+, A^-_j, M^+_{j - 1/2}, M^-_{j - 1/2}) \}_{j \in \Z}  \in  (R^\infty)^2. \]
Then  
\begin{multline}\label{H}
H(x, \varphi^+,\varphi^-)  = \exp\Biggl( \! - \! \sum_{j \in \Z} \Bigl( A^+_j L_j(x,\varphi^+, \varphi^-) + A^-_j J_j(x,\varphi^+, \varphi^-) \\
+ M^+_{j - \frac{1}{2}} G^+_{j - \frac{1}{2}}(x,\varphi^+, \varphi^-) + M^-_{j - \frac{1}{2}} G^-_{j - \frac{1}{2}}(x,\varphi^+, \varphi^-) \Bigr) \Biggr) \cdot (x, \varphi^+, \varphi^-) 
\end{multline}
is $N=2$ superconformal and is the unique formal $N=2$ superconformal power series $H(x, \varphi^+,\varphi^-) = (\tilde{x}, \tilde{\varphi}^+, \tilde{\varphi}^-)$ vanishing at zero satisfying 
\begin{equation}\label{T-condition2}
\frac{\partial}{\partial \varphi^\pm} \tilde{\varphi}^\pm \Bigl|_{(x,\varphi^+, \varphi^-) = (0,0,0)}  = 1 
\end{equation}
such that
\begin{multline}
\varphi^\mp \tilde{\varphi}^\pm = \varphi^\mp \biggl(\varphi^\pm + \sum_{j \in \Z}
\Bigl( E^\pm_{j-\frac{1}{2}} (A^+,A^-,M^+,M^-) x^j \\
+ \varphi^\pm E^\pm_j (A^+,A^-,M^+,M^-) x^j\Bigr) \biggr). 
\end{multline}
\end{prop}

\begin{proof}  
Let 
\begin{multline}\label{define-T}
\mathcal{T} = - \sum_{j \in \Z} \left( A^+_j L_j(x,\varphi^+,\varphi^-) + A^-_j J_j(x,\varphi^+,\varphi^-)  \right. \\
\left.  + M^+_{j - \frac{1}{2}} G^+_{j - \frac{1}{2}}(x,\varphi^+,\varphi^-) + M^-_{j -\frac{1}{2}} G^-_{j - \frac{1}{2}}(x,\varphi^+,\varphi^-) \right) .
\end{multline} 
Then $\mathcal{T} \in (\mbox{Der}(R[[x]][\varphi^+, \varphi^-]))^0$, i.e., $\mathcal{T}$ is even, and thus 
\[e^\mathcal{T} \cdot (x, \varphi^+,\varphi^-) = (e^\mathcal{T} \cdot x, e^\mathcal{T} \cdot \varphi^+,
e^\mathcal{T} \cdot \varphi^-) \]
is in $(R[[x]][\varphi^+, \varphi^-])^0 \oplus (R[[x]][\varphi^+, \varphi^-])^1\oplus (R[[x]][\varphi^+, \varphi^-])^1$.  Let
\[h^\pm(x,\varphi^+,\varphi^-) = \sum_{j \in \Z} \biggl( \Bigl(A^+_j \Bigl(\frac{j + 1}{2} \Bigr) \pm A^-_j \Bigr) (x^j \mp j x^{j-1} \varphi^+\varphi^-) - 2 M^\pm_{j - \frac{1}{2}} j  x^{j - 1} \varphi^\mp \biggr) .\]  
Then $h^\pm \in (R[[x]][\varphi^+,\varphi^-])^0$, i.e., $h^\pm$ is even, and we have
\begin{equation}
 [D^\pm,\mathcal{T}] = h^\pm(x,\varphi^+, \varphi^-)D^\pm . 
 \end{equation} 
Thus 
\begin{equation}
D^\pm e^\mathcal{T} \cdot x = e^{(h^\pm + \mathcal{T})} \cdot D^\pm \cdot x = e^{(h^\pm + \mathcal{T})} \cdot \varphi^\mp . 
\end{equation}
By Proposition \ref{conformalproof},
\begin{equation}\label{thisequation} 
e^{y(h^\pm + \mathcal{T})} \cdot \varphi^\mp = e^{y(h^\pm + \mathcal{T})} \cdot (\varphi^\mp 1) = \left(e^{y\mathcal{T}} \cdot \varphi^\mp \right)\bigl(e^{y(h^\pm + \mathcal{T})} \cdot 1\bigr) .
\end{equation}
But in this case, the coefficient of $y^n$ for a fixed $n \in \mathbb{N}$ has terms with powers of $x$ greater than or equal to $n - 1$.  Thus we can set $y = 1$, and each power series in equation
(\ref{thisequation}) has only a finite number of $x^j$ terms for a given $j \in \mathbb{N}$, i.e., each term is a well-defined power series in $x$.  Therefore  
\begin{equation}
e^{(h^\pm + \mathcal{T})} \cdot \varphi^\mp = \left(e^\mathcal{T} \cdot \varphi^\mp \right)\bigl(e^{(h^\pm + \mathcal{T})} \cdot 1\bigr) .
\end{equation}
Thus writing $H(x, \varphi^+,\varphi^-) = (e^\mathcal{T} \cdot x, e^\mathcal{T} \cdot \varphi^+, e^\mathcal{T} \cdot \varphi^-) = (\tilde{x}, \tilde{\varphi}^+, \tilde{\varphi}^-)$, we have  
\begin{eqnarray}
D^\pm \tilde{x} &=& e^{(h^\pm + \mathcal{T})} \cdot \varphi^\mp \; = \; \left(e^\mathcal{T} \cdot \varphi^\mp \right) \bigl( e^{(h^\pm + \mathcal{T})} \cdot 1\bigr) \\
&=&  \left(e^\mathcal{T} \cdot \varphi^\mp \right)\bigl(e^{(h^\pm + \mathcal{T})} \cdot D^\pm \cdot \varphi^\pm \bigr) \; = \; \left(e^\mathcal{T} \cdot \varphi^\mp \right) \left(D^\pm e^\mathcal{T} \cdot \varphi^\pm \right) \nonumber \\
&=& \tilde{\varphi}^\mp D^\pm \tilde{\varphi}^\pm .  \nonumber
\end{eqnarray}  
In addition,
\begin{equation}
D^\pm \tilde{\varphi}^\mp =  D^\pm e^\mathcal{T} \cdot \varphi^\mp  =   e^{(h^\pm + \mathcal{T})} \cdot D^\pm \cdot \varphi^\mp  =   0 . 
\end{equation}
Thus $H$ satisfies the $N=2$ superconformal conditions (\ref{nice-superconformal-condition}).  

The uniqueness follows {}from Proposition \ref{Ebijection}, and the fact that $H$ is uniquely determined by $\varphi^\mp \tilde{\varphi}^\pm$ and the fact that $\tilde{x}(0,0,0) = 0$. 
\end{proof}

\begin{rema}\label{envelope-again} {\em 
That $\mathcal{T}$ as defined by (\ref{define-T}) is an even superderivation is due to the fact that $\mathcal{T}$ exists in the $R$-envelope of  $\mbox{Der}(\mathbb{C}[[x]][\varphi^+, \varphi^-])$; see Remark \ref{envelope}.    This fact is crucial for the development of the moduli space of $N=2$ super-Riemann spheres with tubes in terms of exponentials of infinitesimal $N=2$ superconformal transformations following  \cite{H-book} and \cite{B-memoir}.  It means that we are working in the envelope of a Lie superalgebra which, by definition, is an ordinary Lie algebra. }
\end{rema}

Now we would like to include in our study of formal $N=2$ superconformal power series vanishing at zero, those with the coefficients of $\varphi^\pm$ in $\tilde{\varphi}^\pm$, respectively, equal to an invertible even element of $R$ but not necessarily equal to one.

For $b \in (R^0)^\times$, we define the linear operators $b^{2x \frac{\partial}{\partial x}}$ and $b^{\varphi^\pm \frac{\partial}{\partial \varphi^\pm}}$ {}from $R[x, x^{-1}$, $\varphi^+, \varphi^-]$ to itself
by    
\begin{eqnarray}
b^{2x \frac{\partial}{\partial x}} \cdot c (\varphi^+)^{j^+} (\varphi^-)^{j^-} x^k &=& c
(\varphi^+)^{j^+} (\varphi^-)^{j^-} b^{2k} x^k\\
b^{\varphi^\pm \frac{\partial}{\partial \varphi^\pm}} \cdot c (\varphi^+)^{j^+} (\varphi^-)^{j^-} x^k &=& c b^{j^\pm} (\varphi^+)^{j^+} (\varphi^-)^{j^-} x^k
\end{eqnarray} 
for $c \in R$, $j^\pm \in \mathbb{Z}_2$, and $k \in \mathbb{Z}$.  Then, recalling the notation (\ref{L-notation}) and (\ref{J-notation}), the operators
\begin{eqnarray}
\qquad b^{-2L_0(x, \varphi^+, \varphi^-)} &=&  b^{\left(2x \frac{\partial}{\partial x} +  \left( \varphi^+ \frac{\partial}{\partial \varphi^+} +  \varphi^- \frac{\partial}{\partial \varphi^-} \right) \right)} \; = \; b^{2x \frac{\partial}{\partial x}} b^{\varphi^+ \frac{\partial}{\partial \varphi^+}}  b^{\varphi^- \frac{\partial}{\partial \varphi^-}}\\
b^{-J_0(x, \varphi^+, \varphi^-)} &=&  b^{ \left( \varphi^+ \frac{\partial}{\partial \varphi^+} -  \varphi^- \frac{\partial}{\partial \varphi^-}  \right)} \; = \;  b^{\varphi^+ \frac{\partial}{\partial \varphi^+}}  b^{-\varphi^- \frac{\partial}{\partial \varphi^-}}
\end{eqnarray} 
are well-defined linear operators on $R [x,x^{-1},\varphi^+, \varphi^-]$.   These operators can be extended to operators on $R[[x,x^{-1}]][\varphi^+, \varphi^-]$ in the obvious way.  We note
that for any formal series $H(x,\varphi^+, \varphi^-) \in R[[x,x^{-1}]][\varphi^+, \varphi^-]$, we have   
\begin{eqnarray}\label{a0property}
\lefteqn{a^{-2L_0(x, \varphi^+, \varphi^-)} \cdot b^{-J_0(x, \varphi^+, \varphi^-)} \cdot H(x, \varphi^+, \varphi^-) }\\
&=& H(a^{-2L_0(x, \varphi^+, \varphi^-)} \cdot b^{-J_0(x, \varphi^+, \varphi^-)} \cdot (x, \varphi^+ \varphi^-)) \nonumber \\ 
&=&  H(a^2x, ab\varphi^+, ab^{-1} \varphi^-).  \nonumber
\end{eqnarray}

If $H$ is of the form (\ref{H}), in order for $H(a^2x, ab\varphi^+, ab^{-1} \varphi^-)$ to correspond to an invertible local coordinate chart vanishing at zero, we must have $a \in (R^0)^\times$, i.e., $a$ must be an invertible even element of the underlying superalgebra $R$.

\begin{rema}\label{square-root-in-R} {\em 
One might wonder why we choose the slightly strange looking expressions $a^{2x \frac{\partial}{\partial x}}$ and $a^{-2L_0(x, \varphi^+, \varphi^-)}$ rather than $a^{x \frac{\partial}{\partial x}}$ and $a^{-L_0(x, \varphi^+, \varphi^-)}$, for $a \in (R^0)^\times$.  The reason is that  for $a \in (R^0)^\times$, we have $a^{-L_0(x, \varphi^+, \varphi^-)} \cdot  H(x,\varphi^+, \varphi^-)= H(a x,\sqrt{a} \varphi^+, \sqrt{a}^{-1} \varphi^-)$, where $(\sqrt{a})^2 = a$.  This forces one to define a square root on $(R^0)^\times$, which is equivalent to choosing a branch cut for the complex logarithm when $R = \bigwedge_*$.   However, it is more natural to give $\sqrt{a} = a'$, an invertible even element of $R$, as the basic data avoiding the need to keep track of a well-defined square root on $(R^0)^\times$.   With this in mind, we note that the operator $a^{-2L_0(x,\varphi^+, \varphi^-)}$ should not be read as $(a^2)^{-L_0(x,\varphi)}$ but as $a^{\left(2x \frac{\partial}{\partial x} +  \left( \varphi^+ \frac{\partial}{\partial \varphi^+} +  \varphi^- \frac{\partial}{\partial \varphi^-} \right) \right)}$ thus retaining the basic data $a$ rather than just $a^2$. }
\end{rema}

\begin{prop}\label{azero}
Let $a,b \in (R^0)^{\times}$.  Then 
\begin{equation}\label{L0}
a^{-2L_0(x,\varphi^+, \varphi^-)} \cdot b^{-J_0(x,\varphi^+, \varphi^-)} \cdot (x,\varphi^+, \varphi^-) = 
\left(a^2 x, a b \varphi^+ , ab^{-1} \varphi^- \right) 
\end{equation}
is $N=2$ superconformal.
\end{prop}
 
\begin{proof} 
Writing $\left(a^2 x, a b \varphi^+ , ab^{-1} \varphi^- \right) = (\tilde{x},\tilde{\varphi}^+, \tilde{\varphi}^-)$, we have that
\begin{eqnarray*}
D^- \tilde{\varphi}^+ = D^+ \tilde{\varphi}^- &=& 0\\
D^+\tilde{x} - \tilde{\varphi}^- D^+ \tilde{\varphi}^+ = \varphi^- a^2 - ab^{-1} \varphi^- ab &=& 0\\
D^-\tilde{x} - \tilde{\varphi}^+ D^-\tilde{\varphi}^- = \varphi^+ a^2 - ab \varphi^+ ab^{-1} &=& 0,
\end{eqnarray*}    
showing that $\left(a^2 x, a b \varphi^+ , ab^{-1} \varphi^- \right)$ satisfies the conditions (\ref{nice-superconformal-condition}).
\end{proof}

\begin{rema}\label{double-cover-remark}
{\em 
For $a,b \in (R^0)^{\times}$, since 
\begin{equation}
\left(a^2 x, a b \varphi^+ , ab^{-1} \varphi^- \right) = \left((-a)^2 x, (-a) (-b) \varphi^+ , (-a)(-b)^{-1} \varphi^- \right)
\end{equation} 
we have that 
\begin{multline}\label{double-cover}
a^{-2L_0(x,\varphi^+, \varphi^-)} \cdot b^{-J_0(x,\varphi^+, \varphi^-)} \cdot (x,\varphi^+, \varphi^-) \\
= (-a)^{-2L_0(x,\varphi^+, \varphi^-)} \cdot (-b)^{-J_0(x,\varphi^+, \varphi^-)} \cdot (x,\varphi^+, \varphi^-) .
\end{multline}
Thus the map $(a,b) \mapsto (a^2x, ab \varphi^+, ab^{-1} \varphi^-)$ {}from $(R^0)^\times \times (R^0)^\times$ to the set of invertible $N=2$ superconformal functions of the form $(a_1x, a_2\varphi^+, a_3\varphi^-)$, for $a_1,a_2,a_3 \in (R^0)^\times$, (i.e. with $a_2a_3 = a_1$) is  two-to-one. }
\end{rema}

For any $(A^+, A^-,M^+,M^-) \in (R^\infty)^2$, we define a map $\tilde{E}$ {}from $(R^\infty)^2$ to 
the set of all formal $N=2$ superconformal power series vanishing at zero and with the coefficient  of $\varphi^\pm$ in $\tilde{\varphi}^\pm$ equal to one, respectively, by defining
\begin{multline}\label{Etilde}
\varphi^\mp \tilde{E}^\pm(A^+, A^-,M^+,M^-)(x,\varphi^+, \varphi^-) = \varphi^\mp \biggl( \varphi^\pm  \\
+ \sum_{j \in \Z} \Bigl( E^\pm_{j-\frac{1}{2}} (A^+, A^-,M^+,M^-) x^j + \varphi^\pm E_j^\pm (A^+, A^-, M^+,M^-) x^j \Bigr) \biggr)
\end{multline}
and letting $\tilde{E}(A^+, A^-,M^+,M^-)(x, \varphi^+, \varphi^-) = (\tilde{x}, \tilde{\varphi}^+, \tilde{\varphi}^-)$ be the unique formal $N=2$ superconformal power series vanishing at zero with even coefficients of $\varphi^\pm$ equal to one such that $\tilde{E}^\pm (A^+, A^-,M^+,M^-)(x,\varphi^+, \varphi^-) = \tilde{\varphi}^\pm$ and such that (\ref{Etilde}) holds. 

For $a_0, b_0 \in (R^0)^\times$, we define a map $\hat{E}$ {}from $((R^0)^\times \times R^\infty)^2$ to the set of all formal $N=2$ superconformal power series vanishing at zero and with invertible leading even coefficients of $\varphi^\pm$ as follows.  Let 
$\tilde{E}(A^+,A^-,M^+, M^-)(x,\varphi^+, $ $\varphi^-) = (\tilde{x}, \tilde{\varphi}^+, \tilde{\varphi}^-)$.  Define 
\begin{equation}\label{define-Ehat}
\hat{E}(a_0,b_0,A^+,A^-,M^+, M^-)(x,\varphi^+,\varphi^-) = (a_0^2 \tilde{x}, a_0b_0\tilde{\varphi}^+, a_0b_0^{-1}\tilde{\varphi}^-).
\end{equation}  
Then $\hat{E}(a_0,b_0,A^+,A^-,M^+, M^-)(x,\varphi^+,\varphi^-)$ is the unique formal $N=2$ 
superconformal power series satisfying (\ref{define-Ehat}) with even coefficient of $\varphi^\pm$ equal to $a_0b_0^{\pm 1}$, respectively.   

\begin{rema}\label{E-hat-double}
{\em By Remark \ref{double-cover-remark}, $\hat{E}$ is not bijective; it is two-to-one.  If, however, we restrict the domain of $\hat{E}$ to $((R^0)^\times)^2/ \langle \pm 1\rangle 
\times (R^\infty)^2$, then $\hat{E}$ is bijective.}
\end{rema}

The following theorem is an immediate consequence of Propositions \ref{Ebijection}, \ref{superconformal}, \ref{azero}, and Remark \ref{E-hat-double}.    

\begin{thm}\label{above}
The map $\hat{E}$ {}from $((R^0)^\times)^2/ \langle \pm 1\rangle \times (R^\infty)^2$ to the set of all formal $N=2$ superconformal power series vanishing at zero with invertible even coefficients of $\varphi^\pm$  is a bijection.

The map $\tilde{E}$ {}from $(R^{\infty})^2$ to the set of formal $N=2$ superconformal power series vanishing at zero with even coefficients of $\varphi^\pm$ equal to one is also a bijection. 

In particular, we have inverses $\tilde{E}^{-1}$ and $\hat{E}^{-1}$.

Specifically, for $(a_0,b_0,A^+,A^-,M^+, M^-) \in ((R^0)^\times)^2/ \langle \pm 1\rangle \times (R^\infty)^2$, we have
\begin{eqnarray}
\lefteqn{\hat{E} (a_0,b_0,A^+,A^-,M^+, M^-) (x, \varphi^+, \varphi^-)  }\\
&=&  \exp \Biggl( - \sum_{j \in \Z} \biggl( A^+_j L_j(x,\varphi^+,\varphi^-) + A^-_j J_j(x,\varphi^+,\varphi^-)   \nonumber \\
& & \quad  + M^+_{j - \frac{1}{2}} G^+_{j - \frac{1}{2}}(x,\varphi^+,\varphi^-) + M^-_{j -\frac{1}{2}} G^-_{j - \frac{1}{2}}(x,\varphi^+,\varphi^-) \biggr) \Biggr)  \cdot \nonumber \\
& & \quad  \cdot a_0^{-2L_0(x,\varphi^+, \varphi^-)}  \cdot b_0^{-J_0(x,\varphi^+, \varphi^-)} \cdot (x, \varphi^+, \varphi^-) \nonumber
\end{eqnarray}
\begin{eqnarray}
\lefteqn{\tilde{E} (A^+,A^-,M^+, M^-) (x, \varphi^+, \varphi^-)  }\\
&=& \exp \Biggl( - \sum_{j \in \Z} \biggl( A^+_j L_j(x,\varphi^+,\varphi^-) + A^-_j J_j(x,\varphi^+,\varphi^-)   \nonumber \\
& & \quad + M^+_{j - \frac{1}{2}} G^+_{j - \frac{1}{2}}(x,\varphi^+,\varphi^-) + M^-_{j -\frac{1}{2}} G^-_{j - \frac{1}{2}}(x,\varphi^+,\varphi^-) \biggr) \Biggr)  \cdot  \nonumber \\
& & \quad \cdot (x, \varphi^+, \varphi^-) \nonumber . 
\end{eqnarray}
\end{thm}

We will use the notation
\begin{eqnarray}
\qquad \mathcal{T}_H (x, \varphi^+, \varphi^-) &=& - \sum_{j \in \Z} \biggl((E^{-1})^+_j (a^+, a^-,m^+,m^-) 
L_j(x,\varphi^+, \varphi^-)  \label{T-notation}\\
& & \qquad + (E^{-1})^-_j (a^+, a^-,m^+,m^-) J_j(x,\varphi^+, \varphi^-)\nonumber \\
& & \qquad + (E^{-1})^+_{j - \frac{1}{2}} (a^+, a^-,m^+, m^-) G^+_{j - \frac{1}{2}} (x,\varphi^+, \varphi^-) \nonumber \\
& & \qquad + (E^{-1})^-_{j - \frac{1}{2}} (a^+, a^-,m^+, m^-) G^-_{j - \frac{1}{2}}(x,\varphi^+, \varphi^-) \biggr) \nonumber
\end{eqnarray} 
for $H(x, \varphi^+, \varphi^-) = (\tilde{x}, \tilde{\varphi}^+, \tilde{\varphi}^-)$ vanishing at zero and formally $N=2$ superconformal satisfying
\begin{equation}\label{Ehat} 
\varphi^\mp \tilde{\varphi}^\pm = \varphi^\mp  a_0b_0^{\pm 1} \Bigl( \varphi^\pm + \sum_{j \in \Z} \left(
m^\pm_{j-\frac{1}{2}} x^j + \varphi^\pm a^\pm_j x^j  \right) \Bigr)
\end{equation} 
for some $(a_0,b_0,a^+, a^-,m^+, m^-) \in ((R^0)^{\times})^2/  \langle \pm 1\rangle \times (R^\infty)^2$.  Thus any such $N=2$ superconformal power series $H(x,\varphi^+, \varphi^-)$ can be written uniquely as  
\begin{equation}\label{Hexpansion} 
H(x, \varphi^+, \varphi^-) = e^{\mathcal{T}_H (x, \varphi^+, \varphi^-)} \cdot a_0^{-2L_0(x,\varphi^+, \varphi^-)} 
\cdot b_0^{-J_0(x,\varphi^+, \varphi^-)} \cdot (x, \varphi^+, \varphi^-) . 
\end{equation}

Recalling (\ref{power-series-at-ith-puncture}), we know that a local $N=2$ superconformal coordinate map $H(w, \rho^+, \rho^-) = (\tilde{w}, \tilde{\rho}^+, \tilde{\rho}^-)$ vanishing at $0 \in \bigwedge_\infty$ is completely determined by $\rho^\mp \tilde{\rho}^\pm = \rho^\mp (\psi^\pm (w) + \rho^\pm g^\pm (w))$, i.e., is completely determined by $g^\pm(w)$ and $\psi^\pm(w)$ where $g^\pm (w)$ can be expanded 
in a power series of the form $a_0^\pm( w + \sum_{j \in \Z} a^\pm_j  w^{j+1})$  with $a_0^\pm \in (\bigwedge_\infty^0)^{\times}$ and  $a^\pm_{j - 1/2} \in \bigwedge_\infty^0$, and $\psi^\pm(w)$ can be expanded in a power series of the form $ \sum_{j \in \Z}  m^\pm_{j-1/2} w^j$ with $m^\pm_{j-1/2} \in \bigwedge_\infty^1$, such that these power series are absolutely convergent to $g^\pm(w)$ and $\psi^\pm(w)$, respectively, in some neighborhood of zero.  Thus we see that formal $N=2$ superconformal power series $H(x,\varphi^+, \varphi^-)$ of the form (\ref{Hexpansion}) can be thought of as the ``local formal $N=2$ superconformal coordinate maps vanishing at zero'' or the ``local formal $N=2$ superconformal coordinate transformations fixing the coordinates of a fixed point to be zero".  {}From 
(\ref{Hexpansion}), we see that the ``local formal $N=2$ superconformal  transformations invertible and vanishing at zero'' are generated  uniquely (except for single-valuedness) by the ``infinitesimal formal $N=2$ superconformal transformations'' of the form   
\begin{multline} 
(\log a_0) 2L_0(x,\varphi^+, \varphi^-) + (\log b_0) J_0(x,\varphi^+, \varphi^-)  + \sum_{j \in \Z} \biggl(
A^+_j L_j(x,\varphi^+, \varphi^-) \\
+ A^-_j J_j(x,\varphi^+, \varphi^-) + M^+_{j - \frac{1}{2}} G^+_{j - \frac{1}{2}}(x,\varphi^+, \varphi^-) + M^-_{j - \frac{1}{2}} G^-_{j - \frac{1}{2}}(x,\varphi^+, \varphi^-) \biggr) ,
\end{multline} 
for $(a_0, b_0) \in ((R^0)^\times)^2/\langle \pm 1 \rangle$, and $\{(A^+_j,A^-_j,M^+_{j - 1/2}, M^-_{j - 1/2})\}_{j \in \Z} \in (R^\infty)^2$.  Theorem \ref{above} states that these ``infinitesimal $N=2$ superconformal transformations'' can be identified with elements in $((R^0)^\times)^2/ \langle \pm 1\rangle \times (R^\infty)^2$.

\begin{prop}\label{auto2}
Let $u, v \in R((x))[\varphi^+, \varphi^-]$; let $(A^+, A^-,M^+, M^-) \in (R^\infty)^2$; and let
\begin{multline}
\mathcal{T} = \mathcal{T}(x, \varphi^+, \varphi^-) = - \sum_{j \in \Z} \biggl(A_j^+ L_j(x,\varphi^+, \varphi^-) + A_j^- J_j(x,\varphi^+, \varphi^-) \label{T} \\
+ M_{j- \frac{1}{2}}^+ G^+_{j - \frac{1}{2}} (x,\varphi^+, \varphi^-) + M^-_{j - \frac{1}{2}} G^-_{j - \frac{1}{2}}(x,\varphi^+, \varphi^-) \biggr) .
\end{multline} 
Then
\begin{equation}\label{specialized-automorphism-property}
e^\mathcal{T} \cdot (uv) = \left( e^\mathcal{T} \cdot u \right)\left( e^\mathcal{T} \cdot v \right).
\end{equation}
\end{prop}

In other words, for $\mathcal{T}$, $u$, and $v$ given above Proposition \ref{conformalproof} holds if $h = 0$ and $y$ is set equal to 1.

\begin{proof}  
Define a $\frac{1}{2} \mathbb{Z}$-grading by weight on $R((x))[\varphi^+, \varphi^-]$ given by $\mathrm{wt} \ a = 0$ for $a \in R$, $\mathrm{wt}\ x^k = k$, for $k \in \mathbb{Z}$, and $\mathrm{wt}\ \varphi^\pm = \frac{1}{2}$.  Then $\mathcal{T} \in \mathrm{Der} (R((x))[\varphi^+, \varphi^-])$ given by (\ref{T}) acting on an element $u' \in R((x))[\varphi^+, \varphi^-]$ of homogeneous weight, raises the weight of $u'$ by $\frac{1}{2}$ and higher.  Thus if $\mathrm{wt} \ u' = k \in \frac{1}{2} \mathbb{Z}$, then $\frac{y^n \mathcal{T}^n}{n!} \cdot u'$ involves terms of weight greater than or equal to $k + \frac{n}{2}$, for $n \in \Z$.   Therefore, for $u \in R((x))[\varphi^+, \varphi^-]$, any term in $e^{y\mathcal{T}} \cdot u$ of homogeneous weight in $x$, $\varphi^+$, and $\varphi^-$ is a polynomial in $y$, and thus we can set $y= 1$.  That is $e^\mathcal{T} \cdot u$ has only a finite number of terms of weight $k \in \frac{1}{2} \mathbb{Z}$, and thus is a well-defined power series in $R((x))[\varphi^+, \varphi^-]$, and similarly for $e^\mathcal{T} \cdot v$.  Therefore (\ref{specialized-automorphism-property}) holds. 
\end{proof} 

\begin{prop}\label{Switch}
Let $\overline{H}(x, \varphi^+, \varphi^-) \in R((x))[\varphi^+, \varphi^-]$, and let $H$ be a formal $N=2$ superconformal power series vanishing at zero given by
\[H(x, \varphi^+, \varphi^-) = e^{\mathcal{T}_H (x, \varphi^+, \varphi^-)} \cdot (x, \varphi^+, \varphi^-).\]  
Then   
\begin{eqnarray}\label{switch}
\overline{H}(H(x, \varphi^+, \varphi^-)) &=& \overline{H}(e^{\mathcal{T}_H (x, \varphi^+, \varphi^-)} \cdot (x, \varphi^+, \varphi^-)) \\
&=& e^{\mathcal{T}_H (x, \varphi^+, \varphi^-)} \cdot \overline{H}(x, \varphi^+, \varphi^-) .  \nonumber
\end{eqnarray}
More generally, if $\overline{H}(x, \varphi^+, \varphi^-) \in R[[x, x^{-1}]][\varphi^+, \varphi^-]$ and the composition $\overline{H} \circ H$ is well defined, then (\ref{switch}) holds.
\end{prop}

\begin{proof}  Write $H(x, \varphi^+, \varphi^-) = e^{\mathcal{T}_H} \cdot (x, \varphi^+, \varphi^-) = (e^{\mathcal{T}_H} \cdot x, e^{\mathcal{T}_H} \cdot \varphi^+,e^{\mathcal{T}_H} \cdot \varphi^- ) = (\tilde{x},\tilde{\varphi}^+, \tilde {\varphi}^-)$. Equation (\ref{switch}) is trivial for $\overline{H}(x, \varphi^+, \varphi^-) = 1$, $\overline{H}(x, \varphi^+, \varphi^-) = x$, and $\overline{H}(x, \varphi^+, \varphi^-) = \varphi^\pm$.  

(i) We prove the result for $\overline{H}(x, \varphi^+, \varphi^-) = x^n$, with $n \in \mathbb{N}$ and $n > 1$, by induction on $n$.  Assume $e^{\mathcal{T}_H} \cdot x^k = (e^{\mathcal{T}_H} \cdot x)^k$ for $k \in \mathbb{N}$, $k < n$. Let $\overline{H}(x, \varphi^+, \varphi^-) = x^n$.  Then by Proposition \ref{auto2}, 
\begin{eqnarray*}
e^{\mathcal{T}_H} \cdot \overline{H}(x, \varphi^+, \varphi^-) \! \! &=& \! \! e^{\mathcal{T}_H} \cdot (x x^{n - 1}) = \left( e^{\mathcal{T}_H} \cdot x \right) \left( e^{\mathcal{T}_H} \cdot x^{n - 1} \right) \\ 
&=& \! \! \left( e^{\mathcal{T}_H} \cdot x \right) \left( e^{\mathcal{T}_H} \cdot x \right)^{n - 1} = \tilde{x} \tilde{x}^{n - 1} = \tilde{x}^n = \overline{H}(H(x, \varphi^+, \varphi^-)) . 
\end{eqnarray*}

(ii) We next prove the result for $\overline{H}(x, \varphi^+, \varphi^-) = x^{-n}$, $n \in \Z$.  Again by Proposition \ref{auto2}, 
\[ 1 = e^{\mathcal{T}_H} \cdot (x x^{-1}) = (e^{\mathcal{T}_H} \cdot x)(e^{\mathcal{T}_H} \cdot x^{-1}) .\]
Thus \[ e^{\mathcal{T}_H} \cdot x^{-1} = (e^{\mathcal{T}_H} \cdot x)^{-1} .\]

Assume $e^{\mathcal{T}_H} \cdot x^{-k} = (e^{\mathcal{T}_H} \cdot x)^{-k}$ for $k \in \mathbb{N}$, $k < n$. Let $\overline{H}(x, \varphi^+, \varphi^-) = x^{-n}$.  Then by Proposition \ref{auto2},   
\begin{eqnarray*}
e^{\mathcal{T}_H} \cdot \overline{H}(x, \varphi^+, \varphi^-) \! &=& \! e^{\mathcal{T}_H} \cdot (x^{-1} x^{-(n - 1)} ) = \left( e^{\mathcal{T}_H} \cdot x^{-1} \right) \left( e^{\mathcal{T}_H} \cdot x^{-(n - 1)} \right) \\ 
&=& \! \left( e^{\mathcal{T}_H} \cdot x \right)^{-1} \left( e^{\mathcal{T}_H} \cdot x \right)^{-(n - 1)} = \tilde{x}^{-1} \tilde{x}^{-(n - 1)} = \tilde{x}^{-n} \\
&=&\! \overline{H}(H(x, \varphi^+, \varphi^-)) . 
\end{eqnarray*}
Thus the result is true for $\overline{H}(x, \varphi^+, \varphi^-) = x^n$, $n \in \mathbb{Z}$.

(iii) For $\overline{H}(x, \varphi^+, \varphi^-) = \varphi^\pm x^n$, $n \in \mathbb{Z}$, we note that by Proposition \ref{auto2}, and the above cases  
\begin{eqnarray*}
e^{\mathcal{T}_H} \cdot \overline{H}(x, \varphi^+, \varphi^-) &=& e^{\mathcal{T}_H} \cdot (\varphi^\pm x^n) =
\left( e^{\mathcal{T}_H} \cdot \varphi^\pm \right) \left( e^{\mathcal{T}_H} \cdot x^n \right) \\ 
&=& \left( e^{\mathcal{T}_H} \cdot \varphi^\pm \right) \left( e^{\mathcal{T}_H} \cdot x \right)^n =
\tilde{\varphi}^\pm  \tilde{x}^n = \overline{H}(H(x, \varphi^+, \varphi^-)) . 
\end{eqnarray*} 

(iv) Finally for $\overline{H}(x, \varphi^+, \varphi^-) = \varphi^+ \varphi^- x^n$, $n \in \mathbb{Z}$, we note that by Proposition \ref{auto2}, and the above cases  
\begin{eqnarray*}
e^{\mathcal{T}_H} \cdot \overline{H}(x, \varphi^+, \varphi^-) &=& e^{\mathcal{T}_H} \cdot (\varphi^+ \varphi^- x^n) = \left( e^{\mathcal{T}_H} \cdot \varphi^+ \right) \left( e^{\mathcal{T}_H} \cdot \varphi^- \right) \left( e^{\mathcal{T}_H} \cdot x^n \right) \\ 
&=& \tilde{\varphi}^+ \tilde{\varphi}^-  \tilde{x}^n = \overline{H}(H(x, \varphi^+, \varphi^-)) . 
\end{eqnarray*} 

Since $e^{\mathcal{T}_H} \cdot  (\varphi^+)^{j^+}  (\varphi^-)^{j^-}  x^n \in R((x))[\varphi^+, \varphi^-]$, for $j^\pm = 0, 1$ and $n \in \mathbb{Z}$, the result follows by linearity. 
\end{proof} 

\begin{prop}\label{inverses}
Any formal $N=2$ superconformal power series $H(x, \varphi^+, \varphi^-)$ vanishing at zero and with even coefficient of $\varphi^\pm$ invertible has a unique inverse $H^{-1}(x, \varphi^+, \varphi^-)$ with respect to composition of formal power series, and this inverse is also $N=2$ superconformal vanishing at zero and with even coefficient of $\varphi^\pm$ invertible.  

If $\hat{E}^{-1} (H(x, \varphi^+, \varphi^-)) = (a_0, b_0, A^+, A^-, M^+, M^-)$, i.e.,
\begin{eqnarray}
\lefteqn{\; \; H(x,\varphi^+, \varphi^-) } \label{explicit-H-at-zero}\\
&=& \exp \Biggl(   - \sum_{j \in \Z} \biggl( A^+_j L_j(x,\varphi^+, \varphi^-) + A^-_j J_j(x,\varphi^+, \varphi^-) \nonumber \\
& & \quad + M^+_{j - \frac{1}{2}} G^+_{j - \frac{1}{2}}(x,\varphi^+, \varphi^-) + M^-_{j - \frac{1}{2}} G^-_{j - \frac{1}{2}}(x,\varphi^+, \varphi^-) \biggr)  \! \Biggr)\nonumber  \\
& & \quad \cdot  a_0^{-2L_0(x,\varphi^+, \varphi^-)}  \cdot  b_0^{-J_0(x,\varphi^+, \varphi^-)} \cdot (x,\varphi^+, \varphi^-) \nonumber \\
&=& \exp \left( \mathcal{T}_H (x,\varphi^+, \varphi^-) \right)  \cdot  a_0^{-2L_0(x,\varphi^+, \varphi^-)}  \cdot  b_0^{-J_0(x,\varphi^+, \varphi^-)} \cdot (x,\varphi^+, \varphi^-) \nonumber 
\end{eqnarray}
then 
\begin{eqnarray}
\lefteqn{\quad H^{-1}(x, \varphi^+, \varphi^-)} \label{inverse}\\
&=&  \exp \Biggl(  \sum_{j \in \Z} \biggl( a_0^{-2j} A^+_j L_j(x,\varphi^+, \varphi^-) + a_0^{-2j} A^-_j J_j(x,\varphi^+, \varphi^-) \nonumber \\
& & \quad + a_0^{-2j +1} b_0 M^+_{j - \frac{1}{2}} G^+_{j - \frac{1}{2}}(x,\varphi^+, \varphi^-) \nonumber \\
& & \quad + a_0^{-2j +1} b_0^{-1} M^-_{j - \frac{1}{2}} G^-_{j - \frac{1}{2}}(x,\varphi^+, \varphi^-) \biggr) 
\! \! \Biggr) \! \cdot \nonumber \\ 
& &  \quad  \cdot a_0^{ 2L_0(x,\varphi^+, \varphi^-)} \cdot b_0^{ J_0(x,\varphi^+, \varphi^-)} \cdot (x,\varphi^+, \varphi^-) 
\nonumber \\    
&=& a_0^{ 2L_0(x,\varphi^+, \varphi^-)} \cdot b_0^{ J_0(x,\varphi^+, \varphi^-)} \cdot \exp  \Biggl( \sum_{j \in \Z} 
\biggl( A^+_j L_j(x,\varphi^+, \varphi^-) \nonumber \\
& & \quad + A^-_j J_j(x,\varphi^+, \varphi^-)  + M^+_{j - \frac{1}{2}} G^+_{j - \frac{1}{2}}(x,\varphi^+, \varphi^-) \nonumber\\
& & \quad + M^-_{j - \frac{1}{2}} G^-_{j - \frac{1}{2}}(x,\varphi^+, \varphi^-)  \biggr) \! \! \Biggr) \! \cdot  (x, \varphi^+, \varphi^-) \nonumber \\
&=&  a_0^{ 2L_0(x,\varphi^+, \varphi^-)} \! \cdot b_0^{ J_0(x,\varphi^+, \varphi^-)} \! \cdot  \exp \left( -\mathcal{T}_H (x, \varphi^+, \varphi^-) \right)  \cdot  (x, \varphi^+, \varphi^-). \nonumber
\end{eqnarray}  
\end{prop}

\begin{proof} Using formula (\ref{inverse}) for $H^{-1}$, {}from Proposition \ref{Switch} and (\ref{a0property}) we have
\begin{eqnarray*}
\lefteqn{H(H^{-1}(x, \varphi^+, \varphi^-))} \\
&=& \left. H(H^{-1}(a_0^2 x_1, a_0 b_0 \varphi^+_1, a_0 b_0^{-1} \varphi^-_1))
\right|_{(x_1, \varphi_1^+, \varphi_1^-) = ( a_0^{-2} x, a_0^{-1} b_0^{-1} \varphi^+, a_0^{-1} b_0 \varphi^-)} \\
&=&  H  \Biggl( \exp \Biggl( \sum_{j \in \Z} \biggl( A_j^+ L_j(x_1,\varphi_1^+, \varphi_1^-)  + A_j^-J_j(x_1,\varphi_1^+, \varphi_1^-)  \\
& & \quad + M^+_{j - \frac{1}{2}} G^+_{j - \frac{1}{2}}(x_1,\varphi_1^+, \varphi_1^-) + M^-_{j -
\frac{1}{2}} G^-_{j - \frac{1}{2}}(x_1,\varphi_1^+, \varphi_1^-)  \biggr) \! \Biggr) \Bigg. \cdot \\
& & \quad \Biggl. \left. (x_1, \varphi_1^+, \varphi_1^-) \Biggr)  \right|_{(x_1, \varphi_1^+, \varphi_1^-) = ( a_0^{-2} x, a_0^{-1} b_0^{-1} \varphi^+, a_0^{-1} b_0 \varphi^-)}\\ 
&=& H(\exp\left( - \mathcal{T}_H (x_1, \varphi_1^+, \varphi_1^-) \right) \cdot  \\
& & \quad \left. (x_1, \varphi_1^+, \varphi_1^-) ) \right|_{(x_1, \varphi_1^+, \varphi_1^-) = ( a_0^{-2} x, a_0^{-1} b_0^{-1} \varphi^+, a_0^{-1} b_0 \varphi^-)} \\
&=&   \exp\left( - \mathcal{T}_H (x_1, \varphi_1^+, \varphi_1^-) \right) \cdot \\
& & \quad \left. H(x_1, \varphi_1^+, \varphi_1^-) \right|_{(x_1, \varphi_1^+, \varphi_1^-) = ( a_0^{-2} x, a_0^{-1} b_0^{-1} \varphi^+, a_0^{-1} b_0 \varphi^-)} \\  
&=& \exp\left( - \mathcal{T}_H (x_1, \varphi_1^+, \varphi_1^-) \right) \cdot \exp\left(\mathcal{T}_H (x_1, \varphi_1^+, \varphi_1^-) \right) \cdot a_0^{-2L_0(x_1,\varphi_1^+, \varphi_1^-)} \cdot\\
& &\left.  \quad  b_0^{-J_0(x_1,\varphi_1^+, \varphi_1^-)}  \cdot (x_1, \varphi_1) ) \right|_{(x_1, \varphi_1^+, \varphi_1^-) = ( a_0^{-2} x, a_0^{-1} b_0^{-1} \varphi^+, a_0^{-1} b_0 \varphi^-)}\\  
&=& \left. (a_0^2 x_1, a_0 b_0 \varphi_1^+, a_0 b_0^{-1} \varphi_1^-)  \right|_{(x_1, \varphi_1^+, \varphi_1^-) = ( a_0^{-2} x, a_0^{-1} b_0^{-1} \varphi^+, a_0^{-1} b_0 \varphi^-)}\\
&=&  (x, \varphi^+, \varphi^-) . 
\end{eqnarray*}

Similarly, {}from Proposition \ref{Switch}, and (\ref{a0property}) we have 
\begin{eqnarray*}
\lefteqn{H^{-1}(H(x, \varphi^+, \varphi^-))} \\
&=&  \left. H^{-1} (H(a_0^{-2} x_1, a_0^{-1}b_0^{-1} \varphi_1^+, a_0^{-1} b_0 \varphi_1^-)) \right|_{(x_1, \varphi_1^+, \varphi_1^-) = (a_0^2 x, a_0 b_0 \varphi^+, a_0 b_0^{-1} \varphi^-)} \\ 
&=& H^{-1} \Bigl( a_0^{ 2L_0(x_1,\varphi_1^+, \varphi_1^-)}  \cdot b_0^{ J_0(x_1,\varphi_1^+, \varphi_1^-)} \cdot \exp \left( \mathcal{T}_H (x_1, \varphi_1^+, \varphi_1^-) \right) \cdot\\
& & \quad  a_0^{ -2L_0(x_1,\varphi_1^+, \varphi_1^-)} \cdot b_0^{ -J_0(x_1,\varphi_1^+, \varphi_1^-)}  \cdot \\
& & \quad \left.  (x_1, \varphi_1^+, \varphi_1^-) \Bigr)  \right|_{(x_1, \varphi_1^+, \varphi_1^-) = (a_0^2 x, a_0 b_0 \varphi^+, a_0 b_0^{-1} \varphi^-)}  \\ 
&=& a_0^{ 2L_0(x_1,\varphi_1^+, \varphi_1^-)}  \cdot b_0^{ J_0(x_1,\varphi_1^+, \varphi_1^-)} \cdot \exp \left( \mathcal{T}_H (x_1, \varphi_1^+, \varphi_1^-) \right) \cdot a_0^{ -2L_0(x_1,\varphi_1^+, \varphi_1^-)} \cdot \\
& & \quad b_0^{ -J_0(x_1,\varphi_1^+, \varphi_1^-)}  \cdot \left. H^{-1} (x_1, \varphi_1^+, \varphi_1^-)  \right|_{(x_1, \varphi_1^+, \varphi_1^-) = (a_0^2 x, a_0 b_0 \varphi^+, a_0 b_0^{-1} \varphi^-)}  \\  
&=& a_0^{ 2L_0(x_1,\varphi_1^+, \varphi_1^-)}  \cdot b_0^{ J_0(x_1,\varphi_1^+, \varphi_1^-)} \cdot \exp \left( \mathcal{T}_H (x_1, \varphi_1^+, \varphi_1^-) \right) \cdot a_0^{ -2L_0(x_1,\varphi_1^+, \varphi_1^-)} \cdot \\
& & \quad b_0^{ -J_0(x_1,\varphi_1^+, \varphi_1^-)}  \cdot  b_0^{ J_0(x_1,\varphi_1^+, \varphi_1^-)}  \cdot a_0^{ 2L_0(x_1,\varphi_1^+, \varphi_1^-)} \cdot \\
& &  \quad \left.   \exp \left(- \mathcal{T}_H (x_1, \varphi_1^+, \varphi_1^-) \right) \cdot (x_1, \varphi_1^+, \varphi_1^-) )  \right|_{(x_1, \varphi_1^+, \varphi_1^-) = (a_0^2 x, a_0 b_0 \varphi^+, a_0 b_0^{-1} \varphi^-)}  \\  
&=& \left. a_0^{ 2L_0(x_1,\varphi_1^+, \varphi_1^-)}  \cdot b_0^{ J_0(x_1,\varphi_1^+, \varphi_1^-)}  \cdot (x_1,\varphi_1^+, \varphi_1^-) \right|_{(x_1, \varphi_1^+, \varphi_1^-) = (a_0^2 x, a_0 b_0 \varphi^+, a_0 b_0^{-1} \varphi^-)} \\ 
&=& (x, \varphi^+, \varphi^-) . 
\end{eqnarray*}

Since the formal composition of two formal $N=2$ superconformal power series is again $N=2$ superconformal, by Propositions \ref{superconformal} and \ref{azero} and the fact that $H^{-1}$ is given by (\ref{inverse}), the formal power series $H^{-1}(x, \varphi^+, \varphi^-)$ is $N=2$ superconformal. 
\end{proof} 

\begin{rema}\label{first-group-remark}
{\em {}From the proposition above, we see that the set of all formal $N=2$ superconformal power series vanishing at zero with even coefficient of $\varphi^\pm$ invertible is a group with composition as the group operation.  This is the group of ``formal $N=2$ superconformal local coordinate transformations vanishing at zero''.  We also see that, the set of all formal $N=2$ superconformal power series 
vanishing at zero with even coefficient of $\varphi^\pm$ equal to one is a subgroup.  This is the group of ``formal $N=2$ superconformal local coordinate transformations vanishing at zero with leading even coefficient of $\varphi^\pm$ equal to one''.}
\end{rema}

Let $H(x,\varphi^+, \varphi^-)$ and $\overline{H}(x, \varphi^+, \varphi^-)$ be two formal $N=2$ superconformal power series vanishing at zero and with even coefficients of $\varphi^\pm$ equal to one, i.e., let $H$ and $\bar{H}$ be of the form (\ref{H}) such that  
\begin{eqnarray}
\tilde{E}^{-1}(H(x, \varphi^+, \varphi^-)) &=& (A^+, A^-,M^+,M^-) \label{for-composition1}\\
\tilde{E}^{-1}(\overline{H}(x, \varphi^+,\varphi^-)) &=& (B^+,B^-,N^+,N^-) \label{for-composition2}
\end{eqnarray}  
for some $(A^+,A^-,M^+,M^-), (B^+,B^-,N^+,N^-) \in (R^\infty)^2$.  We define the {\it composition}  of $(A^+,A^-,M^+,M^-)$ and $(B^+,B^-,N^+,N^-)$, denoted $(A^+,A^-,M^+,M^-) \circ (B^+,B^-,N^+,N^-)$, by
\begin{equation}\label{sequencecompositiondef}
(A^+,A^-,M^+,M^-) \circ (B^+,B^-,N^+,N^-) = \tilde{E}^{-1}((\overline{H} \circ H) (x, \varphi^+,\varphi^-)) 
\end{equation}
where $(\overline{H} \circ H) (x, \varphi^+,\varphi^-)$ is the formal composition of $H$ and $\overline{H}$.

\begin{prop}\label{sequencecompositionprop} 
The set $(R^\infty)^2$ is a group with the operation $\circ$. Let  
\[(A^+,A^-,M^+,M^-) = \{ (A^+_j,A^-_j , M^+_{j - 1/2},M^-_{j - 1/2}) \}_{j \in \Z} \in (R^\infty)^2,\] 
and for any $t \in R^0$, define 
\[t(A^+,A^-,M^+,M^-) = \{ (t A^+_j,tA^-_j, t M^+_{j - 1/2}, tM^-_{j-1/2}) \}_{j \in \Z}.\]  
Then for $t_1,t_2 \in R^0$,   
\begin{equation}\label{composition}
(t_1(A^+,A^-,M^+,M^-)) \circ (t_2(A^+,A^-,M^+,M^-)) = (t_1 + t_2)(A^+,A^-,M^+,M^-) .
\end{equation} 
That is the map $t \mapsto t(A^+,A^-,M^+,M^-)$ is a homomorphism {}from the additive group $R^0$ to $(R^\infty)^2$. 
\end{prop}
\begin{proof}  Since the set of all formal $N=2$ superconformal power series of the form (\ref{H}) is a group with composition as its group operation, it is obvious {}from the definition of $\circ$ that
$R^{\infty}$ is a group with this operation.  Let $H_t (x,\varphi^+, \varphi^-) = \tilde{E} (t(A^+, A^-,M^+, M^-))$ for $t \in R^0$.  By Proposition \ref{Switch} 
\begin{eqnarray*}
H_{t_2} (H_{t_1} (x,\varphi^+, \varphi^-)) &=& e^{\mathcal{T}_{H_{t_1}} (x,\varphi^+, \varphi^-)} \cdot H_{t_2} (x,\varphi^+, \varphi^-) \\ 
&=& e^{\mathcal{T}_{H_{t_1}} (x,\varphi^+, \varphi^-)} \cdot e^{\mathcal{T}_{H_{t_2}} (x,\varphi^+, \varphi^-)} \cdot (x, \varphi^+, \varphi^-) \\
&=& e^{t_1 \mathcal{T}_{H_1} (x,\varphi^+, \varphi^-)} \cdot e^{t_2 \mathcal{T}_{H_1} (x,\varphi^+, \varphi^-)} \cdot (x, \varphi^+, \varphi^-) \\
&=& e^{(t_1 + t_2) \mathcal{T}_{H_1} (x,\varphi^+,\varphi^-)} \cdot (x,\varphi^+, \varphi^-) \\
&=& H_{(t_1 + t_2)} (x,\varphi^+, \varphi^-) \\
&=& \tilde{E} ((t_1 + t_2) (A^+, A^-,M^+, M^-)).
\end{eqnarray*}
Or equivalently,
\[ \tilde{E}^{-1} (H_{t_2} (H_{t_1} (x,\varphi^+, \varphi^-)) = (t_1 + t_2) (A^+, A^-,M^+, M^-) .\]
But then {}from the definition of $t_1 (A^+, A^-,M^+, M^-) \circ t_2 (A^+, A^-,M^+, M^-)$,
\[ \tilde{E}^{-1} (H_{t_2} (H_{t_1} (x,\varphi^+, \varphi^-)) = t_1 (A^+,A^-,M^+, M^-) \circ t_2 (A^+, A^-,M^+,M^-) .\]
Thus we obtain equation (\ref{composition}).  
\end{proof}

We now extend this composition $\circ$ to $((R^0)^\times)^2/\langle \pm 1 \rangle \times (R^\infty)^2$.  Again let $H$ and $\overline{H}$ be two formal $N=2$ superconformal power series of the form (\ref{H}) such that (\ref{for-composition1}) and (\ref{for-composition2}) hold, and let
\begin{eqnarray*}
H_{a^+_0, a^-_0} (x,\varphi^+,\varphi^-) &=& \hat{E} (a_0^+,a_0^-, A^+, A^-, M^+, M^-) (x,\varphi^+, \varphi^-) \\
&=& H((a_0^+)^2x ,a_0^+a_0^- \varphi^+, a_0^+ (a_0^-)^{-1} \varphi^-)\\
\overline{H}_{b_0^+, b_0^-} (x,\varphi^+, \varphi^-) &=& \hat{E}(b_0^+, b_0^-, B^+, B^-,N^+, N^-) (x,\varphi^+, \varphi^-) \\ 
&=& \overline{H}((b_0^+)^2 x,b_0^+ b_0^-\varphi^+, b_0^+ (b_0^-)^{-1} \varphi^-)
\end{eqnarray*}
for $(a_0^+, a_0^-, A^+, A^-, M^+, M^-), (b_0^+, b_0^-, B^+, B^-,N^+,N^-) \in( (R^0)^\times)^2/\langle \pm 1 \rangle  \times (R^\infty)^2$.  
Define
\begin{multline}\label{formal-composition-for-R}
(a_0^+, a_0^-, A^+, A^-, M^+, M^-) \circ (b_0^+, b_0^-, B^+, B^-,N^+,N^-)  \\
= \hat{E}^{-1}((\overline{H}_{b_0^+,b_0^-}\circ H_{a_0^+,a_0^-}) (x,\varphi^+, \varphi^-)) . 
\end{multline}
Then in terms of the composition defined on $(R^\infty)^2$ by  (\ref{sequencecompositiondef}), we have
\begin{multline}\label{composition-with-a}
(a_0^+, a_0^-, A^+, A^-, M^+, M^-) \circ (b_0^+, b_0^-, B^+, B^-,N^+,N^-) \\
=  \Bigl(a_0^+b_0^+, a_0^-b_0^-, (A^+, A^-, M^+, M^-) \circ 
\bigl\{(a_0^+)^{2j} B^+_j, (a_0^+)^{2j} B^-_j, \\
(a_0^+)^{2j - 1} (a_0^-)^{-1}  N^+_{j-\frac{1}{2}}, (a_0^+)^{2j - 1} a_0^- N^-_{j-\frac{1}{2}} \bigr\}_{j \in \Z}\Bigr) .
\end{multline}

\begin{rema}\label{group-remark} {\em 
With the composition operation defined above, $((R^0)^\times)^2/\langle \pm 1 \rangle  \times (R^\infty)^2$ is a group naturally isomorphic to the group of all formal $N=2$ superconformal power series of the form (\ref{explicit-H-at-zero}), i.e., isomorphic to the group of formal $N=2$ superconformal local coordinate transformations vanishing at zero.  See Remark \ref{first-group-remark}.  The subset $(R^\infty)^2$ is a subgroup of $((R^0)^\times)^2/\langle \pm 1 \rangle  \times (R^\infty)^2$ isomorphic to the group of all formal $N=2$ superconformal power series of the form (\ref{H}), i.e., isomorphic to the group of formal $N=2$ superconformal local coordinate transformations vanishing at zero with leading even coefficient of $\varphi^\pm$ equal to one.  In addition, 
\[\{(a_0^+, 1, A^+, A^-, M^+, M^-) \; | \; a_0^+ \in (R^0)^\times, \;   (A^+, A^-, M^+, M^-) \in (R^\infty)^2 \}\] and 
\[\{(1,a_0^-, A^+, A^-, M^+, M^-) \; | \; a_0^- \in (R^0)^\times, \; (A^+, A^-, M^+, M^-) \in (R^\infty)^2 \}\] 
are subgroups of $((R^0)^\times)^2/\langle \pm 1 \rangle \times (R^\infty)^2$.  In Section \ref{NS-algebra-section} we continue our discussion of subgroups of $((R^0)^\times)^2/\langle \pm 1 \rangle \times (R^\infty)^2$, and in Section \ref{moduli-group-section} we relate these group structures to subsets of the moduli space of $N=2$ super-Riemann spheres with one incoming tube and one outgoing tube and a corresponding ``sewing" operation.   }
\end{rema}

We now want to consider the ``formal $N=2$ superconformal coordinate maps vanishing at infinity.''  Let $H(x, \varphi^+, \varphi^-) = (\tilde{x}, \tilde{\varphi}^+, \tilde{\varphi}^-)$ be a formal $N=2$ superconformal powers series with $\tilde{x}, \tilde{\varphi}^\pm  \in x^{-1} R[[x^{-1}]] [\varphi^+, \varphi^-]$  such that
\begin{eqnarray}\label{atinfty}
\varphi^\mp \tilde{\varphi}^\pm &=& \varphi^\mp \biggl( \frac{i\varphi^\pm}{x} + \sum_{j \in \Z}
\left(   m^\pm_{j - \frac{1}{2}} x^{-j} + i \varphi^\pm a^\pm_j x^{-j - 1} \right) \biggr) \\
&=& \varphi^\mp(\psi^\pm(x) + \varphi^\pm g^\pm(x)). \nonumber
\end{eqnarray}
for $a_j^\pm \in R^0$ and $m_{j - 1/2}^\pm \in R^1$, for $j \in \Z$.  That is, $H$ is $N=2$ superconformal with leading coefficient of $\varphi^\pm x^{-1}$ in $\tilde{\varphi}^\pm$ equal to $i$.  Then  if $\lim_{x \rightarrow \infty} H(x, 0, 0) = 0$, $H$ is uniquely determined by (\ref{atinfty}).  That is $H$ is uniquely determined by (\ref{atinfty}) (or equivalently $g^\pm(x)$ and $\psi^\pm(x)$) and the requirement that it be $N=2$ superconformal and vanishing at infinity (cf. Remark \ref{determining-superconformal-remark}).  Explicitly, we have that if $H(x,\varphi^+,\varphi^-) = (\tilde{x}, \tilde{\varphi}^+, \tilde{\varphi}^-)$ is formal $N=2$ superconformal, vanishing at infinity and with leading even coefficient of $\varphi^\pm x^{-1}$ equal to $i$, then $\tilde{x}, \tilde{\varphi}^\pm$ are of the form (\ref{powerseries1-infinity}) -- (\ref{powerseries2-infinity}), respectively, with $(z,\theta^+, \theta^-) = (x, \varphi^+ , \varphi^-)$, and with $a_0^\pm = 1$.
 
Define 
\begin{equation}
I (x, \varphi^+, \varphi^-) = \Bigl(\frac{1}{x}, \frac{i \varphi^+}{x}, \frac{i \varphi^-}{x} \Bigr),
\end{equation}
cf. (\ref{origin-of-I}).   Then $I$ is $N=2$ superconformal, vanishing at infinity and with leading even coefficient of $\varphi^\pm x^{-1}$ equal to $i$.  Note that $I^{-1} =  (1/x, - i \varphi^+/x,$ $- i \varphi^+/x)$ is $N=2$ superconformal, vanishing at infinity and with even coefficient of $\varphi^\pm x^{-1}$ equal to $-i$.  

We now want to use the results we have developed about formal $N=2$ superconformal series vanishing at zero to express any formal $N=2$ superconformal series vanishing at infinity and with even 
coefficient of $\varphi^\pm x^{-1}$ equal to $i$ in terms of superderivations in $\mbox{Der}(R((x^{-1}))[\varphi^+, \varphi^-])$.

Let $H(x, \varphi^+, \varphi^-)$ be $N=2$ superconformal vanishing at infinity and with leading even coefficient of $\varphi^\pm x^{-1}$ equal to $i$, and let 
\begin{equation}\label{define-H-minus-one}
H_{-1}(x, \varphi^+, \varphi^-) = H \circ I^{-1} (x, \varphi^+, \varphi^-).
\end{equation}   
Then $H_{-1}$ is $N=2$ superconformal vanishing at zero and with leading even coefficient of $\varphi^\pm$ equal to one.  Thus $H_{-1}$ is of the form (\ref{Hexpansion}) with $a_0 = b_0 = 1$ and by Proposition \ref{inverses} has a well-defined compositional inverse $H_{-1}^{-1}(x,\varphi^+, \varphi^-)$.  

Note that $H \circ I^{-1} \circ H_{-1}^{-1}(x, \varphi^+, \varphi^-)$, and $I^{-1} \circ H_{-1}^{-1} \circ H(x, \varphi^+, \varphi^-)$ are well-defined formal $N=2$ superconformal series with components in $R[[x]][\varphi^+, \varphi^-]$ and $xR[[x^{-1}]][\varphi^+, \varphi^-]$, respectively.  Moreover, it is clear that the compositional inverse of $H$ is $H^{-1}(x,\varphi^+, \varphi^-) = I^{-1} \circ H_{-1}^{-1}(x, \varphi^+, \varphi^-)$ which has components in $x^{-1} R[[x]][\varphi^+, \varphi^-]$.

Recall the even and odd superderivations introduced in (\ref{L-notation}) -- (\ref{G-notation}).

\begin{prop}\label{Infinity} 
Given $H(x, \varphi^+, \varphi^-)$ $N=2$ superconformal vanishing at infinity and with leading coefficient of $\varphi^\pm x^{-1}$ equal to $i$, we have
\begin{eqnarray}
H(x, \varphi^+, \varphi^-) \! \! \!  &=& \! \! \exp \Biggl( \sum_{j \in \Z} \biggl( A^+_j L_{-j}(x,\varphi^+, \varphi^-)   - A^-_j J_{-j}(x,\varphi^+, \varphi^-)  \biggr. \Biggr.    \label{infty}\\
& &  \quad +  iM^+_{j - \frac{1}{2}} G^+_{-j +\frac{1}{2}}(x,\varphi^+, \varphi^- )  \nonumber \\
& & \quad \biggl. \Biggl. +  i M^-_{j - \frac{1}{2}}G^-_{-j +\frac{1}{2}}(x,\varphi^+, \varphi^- ) \biggr) \! \Biggr) \cdot  \Bigl(\frac{1}{x}, \frac{i \varphi^+}{x},\frac{i \varphi^-}{x}  \Bigr) \nonumber \\
&=& \! \! \exp \biggl( \mathcal{T}_{H_{-1}} \Bigl(\frac{1}{x}, \frac{i \varphi^+}{x}, \frac{i \varphi^-}{x} \Bigr) \biggr) \cdot \ \Bigl(\frac{1}{x}, \frac{i \varphi^+}{x}, \frac{i \varphi^-}{x} \Bigr) . \nonumber
\end{eqnarray} 
for some $(A^+, A^-, M^+, M^-) \in (R^\infty)^2$, and for $H_{-1} = H \circ I^{-1}$. The inverse of $H$ with respect to composition is given by  $H^{-1}(x, \varphi^+, \varphi^-) = I^{-1} \circ H^{-1}_{-1} (x, \varphi^+, \varphi^-)$, and   
\begin{eqnarray}
 H^{-1} \circ I (x, \varphi^+, \varphi^-)  &=& H^{-1} \Bigl(\frac{1}{x}, \frac{i \varphi^+}{x}, \frac{i \varphi^-}{x} \Bigr)  \label{forsewing} \\
&=&  \exp \biggl(- \mathcal{T}_{H_{-1}} \Bigl(\frac{1}{x}, \frac{i \varphi^+}{x}, \frac{i \varphi^-}{x}\Bigr) \biggr) \cdot (x, \varphi^+, \varphi^-) .\nonumber
\end{eqnarray} 
\end{prop}

\begin{proof} Since $H$ is $N=2$ superconformal vanishing at infinity and with leading coefficient of $\varphi^\pm x^{-1}$ equal to $i$, the power series 
\[H_{-1} (x, \varphi^+, \varphi^-) = H \circ I^{-1} (x,\varphi^+, \varphi^-)\] 
has components in $R[[x]][\varphi^+, \varphi^-]$, is vanishing at zero and has leading coefficients of $\varphi^\pm$ equal to one.   Thus by Proposition \ref{superconformal}, we have
\begin{eqnarray*}
\lefteqn{H_{-1} (x, \varphi^+, \varphi^-) }\\
&=&  \exp\Biggl( \! - \! \sum_{j \in \Z} \Bigl( A^+_j L_j(x,\varphi^+, \varphi^-) + A^-_j J_j(x,\varphi^+, \varphi^-) \\
& & \quad + M^+_{j - \frac{1}{2}} G^+_{j - \frac{1}{2}}(x,\varphi^+, \varphi^-) + M^-_{j - \frac{1}{2}} G^-_{j - \frac{1}{2}}(x,\varphi^+, \varphi^-) \Bigr) \Biggr) \cdot (x, \varphi^+, \varphi^-) \\
&=&  \exp ( \mathcal{T}_{H_{-1}} (x, \varphi^+, \varphi^-)) \cdot (x, \varphi^+, \varphi^-) 
\end{eqnarray*}
for some $(A^+, A^-, M^+, M^-) \in (R^\infty)^2$. 
Write $I(x, \varphi^+, \varphi^-) = (1/x, i \varphi^+/x,i \varphi^-/x) = (\tilde{x}, \tilde{\varphi}^+, \tilde{\varphi}^-)$.   By the chain rule
\[ ix^{-1} \frac{\partial}{\partial \tilde{\varphi}^\pm} = \frac{\partial \tilde{\varphi}^\pm}{\partial \varphi^\pm} \frac{\partial}{\partial \tilde{\varphi}^\pm} = \frac{\partial}{\partial \varphi^\pm} -  \frac{\partial
\tilde{x}}{\partial \varphi^\pm} \frac{\partial}{\partial \tilde{x}}  -   \frac{\partial \tilde{\varphi}^\mp}{\partial \varphi^\pm} \frac{\partial}{\partial \tilde{\varphi}^\mp} = \frac{\partial}{\partial \varphi^\pm} \]
and thus
\begin{eqnarray*}
-x^{-2} \frac{\partial}{\partial \tilde{x}} &=& \frac{\partial \tilde{x}}{\partial x} \frac{\partial}{\partial \tilde{x}} =
\frac{\partial}{\partial x} - \Bigl(- \frac{\partial
\tilde{\varphi}^+}{\partial x} \frac{\partial}{\partial \tilde{\varphi}^+} \Bigr) - \Bigl(- \frac{\partial
\tilde{\varphi}^-}{\partial x} \frac{\partial}{\partial \tilde{\varphi}^-} \Bigr)\\
&=& \frac{\partial}{\partial x} - i \varphi^+ x^{-2} \frac{\partial}{\partial \tilde{\varphi}^+}  - i \varphi^- x^{-2}  \frac{\partial}{\partial \tilde{\varphi}^-}\\
&=& \frac{\partial}{\partial x} + \varphi^+ x^{-1} \frac{\partial}{\partial \varphi^+} + \varphi^- x^{-1} \frac{\partial}{\partial \varphi^-} .
\end{eqnarray*}
Therefore 
\begin{eqnarray*}
\lefteqn{ H(x, \varphi^+, \varphi^-) = H \circ I^{-1} \circ I (x, \varphi^+, \varphi^-) = H \circ I^{-1} (\tilde{x}, \tilde{\varphi}^+, \tilde{\varphi}^-) = H_{-1} (\tilde{x}, \tilde{\varphi}^+, \tilde{\varphi}^-) }  \\
&=& \! \! \exp ( \mathcal{T}_{H_{-1}} (\tilde{x}, \tilde{\varphi}^+, \tilde{\varphi}^-)) \cdot (\tilde{x},
\tilde{\varphi}^+, \tilde{\varphi}^-) \\ 
&=& \! \! \exp\Biggl( \! - \! \sum_{j \in \Z} \Bigl( A^+_j L_j(\tilde{x},\tilde{\varphi}^+, \tilde{\varphi}^-) + A^-_j J_j(\tilde{x},\tilde{\varphi}^+, \tilde{\varphi}^-) \\
& & \ \  + M^+_{j - \frac{1}{2}} G^+_{j - \frac{1}{2}}(\tilde{x},\tilde{\varphi}^+, \tilde{\varphi}^-) + M^-_{j - \frac{1}{2}} G^-_{j - \frac{1}{2}}(\tilde{x},\tilde{\varphi}^+, \tilde{\varphi}^-) \Bigr)\!  \Biggr) \cdot (\tilde{x},\tilde{\varphi}^+, \tilde{\varphi}^-) \hspace{.5in}\\
&=& \! \!  \exp\Biggl( \sum_{j \in \Z} \Bigl( A^+_j \biggl(\tilde{x}^{j + 1} \frac{\partial}{\partial \tilde{x}} + \Bigl(\frac{j + 1}{2}\Bigr) \tilde{x}^j \Bigl( \tilde{\varphi}^+ \frac{\partial}{\partial \tilde{\varphi}^+} + \tilde{\varphi}^- \frac{\partial}{\partial \tilde{\varphi}^-}\Bigr) \! \biggr) \\
& & \ \  + A^-_j \tilde{x}^j\Bigl(\tilde{\varphi}^+\frac{\partial}{\partial \tilde{\varphi}^+}  - \tilde{\varphi}^- \frac{\partial}{\partial \tilde{\varphi}^-}\Bigr)  + M^+_{j - \frac{1}{2}} \biggl( \tilde{x}^j \Bigl( \frac{\partial}{\partial \tilde{\varphi}^+} - \tilde{\varphi}^- \frac{\partial}{\partial \tilde{x}}\Bigr) \\
& & \ \  + j\tilde{x}^{j-1} \tilde{\varphi}^+ \tilde{\varphi}^- \frac{\partial}{\partial \tilde{\varphi}^+} \biggr) + M^-_{j - \frac{1}{2}} \biggl(\tilde{x}^j \Bigl( \frac{\partial}{\partial \tilde{\varphi}^-} - \tilde{\varphi}^+ \frac{\partial}{\partial \tilde{x}}\Bigr) -  j\tilde{x}^{j-1} \tilde{\varphi}^+ \tilde{\varphi}^- \frac{\partial}{\partial \tilde{\varphi}^-} \Bigr) \! \biggr)\!  \Biggr) \cdot \\
& & \ \ (\tilde{x},\tilde{\varphi}^+, \tilde{\varphi}^-) \\
\end{eqnarray*}
\begin{eqnarray*}
&=& \!\! \exp\Biggl( \sum_{j \in \Z} \Bigl( A^+_j \biggl(x^{-j -1} (-x^2) \Bigl(  \frac{\partial}{\partial x} + \varphi^+ x^{-1} \frac{\partial}{\partial \varphi^+} + \varphi^- x^{-1} \frac{\partial}{\partial \varphi^-} \Bigr) \\
& & \ \  + \Bigl(\frac{j + 1}{2}\Bigr) x^{-j} \Bigl( i\varphi^+ x^{-1} (-ix) \frac{\partial}{\partial \varphi^+} + i \varphi^- x^{-1} (-ix) \frac{\partial}{\partial \varphi^-}\Bigr) \! \biggr) \\
& & \ \  + A^-_j x^{-j} \Bigl(i \varphi^+ x^{-1} (-ix)\frac{\partial}{\partial \varphi^+} - i \varphi^-x^{-1}(-ix) \frac{\partial}{\partial \varphi^-}\Bigr)  + M^+_{j - \frac{1}{2}} \biggl(x^{-j} \Bigl(-ix  \frac{\partial}{\partial \varphi^+} \\
& & \ \   - i\varphi^- x^{-1} (-x^2) \Bigl(  \frac{\partial}{\partial x} + \varphi^+ x^{-1} \frac{\partial}{\partial \varphi^+} + \varphi^- x^{-1} \frac{\partial}{\partial \varphi^-} \Bigr) \Bigr) \\
& & \ \ - j x^{-j-1} \varphi^+ \varphi^- (-ix) \frac{\partial}{\partial \varphi^+} \biggr)  + M^-_{j - \frac{1}{2}} \biggl(x^{-j} \Bigl(-ix \frac{\partial}{\partial \varphi^-} - i \varphi^+ x^{-1} (-x^2) \Bigl(  \frac{\partial}{\partial x} \\
& & \ \ + \varphi^+ x^{-1} \frac{\partial}{\partial \varphi^+} + \varphi^- x^{-1} \frac{\partial}{\partial \varphi^-} \Bigr) \Bigr) + j x^{-j-1} \varphi^+ \varphi^-(-ix) \frac{\partial}{\partial \varphi^-} \Bigr)\!  \biggr) \! \Biggr) \cdot \\
& & \ \ \left(\frac{1}{x},\frac{i\varphi^+}{x}, \frac{i\varphi^-}{x} \right) \\
&=& \! \! \exp\Biggl( - \sum_{j \in \Z} \Bigl( A^+_j \biggl(x^{-j +1}    \frac{\partial}{\partial x} +  \Bigl(\frac{-j + 1}{2}\Bigr) x^{-j} \Bigl( \varphi^+ \frac{\partial}{\partial \varphi^+} +  \varphi^-  \frac{\partial}{\partial \varphi^-}\Bigr) \! \biggr) \\
& & \ \  - A^-_j x^{-j} \Bigl(\varphi^+ \frac{\partial}{\partial \varphi^+} - \varphi^- \frac{\partial}{\partial \varphi^-}\Bigr)  + i M^+_{j - \frac{1}{2}} \biggl(x^{-j+1} \Bigl(  \frac{\partial}{\partial \varphi^+} - \varphi^-   \frac{\partial}{\partial x} \Bigr)  \\
& & \ \  + (-j+1) x^{-j} \varphi^+ \varphi^- \frac{\partial}{\partial \varphi^+} \biggr)  +i M^-_{j - \frac{1}{2}} \biggl(x^{-j+1} \Bigl( \frac{\partial}{\partial \varphi^-} -  \varphi^+  \frac{\partial}{\partial x}  \Bigr)  \\
& & \ \  -(- j+1) x^{-j} \varphi^+ \varphi^- \frac{\partial}{\partial \varphi^-} \Bigr) \! \biggr) \!  \Biggr) \!  \cdot \! \left(\frac{1}{x},\frac{i\varphi^+}{x}, \frac{i\varphi^-}{x} \right) \\
&=&\! \! \exp\Biggl(  \sum_{j \in \Z} \Bigl( A^+_j L_{-j}(x, \varphi^+, \varphi^-) - A^-_j  J_{-j}(x, \varphi^+, \varphi^-)  \\
& & \ \  + i M^+_{j - \frac{1}{2}} G^+_{-j+ \frac{1}{2}}(x, \varphi^+, \varphi^-)  +i M^-_{j - \frac{1}{2}} G^-_{-j+ \frac{1}{2}}(x, \varphi^+, \varphi^-)  \Biggr) \!  \cdot \!  \left(\frac{1}{x},\frac{i\varphi^+}{x}, \frac{i\varphi^-}{x} \right) 
\end{eqnarray*}  
which gives (\ref{infty}).

By Proposition \ref{inverses}, we know that $H_{-1} (x, \varphi^+, \varphi^-)$ has a unique inverse, namely  $H_{-1}^{-1} (x, \varphi^+, \varphi^-)$, with
\begin{eqnarray*}
H_{-1}^{-1} (x, \varphi^+, \varphi^-) &=& \exp(\mathcal{T}_{H_{-1}^{-1}} (x, \varphi^+, \varphi^-) ) \cdot
(x, \varphi^+, \varphi^-) \\
&=& \exp(- \mathcal{T}_{H_{-1}} (x, \varphi^+, \varphi^-) ) \cdot (x, \varphi^+, \varphi^-) .
\end{eqnarray*}

Setting $H^{-1} (x, \varphi^+, \varphi^-) = I^{-1} \circ H_{-1}^{-1} (x, \varphi^+, \varphi^-)$, and since $H (x, \varphi^+, \varphi^-) = H_{-1} \circ I (x, \varphi^+, \varphi^-)$,
we have 
\[H \circ H^{-1} (x, \varphi^+, \varphi^-) = H_{-1} \circ I \circ I^{-1} \circ H_{-1}^{-1} (x, \varphi^+, \varphi^-) =   (x, \varphi^+, \varphi^-) \]
and
\[H^{-1} \circ H (x, \varphi^+, \varphi^-) = I^{-1} \circ H_{-1}^{-1} \circ H_{-1} \circ I (x, \varphi^+, \varphi^-) = (x, \varphi^+, \varphi^-) .\]
Moreover, by Proposition \ref{Switch}, with $\overline{H}  = I^{-1}$ and $H$ in the Proposition replaced by $H_{-1}^{-1}$, we have
\begin{eqnarray*}
\lefteqn{H^{-1} \circ I (x, \varphi^+, \varphi^-) }\\
&=& I^{-1} \circ H_{-1}^{-1} \circ I (x, \varphi^+, \varphi^-) \\
&=&  \left.  I^{-1}  \left(\exp (-\mathcal{T}_{H_{-1}} (x, \varphi^+ , \varphi^-)) \cdot (x, \varphi^+, \varphi^-) \right) \right|_{(x, \varphi^+ , \varphi^- ) = I(x,\varphi^+, \varphi^-)}\\
&=& \left. \exp (- \mathcal{T}_{H_{-1}} (x, \varphi^+ , \varphi^-) ) 
\cdot I^{-1}(x, \varphi^+ , \varphi^-)   \right|_{(x, \varphi^+ , \varphi^- ) = I(x,\varphi^+, \varphi^-)}\\\  
&=& \exp \Bigl(- \mathcal{T}_{H_{-1}} \Bigl(\frac{1}{x}, \frac{i \varphi^+}{x} , \frac{i \varphi^-}{x} 
\Bigr) \Bigr) \cdot (x, \varphi^+, \varphi^-) 
\end{eqnarray*}
which gives (\ref{forsewing}).
\end{proof}

\begin{rema} {\em The formal $N=2$ superconformal power series of the form (\ref{infty}) can be thought of as the ``formal $N=2$ superconformal local coordinate maps vanishing at $\infty = (\infty, 0,0)$ with leading even coefficient of $\varphi^\pm x^{-1}$ equal to $i$''. }
\end{rema}

The following two propositions are analogous to Proposition \ref{auto2} and Proposition \ref{Switch}, respectively.

\begin{prop}\label{auto3}
Let $u, v \in R((x^{-1}))[\varphi^+, \varphi^-]$; let $(B^+, B^-,N^+,N^-) \in (R^\infty)^2$; and let  
\begin{eqnarray}
 \bar{\mathcal{T}} &=& -  \sum_{j \in \Z} \Bigl( B^+_j L_{-j}(x,\varphi^+, \varphi^-) +  B^-_j J_{-j}(x,\varphi^+, \varphi^-)  \label{T-bar} \\
& & \quad + N^+_{j - \frac{1}{2}} G^+_{-j + \frac{1}{2}}(x,\varphi^+, \varphi^+ ) + N^-_{j -
\frac{1}{2}} G^-_{-j + \frac{1}{2}}(x,\varphi^+, \varphi^+ ) \Bigr) . \nonumber 
\end{eqnarray}
Then 
\begin{equation} \label{specialized-automorphism-property-at-infinity} 
e^{\bar{\mathcal{T}}} \cdot (uv) = \left( e^{\bar{\mathcal{T}}} \cdot u \right)\left( e^{\bar{\mathcal{T}}} \cdot v  \right). 
\end{equation}
\end{prop}

In other words, for $\bar{\mathcal{T}}$, $u$, and $v$ given above, Proposition \ref{conformalproof} holds if $h=0$ and $y$ is set equal to 1. 

\begin{proof} The proof is analogous to the proof of Proposition \ref{auto2}.   Again, we use a $\frac{1}{2} \mathbb{Z}$ grading by weight on $R((x^{-1}))[\varphi^+, \varphi^-]$ given by $\mathrm{wt} \ a = 0$ if $a \in R$, $\mathrm{wt}\ x^k = k$, for $k \in \mathbb{Z}$, and $\mathrm{wt}\ \varphi^\pm = \frac{1}{2}$.  Then $\bar{\mathcal{T}} \in \mathrm{Der} (R((x^{-1}))[\varphi^+, \varphi^-])$ given by (\ref{T-bar}) acting on an element $u' \in R((x^{-1}))[\varphi^+, \varphi^-]$ of homogeneous weight, lowers the weight of $u'$ by at least $\frac{1}{2}$.  Thus if $\mathrm{wt} \ u' = k \in \frac{1}{2} \mathbb{Z}$, then $\frac{y^n \bar{\mathcal{T}}^n}{n!} \cdot u'$ involves terms of weight lower than or equal to $k - \frac{n}{2}$, for $n \in \Z$.   Therefore, for $u \in R((x^{-1}))[\varphi^+, \varphi^-]$, any term in $e^{y\bar{T}} \cdot u$ of homogeneous weight in $x$, $\varphi^+$, and $\varphi^-$ is a polynomial in $y$, and thus we can set $y= 1$.  That is $e^{\bar{\mathcal{T}}} \cdot u$ has only a finite number of terms of weight $k \in \frac{1}{2} \mathbb{Z}$, and thus is a well-defined power series in $R((x^{-1}))[\varphi^+, \varphi^-]$, and similarly for $e^{\bar{\mathcal{T}}} \cdot v$.  Therefore (\ref{specialized-automorphism-property-at-infinity}) holds. 
 \end{proof}

\begin{prop}\label{Switch2}
Let $H(x, \varphi^+, \varphi^-) = e^{\bar{\mathcal{T}}} \cdot (x, \varphi^+, \varphi^-)$ with $\bar{\mathcal{T}}$ given by (\ref{T-bar}), and let $\overline{H}(x, \varphi^+, \varphi^-) \in R((x^{-1}))[\varphi^+, \varphi^-]$.  Then   
\begin{equation}\label{switch2}
\overline{H}(H(x, \varphi^+, \varphi^-)) = \overline{H}(e^{\bar{\mathcal{T}}} \cdot (x, \varphi^+, \varphi^-)) = e^{\bar{\mathcal{T}}} \cdot \overline{H}(x,\varphi^+, \varphi^-) .  
\end{equation}
More generally, if $\overline{H}(x, \varphi^+, \varphi^-) \in R[[x,x^{-1}]][\varphi^+, \varphi^-]$ and the composition $\overline{H} \circ H$ is well defined, then (\ref{switch2}) holds. 
\end{prop}

\begin{proof}  The first part of the proof is identical to steps (i) -- (iv) in the proof of Proposition \ref{Switch}.  To finish the proof, we only need note that since $e^{\bar{\mathcal{T}}} \cdot (\varphi^+)^{j^+} (\varphi^-)^{j^-} x^n \in R((x^{-1}))[\varphi^+, \varphi^-]$ for $j^\pm = 0, 1$ and $n \in \mathbb{Z}$, the result follows by linearity.  \end{proof}   

Let $H(x,\varphi^+, \varphi^-)$ and $\overline{H}(x, \varphi^+, \varphi^-)$ be two formal $N=2$ superconformal power series vanishing at infinity and with even coefficients of $\varphi^\pm x^{-1}$ equal to $i$, i.e., let $H$ and $\bar{H}$ be of the form (\ref{infty}) such that  
\begin{eqnarray}
\tilde{E}^{-1}(H\circ I^{-1} (x, \varphi^+, \varphi^-)) &=& (A^+, -A^-,-iM^+,-iM^-) \label{for-composition1-infinity}\\
\tilde{E}^{-1}(\overline{H} \circ I^{-1} (x, \varphi^+,\varphi^-)) &=& (B^+,-B^-,-iN^+,-iN^-) \label{for-composition2-infinity}
\end{eqnarray}  
for some $(A^+,A^-,M^+,M^-), (B^+,B^-,N^+,N^-) \in (R^\infty)^2$.  Define the {\em composition at infinity} of 
$(A^+,A^-,M^+,M^-)$ and $(B^+,B^-,N^+,N^-)$, denoted 
\begin{equation}
(A^+,A^-,M^+,M^-) \circ_\infty (B^+,B^-,N^+,N^-),
\end{equation} 
by
\begin{equation}\label{sequencecompositiondef-infty}
(A^+,A^-,M^+,M^-) \circ_\infty (B^+,B^-,N^+,N^-) = (C^+, C^-, P^+, P^-) 
\end{equation}
for
\begin{equation}\label{sequencecompositiondef-infty2}
(C^+, -C^-,-i P^+, -iP^-) =   \tilde{E}^{-1}((H \circ I^{-1} \circ \overline{H} \circ I^{-1} ) (x, \varphi^+,\varphi^-)) \in (R^\infty)^2,
\end{equation}
where $(H \circ I^{-1} \circ \overline{H} \circ I^{-1} ) (x, \varphi^+,\varphi^-)$ is the formal composition of $H$, $I^{-1}$  and $\overline{H}$.  

We have the following corollary in analogy to Proposition \ref{sequencecompositionprop}.

\begin{cor} \label{sequence-composition-cor-infty}
The set $(R^\infty)^2$ is a group with the operation $\circ_\infty$, and letting
\begin{equation}
(C^+, -C^-,-i P^+, -iP^-) =  (B^+,-B^-,-iN^+,-iN^-) \circ (A^+,-A^-,-iM^+,-iM^-),
\end{equation}
where the composition $\circ$ on $(R^\infty)^2$ is that defined by (\ref{sequencecompositiondef}), we have
\begin{equation}
(A^+,A^-,M^+,M^-) \circ_\infty (B^+,B^-,N^+,N^-)  = (C^+, C^-,  P^+,  P^-).  
\end{equation}

Furthermore, letting  
\[(A^+,A^-,M^+,M^-) = \{ (A^+_j,A^-_j , M^+_{j - 1/2},M^-_{j - 1/2}) \}_{j \in \Z} \in (R^\infty)^2,\] 
and for any $t \in R^0$, if we define 
\[t(A^+,A^-,M^+,M^-) = \{ (t A^+_j,tA^-_j, t M^+_{j - 1/2}, tM^-_{j-1/2}) \}_{j \in \Z}, \]  
then for $t_1,t_2 \in R^0$,   
\begin{multline}\label{composition-infty}
(t_1(A^+,A^-,M^+,M^-)) \circ_\infty (t_2(A^+,A^-,M^+,M^-)) \\
= (t_1 + t_2)(A^+,A^-,M^+,M^-) .
\end{multline} 
That is the map $t \mapsto t(A^+,A^-,M^+,M^-)$ is a homomorphism {}from the additive group $R^0$ to the group  $(R^\infty)^2$ with group operation $\circ_\infty$. 
\end{cor}

\begin{proof}  Let $H \circ I^{-1} = \tilde{E} (A^+,-A^-,-iM^+,-iM^-)$ and $\overline{H} \circ I^{-1} =  \tilde{E} (B^+,-B^-,-iN^+,$ $-iN^-)$.  Then by definition 
\begin{multline*}
(B^+,-B^-,-iN^+,-iN^-) \circ (A^+,-A^-,-iM^+,-iM^-)\\
 =  \tilde{E}^{-1} ( H\circ I^{-1} \circ \overline{H} \circ I^{-1}(x, \varphi^+ , \varphi^-) )
 \end{multline*}
proving the first statement.  The rest of the corollary then follows {}from Proposition  \ref{sequencecompositionprop}.
\end{proof}

\begin{rema}\label{group-remark-infty} {\em 
With the composition operation $\circ_\infty$ defined above, $(R^\infty)^2$ can be thought of as the group of formal $N=2$ superconformal local coordinate transformations vanishing at infinity and with leading even coefficient of $\varphi^\pm x^{-1}$ equal to $i$.   In Section \ref{NS-algebra-section} we discuss subgroups of $(R^\infty)^2$ under the group operations $\circ_\infty$ and $\circ$, and in Section \ref{moduli-group-section} we relate these group structures to subsets of the moduli space of $N=2$ super-Riemann spheres with one incoming tube and one outgoing tube and a corresponding ``sewing" operation.   }
\end{rema}

\begin{rema} \label{general-formal-infinitesimal}
{\em {}From Theorem \ref{above} and Proposition \ref{Infinity},  we see that in general the ``formal infinitesimal $N=2$ superconformal transformations'' are of the form 
\begin{multline}\label{infinitesimal}
\sum_{j \in \Z} \biggl(B^+_j L_{-j}(x,\varphi^+, \varphi^-) + B^-_j J_{-j}(x,\varphi^+, \varphi^-) + N^+_{j - \frac{1}{2}} G^+_{-j + \frac{1}{2}} (x, \varphi^+, \varphi^-) \\
+ N^-_{j - \frac{1}{2}} G^-_{-j + \frac{1}{2}} (x, \varphi^+, \varphi^-) \biggr)  + (\log a_0^+) 2L_0(x,\varphi^+, \varphi^-)  + (\log a_0^-) J_0(x,\varphi^+, \varphi^-) \\
+  \sum_{j \in \Z} \biggl( A_j^+ L_j(x,\varphi^+, \varphi^-)  + A_j^- J_j(x,\varphi^+, \varphi^-)  + M^+_{j - \frac{1}{2}} G^+_{j - \frac{1}{2}}(x,\varphi^+, \varphi^- ) \\
+ M^-_{j - \frac{1}{2}} G^-_{j - \frac{1}{2}}(x,\varphi^+, \varphi^- )\biggr)  
\end{multline}  
for $(a_0^+, a_0^-) \in ((R^0)^\times)^2/\langle \pm 1 \rangle$, $A^\pm_j, B^\pm_j \in R^0$, and $M^\pm_{j - 1/2}, N^\pm_{j - 1/2} \in R^1$.  }
\end{rema}

\section[The $N=2$ Neveu-Schwarz algebra]
{The $N=2$ Neveu-Schwarz algebra and the group of $N=2$  superprojective transformations}\label{NS-algebra-section}

In this section, we define the $N = 2$ Neveu-Schwarz algebra \cite{DPYZ} and point out that the superderivations we used in Section \ref{infinitesimals-section} give a representation of the $N=2$ Neveu-Schwarz algebra with central charge zero \cite{Ki}.  This shows that the $N=2$ Neveu-Schwarz algebra is the algebra of infinitesimal  $N=2$ superconformal transformations.  We point out several subalgebras of the $N=2$ Neveu-Schwarz algebra and their corresponding Lie supergroups in connection with our development of the groups corresponding to $N=2$ superconformal local coordinates vanishing at zero and at infinity as developed in Section \ref{infinitesimals-section}.  We discuss the subalgebra of the $N=2$ Neveu-Schwarz algebra consisting of infinitesimal global $N=2$ superconformal transformations (i.e., infinitesimal $N=2$ superprojective transformations), and derive the action of the corresponding Lie supergroup of $N=2$ superprojective transformations on the $N=2$ super-Riemann sphere $S^2\hat{\mathbb{C}}$.  In addition, we point out several errors in the literature concerning the presentation of these $N=2$ superprojective transformations.  

Let $\mathfrak{ns}_2$ denote the $N=2$ Neveu-Schwarz Lie superalgebra with central charge $d$, basis consisting of the central element $d$, even elements $L_n$ and $J_n$ and odd elements $G^\pm_{n + 1/2}$, for $n \in \mathbb{Z}$, and commutation relations  
\begin{eqnarray}
\left[L_m ,L_n \right] &=& (m - n)L_{m + n} + \frac{1}{12} (m^3 - m) \delta_{m + n , 0} \; d , \label{Virasoro-relation} \\
\left[J_m , J_n \right] &=& \frac{1}{3} m \delta_{m + n , 0} \; d , \label{NS-relation1} 
\end{eqnarray}
\begin{eqnarray}
\left[L_m , J_n \right] &=&  -nJ_{m+n} , \label{NS-relation2} \\
\bigl[L_m,G^\pm_{n + \frac{1}{2}}\bigr] &=& \Bigl(\frac{m}{2} - n - \frac{1}{2} \Bigr) G^\pm_{m + n + \frac{1}{2}} ,  \label{NS-relation3} \\ 
\bigl[ J_m , G^\pm_{n + \frac{1}{2}}\bigr] &=& \pm G^\pm_{m + n + \frac{1}{2}} , \label{NS-relation4} \\ 
\bigl[ G^+_{m + \frac{1}{2}} , G^+_{n + \frac{1}{2}} \bigr] &=& \bigl[ G^-_{m + \frac{1}{2}} , G^-_{n + \frac{1}{2}} \bigr] \; = \; 0, \label{NS-relation5} \\
\hspace{.4in} \bigl[ G^+_{m + \frac{1}{2}} , G^-_{n - \frac{1}{2}} \bigr] &=& 2L_{m + n} + (m-n+1) J_{m+n} \label{NS-relation6} \\
& &  \hspace{1.3in} + \; \frac{1}{3}(m^2 + m) \delta_{m + n , 0} \; d . \nonumber
\end{eqnarray}

It is easy to check that the superderivations in $\mbox{Der}(\mathbb{C} [x, x^{-1},\varphi^+, \varphi^-])$ given by (\ref{L-notation}) -- (\ref{G-notation}),  i.e., the superderivations  
\begin{eqnarray}
L_n(x,\varphi^+,\varphi^-) &=& - \Bigl( x^{n + 1} \frac{\partial}{\partial x} + \bigl(\frac{n + 1}{2}\bigr) x^n \Bigl( \varphi^+ \frac{\partial}{\partial \varphi^+} + \varphi^- \frac{\partial}{\partial \varphi^-} \Bigr) \Bigr) 
\label{L(n)}\\  
J_n(x,\varphi^+,\varphi^-) &=& - x^n\Bigl(\varphi^+\frac{\partial}{\partial \varphi^+} - \varphi^- \frac{\partial}{\partial \varphi^-}\Bigr) \label{J(n)}\\
\qquad G^\pm_{n -\frac{1}{2}}(x,\varphi^+,\varphi^-) &=& - \Bigl( x^n \Bigl( \frac{\partial}{\partial \varphi^\pm} - \varphi^\mp \frac{\partial}{\partial x} \Bigr) \pm nx^{n-1} \varphi^+ \varphi^- \frac{\partial}{\partial \varphi^\pm} \Bigr) \label{G(n)}
\end{eqnarray}
for $n \in \mathbb{Z}$, satisfy the $N = 2$ Neveu-Schwarz relations (\ref{Virasoro-relation}) -- (\ref{NS-relation6}) with central charge zero (cf. \cite{Ki}).

We would now like to discuss some subalgebras of $\mathfrak{ns}_2$ and their interpretation in terms of the group structures discussed in Proposition \ref{sequencecompositionprop} and Corollary \ref{sequence-composition-cor-infty}.  We first note that letting $Z = \pm \Z$, the following are subalgebras of $\mathfrak{ns}_2$.
\begin{equation}
\qquad \ \  \{ L_n \ | \ n \in Z \}, \quad \{J_n \ | \ n \in Z \}, \quad \{L_m, J_n \ | \ m,n \in Z \} \label{subalgebra1}
\end{equation}
\vspace{-.2in}
\begin{eqnarray}
\{ G^+_{n + \frac{1}{2}} \ | \ n \in Z \}, &\quad& \{ G^-_{n + \frac{1}{2}} \ | \ n \in Z \},\\
\{L_m, G^+_{n + \frac{1}{2}} \ | \ m,n \in Z \}, &\quad& \{L_m, G^-_{n + \frac{1}{2}} \ | \ m,n \in Z \}, \\
\{J_m, G^+_{n + \frac{1}{2}} \ | \ m,n \in Z \}, &\quad& \{J_m, G^-_{n + \frac{1}{2}} \ | \ m,n \in Z \}, \\
\{L_k, J_m, G^+_{n + \frac{1}{2}} \ | \ k,m,n \in Z \}, & \quad& \{L_k, J_m,  G^-_{n + \frac{1}{2}} \ | \ k,m,n \in Z \}, \label{subalgebra-next-to-last}  
\end{eqnarray}
\vspace{-.2in}
\begin{equation}
\{L_k, J_l,  G^+_{m + \frac{1}{2}}  G^-_{n + \frac{1}{2}} \ | \ k,l,m,n \in Z \} .\label{subalgebra-last}
\end{equation}
(This is of course not an exhaustive list.)   Interpreting this in terms of the groups $((R^\infty)^2, \circ)$ and $((R^\infty)^2, \circ_\infty)$, respectively, defined in Section \ref{infinitesimals-section} for a superalgebra $R$, we have the following corollary.

\begin{cor}\label{subgroups-corollary}
The following are subgroups of $((R^\infty)^2, \circ)$ and of $((R^\infty)^2, \circ_\infty)$ in correspondence with the subalgebras given in (\ref{subalgebra1}) -- (\ref{subalgebra-next-to-last}) for $Z = \pm \Z$, respectively.  \footnote{ In \cite{B-memoir}, Proposition 3.15 (which is the analogue of Proposition \ref{sequencecompositionprop} for the $N=1$ superconformal case) states that $(R^1)^\infty$ is a subgroup of $R^\infty$ under composition in $R^\infty$ defined using formal $N=1$ superconformal local coordinates vanishing at zero.  This is false; $(R^0)^\infty$ is a subgroup, but $(R^1)^\infty$ is not.}
\begin{equation}
\{ (A^+,\mathbf{0}, \mathbf{0}, \mathbf{0}) \; | \; A^+ \in (R^0)^\infty\}, \   \{ (\mathbf{0},A^-, \mathbf{0}, \mathbf{0}) \; | \; A^- \in (R^0)^\infty\}, \   ((R^0)^\infty)^2, 
\end{equation}
\vspace{-.25in}
\begin{eqnarray}
\{ (\mathbf{0}, \mathbf{0},M^+, \mathbf{0}) \; | \; M^+ \in (R^1)^\infty\},  &  \{ (\mathbf{0}, \mathbf{0}, \mathbf{0}, M^-) \; | \; M^- \in (R^1)^\infty\},  \\
\quad \quad \ \  \{ (A^+,\mathbf{0},M^+, \mathbf{0}) \; | \; (A^+,M^+) \in R^\infty\}, &   \{ (A^+,\mathbf{0}, \mathbf{0}, M^-) \; | \; (A^+,M^-) \in R^\infty\}
\end{eqnarray}
\begin{eqnarray}
\qquad \ \  \{ (\mathbf{0},A^-, M^+, \mathbf{0}) \; | \; (A^-,M^+) \in R^\infty\}, &  \{ (\mathbf{0},A^-, \mathbf{0}, M^-) \; | \; (A^-,M^-) \in R^\infty\}
\end{eqnarray}
\vspace{-.3in}
\begin{eqnarray}
\{ (A^+,A^-,M^+, \mathbf{0}) \; | \; A^+,A^- \in (R^0)^\infty, \; M^+ \in (R^1)^\infty\}, \\
\{ (A^+,A^-, \mathbf{0}, M^-) \; | \;  A^+,A^- \in (R^0)^\infty, \; M^- \in (R^1)^\infty\}, 
\end{eqnarray}
where $\mathbf{0}$ denotes the sequence in $(R^j)^\infty$ consisting of all zeros, for $j = 0,1$.
The subalgebras given by (\ref{subalgebra-last}) for $Z = \pm \Z$, respectively, correspond to the groups $((R^\infty)^2, \circ)$ and $((R^\infty)^2, \circ_\infty)$, respectively.
\end{cor}

We now focus on another subalgebra of $\mathfrak{ns}_2$, that  given by 
\begin{equation}
\mbox{span}_{\mathbb{C}} \{L_{\pm1}, L_0,  J_0, G^+_{\pm 1/2}, G^-_{\pm 1/2} \}.
\end{equation} 

Let $W$ be a $\mathbb{Z}_2$-graded vector space over $\mathbb{C}$ such that $\mathrm{dim} \ W^0 =  \mathrm{dim} \ W^1 = 2$.  Recall the classical Lie superalgebra $\mathfrak{osp}_{\mathbb{C}}(2|2)$ (cf. \cite{Ka}) the orthogonal-symplectic superalgebra  
\begin{eqnarray*}
\mathfrak{osp}_{\mathbb{C}}(2|2) = \left\{ \biggl. \left(\begin{array}{cccc}   
                                         e & 0 & p & q \\
                                         0 & -e & r & s \\
                                         s & q  & a & b \\
                                         -r &-p & c & -a 
                                            \end{array} \right) \in
\mathfrak{gl}_{\mathbb{C}}(2|2) \; \biggr| \; a,b,c,e,p,q, r, s\in \mathbb{C} \right\} 
\end{eqnarray*}
which is the subalgebra of $\mathfrak{gl}_{\mathbb{C}}(2|2)$ leaving the non-degenerate form  $\beta$ on $W$ given by 
\[\beta = \left(\begin{array}{cccc}  
                0 & 1 & 0 & 0\\
                1 & 0 & 0 & 0\\
                0 & 0 & 0 & 1\\
                0 & 0 & -1 & 0 \end{array}
\right)\]
invariant, meaning $\beta (Xu,v) + (-1)^{\eta(X) \eta(u)} 
\beta(u,Xv) = 0$ for $X \in \mathfrak{osp}_{\mathbb{C}}(2|2)$, $u,v \in W$, and $X$ and $u$ homogeneous. 

The subalgebra of $\mathfrak{ns}_2$ given by $\mbox{span}_{\mathbb{C}} \{L_{\pm1}, L_0,  J_0, G^+_{\pm 1/2}, G^-_{\pm 1/2} \}$ is isomorphic to $\mathfrak{osp}_{\mathbb{C}}(2|2)$. The
correspondence 
\begin{eqnarray*}
\left(\begin{array}{cccc}  
        0 & 0 & 0 & 0 \\ 
        0 & 0 & 0 & 0 \\
        0 & 0 & 0 & 1 \\
        0 & 0 & 0 & 0 \end{array} \right) \longleftrightarrow & \hspace{-.1in}  -  \frac{\partial}{\partial x}  & \hspace{-.1in} = L_{-1}(x, \varphi^+, \varphi^-),  \\
\\
\frac{1}{2} \left( \begin{array}{cccc}  
         0 & 0 & 0 & 0 \\  
         0 & 0 & 0 & 0 \\
         0 & 0 & 1 & 0 \\
         0 & 0 & 0 & -1 \end{array} \right)  \longleftrightarrow& \hspace{-.1in}   - \left( x \frac{\partial}{\partial x} + \frac{1}{2} \left( \varphi^+ \frac{\partial}{\partial \varphi^+ } +  \varphi^- \frac{\partial}{\partial \varphi^- } \right)  \right) & \hspace{-.1in} = L_0(x, \varphi^+, \varphi^-), \\ 
\\
\left(\begin{array}{cccc}   
0 & 0 & 0 & 0 \\
0 & 0 & 0 & 0 \\
0 & 0 & 0 & 0 \\
0 & 0 & -1 & 0 \end{array} \right) \longleftrightarrow& \hspace{-.1in}   - \left( x^2 \frac{\partial}{\partial x} +  x \left( \varphi^+ \frac{\partial}{\partial \varphi^+ } +  \varphi^- \frac{\partial}{\partial \varphi^- } \right)\right) & \hspace{-.1in}  = L_1(x, \varphi^+, \varphi^-), \\ 
\\
\left(\begin{array}{cccc}   
1 & 0 & 0 & 0 \\
0 & -1 & 0 & 0 \\
0 & 0 & 0 & 0 \\
0 & 0 & 0 & 0 \end{array} \right) \longleftrightarrow& \hspace{-.1in}    - \left(  \varphi^+ \frac{\partial}{\partial \varphi^+ } -  \varphi^- \frac{\partial}{\partial \varphi^- } \right) & \hspace{-.1in} = J_0(x, \varphi^+, \varphi^-), \\ 
\end{eqnarray*}
\begin{eqnarray*}
\left(\begin{array}{cccc}
0 & 0 & 0 & 1 \\   
0 & 0 & 0 & 0 \\
0 & 1 & 0 & 0 \\
0 & 0 & 0 & 0 \end{array} \right)  \longleftrightarrow& \hspace{-.1in}   - \left(  \frac{\partial}{\partial \varphi^+} -  \varphi^- \frac{\partial}{\partial x} \right) & \hspace{-.1in}  = G^+_{-\frac{1}{2}} (x, \varphi^+, \varphi^-), \\
\\
\left(\begin{array}{cccc}    
0 & 0 & 1 & 0 \\
0 & 0 & 0 & 0 \\
0 & 0 & 0 & 0 \\
0 & -1 & 0 & 0 \end{array} \right) \longleftrightarrow&  \hspace{-.1in}  - \left( x\left( \frac{\partial}{\partial \varphi^+ } -  \varphi^-  \frac{\partial}{\partial x} \right) + \varphi^+ \varphi^- \frac{\partial}{\partial \varphi^+} \right) & \hspace{-.1in} = G^+_{\frac{1}{2}} (x, \varphi^+ \varphi^-),  \\
\\
\left(\begin{array}{cccc}
0 & 0 & 0 & 0 \\   
0 & 0 & 0 & 1 \\
1 & 0 & 0 & 0 \\
0 & 0 & 0 & 0 \end{array} \right)  \longleftrightarrow& \hspace{-.1in}  - \left(  \frac{\partial}{\partial \varphi^-} -  \varphi^+ \frac{\partial}{\partial x} \right) & \hspace{-.1in} = G^-_{-\frac{1}{2}} (x, \varphi^+, \varphi^-), \\
\\
\left(\begin{array}{cccc}    
0 & 0 & 0 & 0 \\
0 & 0 & 1 & 0 \\
0 & 0 & 0 & 0 \\
-1 & 0 & 0 & 0 \end{array} \right) \longleftrightarrow&  \hspace{-.1in} - \left( x\left( \frac{\partial}{\partial \varphi^- } -  \varphi^+  \frac{\partial}{\partial x} \right) - \varphi^+ \varphi^- \frac{\partial}{\partial \varphi^-} \right) & \hspace{-.1in} = G^-_{\frac{1}{2}} (x, \varphi^+ \varphi^-),  
\end{eqnarray*} 
defines a Lie superalgebra isomorphism between $\mathfrak{osp}_{\mathbb{C}}(2|2)$ and the Lie superalgebra of infinitesimal $N=2$ superconformal transformations generated by $L_{\pm1}(x,\varphi^+, \varphi^-)$, $L_{0}(x,\varphi^+, \varphi^-),  J_{0}(x,\varphi^+, \varphi^-), G^+_{\pm 1/2}(x,\varphi^+, \varphi^-)$, and  $G^-_{\pm 1/2}(x,\varphi^+, \varphi^-)$.   Note that this isomorphism is not unique; for example, the automorphism of $\mathfrak{osp}_{\mathbb{C}}(2|2)$  given by $J_0 \mapsto -J_0$, $G^\pm_{1/2} \mapsto G^\mp_{1/2}$, $G^\pm_{-1/2} \mapsto G^\mp_{-1/2}$, and $L_j \mapsto L_j$ for $j = 0,\pm1$ gives another identification.

Let $y$ be an even formal variable and $\xi$ an odd formal variable.  Letting $X$ denote each of the eight matrices above, we observe that
\begin{eqnarray*}
e^{-yX} = \left(\begin{array}{cccc}  
                  1 & 0 & 0 & 0 \\ 
                  0& 1 & 0 & 0 \\
                  0 & 0 & 1 & -y \\
                  0 & 0 & 0 & 1 \end{array} \right) , \quad
\left(\begin{array}{cccc}    
1 & 0 & 0 & 0 \\
0 & 1 & 0 & 0 \\
0 & 0 & e^{-\frac{y}{2}} & 0 \\
0 & 0 & 0 & e^{\frac{y}{2}}  \end{array} \right) , \quad
\left(\begin{array}{cccc}   
1 & 0 & 0 & 0 \\ 
0 & 1 & 0 & 0 \\
0 & 0 & 1 & 0 \\
0 & 0 & y & 1 \end{array} \right) 
\end{eqnarray*}
\begin{equation*}
\left(\begin{array}{cccc}    
e^{-y} & 0 & 0 & 0 \\
0 & e^y & 0 & 0 \\
0 & 0 & 1 & 0 \\
0 & 0 & 0 & 1 \end{array} \right) 
\end{equation*} 
\begin{eqnarray*}
e^{-\xi X} =  \left(\begin{array}{cccc}   
                      1 & 0 & 0 & -\xi \\
                       0 & 1 & 0 & 0 \\
                      0 & -\xi & 1 & 0 \\
                      0 &  0 & 0 & 1 \end{array} \right) , \quad 
\left(\begin{array}{cccc}   
1 & 0 & -\xi & 0 \\
0 & 1 & 0 & 0 \\
0 & 0 & 1 & 0 \\
0 & \xi & 0 & 1 \end{array}\right), \quad
\left(\begin{array}{cccc}   
                      1 & 0 & 0 & 0 \\
                       0 & 1 & 0 & -\xi  \\
                      -\xi & 0  & 1 & 0 \\
                      0 &  0 & 0 & 1 \end{array} \right) , 
\end{eqnarray*}
\begin{equation*}
\left(\begin{array}{cccc}   
1 & 0 & 0 & 0 \\
0 & 1 & -\xi  & 0 \\
0 & 0 & 1 & 0 \\
\xi & 0 & 0 & 1 \end{array}\right),
\end{equation*}
respectively.  These are all elements in the connected component of the Lie supergroup $OSP(2|2)$ containing the identity with matrix elements in $\mathbb{C}[[y]][\xi]$.  They have superdeterminant 1, 
where the superdeterminant is defined as  
\[\mbox{sdet}\left(\begin{array}{cc}  A & B \\
                                      C & D \end{array}\right) =
\mbox{det}(A - BD^{-1}C)(\mbox{det}D)^{-1} .\]  
(In our case, $A, B , C$ and $D$ are all two-by-two matrices.)  In fact, 
for $R$ a superalgebra with $y \in R^0$ and $\xi \in R^1$, the eight matrices above generate the connected component of $OSP_R(2|2)$ containing the identity (cf. \cite{D}, \cite{Var}).  Denote this group by $\mathcal{G}$.

$\mathcal{G}$ acts on an even and an odd formal variable by the $N=2$ superprojective transformations, i.e,  for $g = e^{-yX}$ and $g = e^{-\xi X}$ above, we have
\begin{eqnarray}\label{group action}
\lefteqn{g \cdot (x,\varphi^+, \varphi^-)} \\
&=& (x + y, \varphi^+, \varphi^-), \quad (e^yx, e^{\frac{y}{2}} \varphi^+, e^{\frac{y}{2}} \varphi^-), \nonumber \\
& & \Bigl(\frac{x}{1 - yx}, \varphi^+ \frac{1}{1- yx},  \varphi^- \frac{1}{1- yx}\Bigr), \quad (x, e^y\varphi^+, e^{-y} \varphi^-)  \nonumber \\
& & (x + \varphi^- \xi, \xi + \varphi^+, \varphi^-), \quad (x + \varphi^- \xi x, \xi  x + \varphi^+ + \varphi^+ \varphi^- \xi, \varphi^- ) , \nonumber\\
& & (x + \varphi^+ \xi,  \varphi^+, \xi + \varphi^-), \quad (x + \varphi^+ \xi x,\varphi^+,  \xi  x + \varphi^- - \varphi^+ \varphi^- \xi) , \nonumber
\end{eqnarray}
respectively.  These generate the supergroup of $N=2$ superprojective transformations which is the group of global $N=2$ superconformal automorphisms of the super-Riemann sphere studied in Section \ref{moduli-section}.   Thus, $\mathfrak{osp}_{\mathbb{C}}(2|2)$ is the Lie superalgebra of infinitesimal $N=2$ superprojective transformations.  Note that for the representative elements $L_{\pm 1} (x, \varphi^+, \varphi^-)$, $L_0(x, \varphi^+, \varphi^-)$, $J_0(x, \varphi^+, \varphi^-)$, $G^+_{\pm 1/2}(x, \varphi^+, \varphi^-)$, and $G^-_{\pm 1/2}(x, \varphi^+, \varphi^-)$ in $\mathrm{Der} (\mathbb{C} [x, x^{-1}, \varphi^+, \varphi^-])$, we have 
\begin{eqnarray}
e^{-yL_{-1}(x, \varphi^+, \varphi^-)} \cdot (x,\varphi^+,\varphi^-) &=&  (x + y, \varphi^+,\varphi^-),  \label{group-action1} \\
e^{-yL_0(x, \varphi^+, \varphi^-)} \cdot (x,\varphi^+,\varphi^-) &=&  (e^y x,e^{\frac{y}{2}}\varphi^+, e^{\frac{y}{2}} \varphi^-) ,\\
e^{-yL_1(x, \varphi^+, \varphi^-)} \cdot (x,\varphi^+,\varphi^-) &=&  \left(\frac{x}{1-yx}, \varphi^+ \frac{1}{1-yx}, \varphi^- \frac{1}{1-yx}\right), \label{group-action3} \\
e^{-yJ_0(x, \varphi^+, \varphi^-)} \cdot (x,\varphi^+,\varphi^-) &=&  (x,e^y \varphi^+, e^{-y} \varphi^-) ,\\
\qquad e^{-\xi G^+_{-1/2}(x, \varphi^+, \varphi^-)} \cdot (x,\varphi^+,\varphi^-) &=& (x + \varphi^-  \xi , \xi + \varphi^+, \varphi^-), \\ 
e^{-\xi G^+_{1/2}(x, \varphi^+, \varphi^-)} \cdot (x,\varphi^+,\varphi^-) &=& (x + \varphi^- \xi x, \xi x + \varphi^+ + \varphi^+ \varphi^- \xi, \varphi^-),
\\ e^{-\xi G^-_{-1/2}(x, \varphi^+, \varphi^-)} \cdot (x,\varphi^+,\varphi^-)&=& (x + \varphi^+ \xi, \varphi^+, \xi + \varphi^-),
\\ e^{-\xi G^-_{1/2}(x, \varphi^+, \varphi^-)} \cdot (x,\varphi^+,\varphi^-) &=& (x + \varphi^+ \xi x, \varphi^+, \xi x +  \varphi^- - \varphi^+ \varphi^- \xi), \label{group-action-last}
\end{eqnarray}
as expected. 

In general, a formal $N=2$ superprojective transformation $T(x, \varphi^+, \varphi^-) = (\tilde{x}, \tilde{\varphi}^+, \tilde{\varphi}^-)$  has the form
\begin{eqnarray}
\tilde{x} &=& \frac{ax + b}{cx + d} + \varphi^+ \frac{e^+(\gamma^- x + \delta^-)}{(cx + d)^2} + \varphi^+ \frac{(f^+ x + h^+) (\gamma^- x + \delta^-)}{(cx + d)^3}  \label{transform1-x}\\
& &  + \varphi^- \frac{e^-(\gamma^+ x + \delta^+)}{(cx + d)^2} +  \varphi^- \frac{(f^-x + h^-)(\gamma^+ x + \delta^+)}{(cx + d)^3}  \nonumber\\
& & + \varphi^+ \varphi^- \frac{2\gamma^+ \gamma^-dx  - (\gamma^+ \delta^- + \delta^+ \gamma^-)(cx - d) - 2\delta^+ \delta^-c}{(cx + d)^3} \nonumber 
\end{eqnarray}
\begin{eqnarray}
\qquad \tilde{\varphi}^+ &=& \frac{\gamma^+ x + \delta^+}{cx + d} + \varphi^+ \frac{e^+}{cx+d} +  \varphi^+ \frac{f^+ x+ h^+}{(cx+d)^2} + \varphi^+ \varphi^- \frac{\gamma^+ d - \delta^+ c}{(cx+d)^2}\\
\tilde{\varphi}^- &=& \frac{\gamma^- x + \delta^-}{cx + d} + \varphi^- \frac{e^-}{cx+d} +  \varphi^- \frac{f^- x + h^-}{(cx+d)^2} - \varphi^+ \varphi^- \frac{\gamma^- d - \delta^-c}{(cx+d)^2} \label{transform3-x}
\end{eqnarray}
for $a,b,c,d,e^\pm, f^\pm, h^\pm \in \bigwedge_\infty^0$ and $\gamma^\pm, \delta^\pm \in \bigwedge_\infty^1$ satisfying
\begin{eqnarray}
ad - bc &=& 1, \label{transform-condition1-x}\\
e^+ e^-  &=& 1- \gamma^+ \delta^- + \delta^+ \gamma^- \label{transform-condition2-x}\\
f^\pm &=& \mp  e^\pm \gamma^+ \gamma^-d \label{transform-condition3-x} \\
h^\pm &=& \pm e^\pm ( \delta^+ \delta^-c - (\gamma^+ \delta^- + \delta^+ \gamma^-)d  \mp \delta^+ \delta^- \gamma^+ \gamma^- d)  \label{transform-condition4-x}, 
\end{eqnarray} 
or equivalently of the form (\ref{short-transform1}) -- (\ref{short-transform3}) (with $(w, \rho^+, \rho^-) = (x,\varphi^+, \varphi^-)$) satisfying (\ref{transform-condition1-x}) and (\ref{transform-condition2-x}).
This can be seen by exponentiating a general element of $\mathfrak{osp}_{\mathbb{C}}(2|2)$ using the above action, or equivalently exponentiating a general infinitesimal $N=2$ superprojective transformation given by a linear combination of $L_{\pm 1} (x, \varphi^+, \varphi^-)$, $L_0(x, \varphi^+, \varphi^-)$, $J_0(x, \varphi^+, \varphi^-)$, $G^+_{\pm 1/2}(x, \varphi^+, \varphi^-)$, and $G^-_{\pm 1/2}(x, \varphi^+$, $\varphi^-)$. 

We note that this general form for the $N=2$ superprojective transformations does not agree with that presented in for instance \cite{BL},\cite{C}, \cite{Ki}, or \cite{Nogueira}  or any other sources in the literature that we have encountered.  Below we give a concrete example of an $N=2$ superprojective transformation that is not of the form given in all the descriptions we have found in the literature to date.   This includes presentations in \cite{Melzer} and \cite{Sc} in the ``nonhomogeneous" coordinate frame (see Section  \ref{nonhomo-section}).

\begin{ex}\label{ex1}
Let $A_1, A_{-1} \in \bigwedge_\infty^0$, and $M^\pm_{\frac{1}{2}} \in \bigwedge_\infty^1$.  Consider the following formal $N=2$ superprojective transformation.
\begin{equation} \label{ex-transformation}
T(x, \varphi^+, \varphi^-) \hspace{3.5in}
\end{equation}
\vspace{-.3in}
\begin{eqnarray*}
&=& e^{- A_{-1} L_{-1} (x,\varphi^+,\varphi^-)   } e^{ - \left(  A_{1} L_{1} (x,\varphi^+,\varphi^-)  + M^+_{  \frac{1}{2}} G^+_{  \frac{1}{2}}(x,\varphi^+,\varphi^-)   +  M^-_{ \frac{1}{2}} G^-_{  \frac{1}{2}}(x,\varphi^+,\varphi^-)   \right) } \cdot \\
& & \ \ (x, \varphi^+, \varphi^-) \\
&=& \! \! e^{- A_{-1} L_{-1} (x,\varphi^+,\varphi^-)   } e^{ - \left(  M^+_{  \frac{1}{2}} G^+_{  \frac{1}{2}}(x,\varphi^+,\varphi^-)   +  M^-_{ \frac{1}{2}} G^-_{  \frac{1}{2}}(x,\varphi^+,\varphi^-)   \right) } e^{ - A_{1} L_{1} (x,\varphi^+,\varphi^-)  }  \cdot  \\
& & \ \  (x, \varphi^+, \varphi^-) \\
&=& \! \! e^{ A_{-1}   \frac{\partial}{\partial x}    }  e^{ \left(  M^+_{  \frac{1}{2}} \Bigl( x \Bigl( \frac{\partial}{\partial \varphi^+} - \varphi^- \frac{\partial}{\partial x} \Bigr)  +  \varphi^+ \varphi^- \frac{\partial}{\partial \varphi^+} \Bigr)  +  M^-_{ \frac{1}{2}} \Bigl( x \Bigl( \frac{\partial}{\partial \varphi^-} - \varphi^+ \frac{\partial}{\partial x} \Bigr) - \varphi^+ \varphi^- \frac{\partial}{\partial \varphi^-} \Bigr) \right)}\cdot \\
& & \ \  e^{ A_{1} \Bigl( x^2 \frac{\partial}{\partial x} +  x \Bigl( \varphi^+ \frac{\partial}{\partial \varphi^+} + \varphi^- \frac{\partial}{\partial \varphi^-} \Bigr) \Bigr) }   \cdot (x, \varphi^+, \varphi^-)  \\
&=&  \! \!  \biggl(\frac{x + A_{-1}}{1- A_1(x + A_{-1} )} + \varphi^+ \frac{M^-_\frac{1}{2}  (x +  A_{-1})  }{1 -A_1(x + A_{-1})}   + \varphi^+ \frac{A_1M^-_\frac{1}{2}  (x +  A_{-1})^2 }{(1 -A_1(x + A_{-1}) )^2}  \\
& & \ \ + \varphi^- \frac{M^+_\frac{1}{2}  (x +   A_{-1})  }{1 -A_1(x + A_{-1}) } + \varphi^- \frac{A_1 M^+_\frac{1}{2}  (x +   A_{-1})^2  }{(1 -A_1(x + A_{-1}) )^2} \\
& & \ \ + \varphi^+ \varphi^- \frac{2M^+_\frac{1}{2} M^-_\frac{1}{2} ( x  + A_{-1} )}{1-A_1(x +  A_{-1}) } + \varphi^+ \varphi^- \frac{4A_1M^+_\frac{1}{2} M^-_\frac{1}{2} ( x  + A_{-1} )^2}{(1-A_1(x +  A_{-1}) )^2} \\
& & \ \ + \varphi^+ \varphi^- \frac{2A_1^2M^+_\frac{1}{2} M^-_\frac{1}{2} ( x  + A_{-1} )^3}{(1-A_1(x +  A_{-1}) )^3} ,  \ \frac{M^+_\frac{1}{2} ( x +  A_{-1} )}{1-A_1(x +A_{-1}) }  
\end{eqnarray*}
\begin{eqnarray*}
& &  \ \ + \varphi^+ \frac{1-M^+_\frac{1}{2}  M^-_\frac{1}{2} (x +   A_{-1}) }{1-A_1(x+A_{-1}) } -  \varphi^+ \frac{ A_1 M^+_\frac{1}{2}  M^-_\frac{1}{2} (x +   A_{-1})^2 }{(1-A_1(x+A_{-1}) )^2}  \\
& & \ \ +  \varphi^+ \varphi^- \frac{M^+_\frac{1}{2} }{1-A_1(x+A_{-1} )} + \varphi^+ \varphi^- \frac{A_1 M^+_\frac{1}{2} (x + A_{-1})}{(1-A_1(x+A_{-1}) )^2}, \\
& & \ \      \frac{M^-_\frac{1}{2}  (x +  A_{-1} )}{1-A_1(x + A_{-1}) } + \varphi^- \frac{1+  M^+_\frac{1}{2}  M^-_\frac{1}{2}(x +   A_{-1}) }{1-A_1(x+A_{-1}) }    \\
& & \ \ +  \varphi^- \frac{ A_1 M^+_\frac{1}{2}  M^-_\frac{1}{2}(x +   A_{-1})^2 }{(1-A_1(x+ A_{-1}) )^2} - \varphi^+ \varphi^- \frac{M^-_\frac{1}{2}  }{1-A_1(x+A_{-1}) } \\
& & \ \  - \varphi^+ \varphi^- \frac{A_1 M^-_\frac{1}{2}(x + A_{-1} )  }{(1-A_1(x+A_{-1}) )^2} \biggr)  \\
&=&  \! \!  \biggl(\frac{x + A_{-1}}{-A_1x + 1- A_1 A_{-1} } + \varphi^+ \frac{M^-_\frac{1}{2}  (x +  A_{-1})  }{(-A_1x + 1- A_1 A_{-1} )^2} \\
& & \ \ + \varphi^- \frac{M^+_\frac{1}{2}  (x +   A_{-1})  }{(-A_1x + 1- A_1 A_{-1} )^2} + \varphi^+ \varphi^- \frac{2M^+_\frac{1}{2} M^-_\frac{1}{2} ( x  + A_{-1} )}{(-A_1x + 1- A_1 A_{-1} )^3},  \\
& &  \ \ \frac{M^+_\frac{1}{2} ( x +  A_{-1} )}{-A_1x + 1- A_1 A_{-1} } + \varphi^+ \frac{1}{-A_1x+1- A_1 A_{-1} } +  \varphi^+ \frac{ -  M^+_\frac{1}{2}  M^-_\frac{1}{2} (x +   A_{-1}) }{(-A_1x+1- A_1 A_{-1} )^2}  \\
& & \ \ + \varphi^+ \varphi^- \frac{M^+_\frac{1}{2} }{(-A_1x+1- A_1 A_{-1} )^2}, \     \frac{M^-_\frac{1}{2}  (x +  A_{-1} )}{-A_1x + 1- A_1 A_{-1} } \\
& & \ \ + \varphi^- \frac{1}{-A_1x+1- A_1 A_{-1} }    +  \varphi^- \frac{ M^+_\frac{1}{2}  M^-_\frac{1}{2}(x +   A_{-1}) }{(-A_1x+1- A_1 A_{-1} )^2}  \\
& & \ \ - \varphi^+ \varphi^- \frac{M^-_\frac{1}{2}  }{(-A_1x+1- A_1 A_{-1} )^2} \biggr) 
\end{eqnarray*}
which we need to rewrite as
\begin{equation}\label{ex-transformation2}
T(x, \varphi^+ , \varphi^-) \hspace{3.5in}
\end{equation}
\begin{eqnarray*}
\lefteqn{ = \biggl(\frac{x + A_{-1}}{-A_1x + 1- A_1 A_{-1} } + \varphi^+ \frac{M^-_\frac{1}{2}  (x +  A_{-1})  }{(-A_1x + 1- A_1 A_{-1} )^2}   + \varphi^- \frac{M^+_\frac{1}{2}  (x +   A_{-1})  }{(-A_1x + 1- A_1 A_{-1} )^2}   }\\
& &\! \! \! \! + \varphi^+ \varphi^- \frac{2M^+_\frac{1}{2} M^-_\frac{1}{2} ( x  + A_{-1} )}{(-A_1x + 1- A_1 A_{-1} )^3},  \  \frac{M^+_\frac{1}{2} ( x +  A_{-1} )}{-A_1x + 1- A_1 A_{-1} } \nonumber\\
& & \! \!  \! \! + \varphi^+ \frac{1 + M^+_\frac{1}{2} M^-_\frac{1}{2} A_{-1}}{-A_1x+1- A_1 A_{-1} }  +  \varphi^+ \frac{ M^+_\frac{1}{2}  M^-_\frac{1}{2} (-  ( 1 -  A_1 A_{-1} ) x  - 2 A_{-1} + A_1 A_{-1}^2 )}{(-A_1x+1- A_1 A_{-1} )^2}  \nonumber\\
& & \! \! \! \! + \varphi^+ \varphi^- \frac{M^+_\frac{1}{2} }{(-A_1x+1- A_1 A_{-1} )^2}, \   \frac{M^-_\frac{1}{2}  (x +  A_{-1} )}{-A_1x + 1- A_1 A_{-1} } + \varphi^- \frac{1 - M^+_\frac{1}{2} M^-_\frac{1}{2} A_{-1}}{-A_1x+1- A_1 A_{-1} }   \nonumber \\ 
& & \! \! \! \! +  \varphi^- \frac{ M^+_\frac{1}{2}  M^-_\frac{1}{2} (  (1- A_1 A_{-1} ) x  + 2A_{-1} - A_1 A_{-1}^2 )}{(-A_1x+1- A_1 A_{-1} )^2}  - \varphi^+ \varphi^- \frac{M^-_\frac{1}{2}  }{(-A_1x+1- A_1 A_{-1} )^2} \biggr) \nonumber
\end{eqnarray*}
to be of the form (\ref{transform1-x}) -- (\ref{transform3-x}) satisfying (\ref{transform-condition1-x}) -- (\ref{transform-condition4-x}), where here $a = 1$, $b = A_{-1}$, $c = -A_1$, $d = 1- A_1A_{-1}$, $e^\pm = 1 \pm M^+_\frac{1}{2} M^-_\frac{1}{2} A_{-1}$, $\gamma^\pm = M^\pm_\frac{1}{2}$, and $\delta^\pm = M^\pm_\frac{1}{2} A_{-1}$.  (And thus by (\ref{transform-condition3-x}) and (\ref{transform-condition4-x}), we have $f^\pm =  \mp  M^+_\frac{1}{2}M^-_\frac{1}{2}(1- A_1A_{-1})$ and   $h^\pm  =  \pm  ( - M^+_\frac{1}{2} M^-_\frac{1}{2} A_{-1}^2 A_1 - 2M^+_\frac{1}{2} M^-_\frac{1}{2} A_{-1} ( 1- A_1A_{-1}) )$.)  

Note that the calculation of $T$ above can either be done by direct expansion of the exponentials or more easily by observing that letting
\begin{eqnarray}
T_1 (x, \varphi^+, \varphi^-) &=& e^{-A_1L_1(x, \varphi^+, \varphi^-)} \cdot (x, \varphi^+, \varphi^-)\\
\quad T_2 (x, \varphi^+, \varphi^-) &=& e^{-\left(M^+_\frac{1}{2} G^+_\frac{1}{2} (x, \varphi^+, \varphi^-)+ M^-_\frac{1}{2} G^-_\frac{1}{2} (x, \varphi^+, \varphi^-)\right)}\cdot (x,\varphi^+, \varphi^-) \\
T_3 (x, \varphi^+, \varphi^-) &=& e^{-A_{-1} L_{-1}(x, \varphi^+, \varphi^-)} \cdot (x, \varphi^+, \varphi^-),
\end{eqnarray}
then by (\ref{group-action1}), direct expansion, and (\ref{group-action3}), respectively, we have
\begin{eqnarray}
T_1 (x, \varphi^+, \varphi^-) &=& \left( \frac{x}{1 - A_1 x}, \varphi^+ \frac{1}{1-A_1x}, \varphi^- \frac{1}{1 - A_1x} \right) \\
T_2 (x, \varphi^+, \varphi^-) &=& \Bigl(x + \varphi^+ M^-_\frac{1}{2} x + \varphi^- M^+_\frac{1}{2} x + 2 \varphi^+ \varphi^- M^+_\frac{1}{2} M^-_\frac{1}{2} x, \label{T2-for-example} \\
& & M^+_\frac{1}{2} x + \varphi^+  - \varphi^+ M^+_\frac{1}{2} M^-_\frac{1}{2} x + \varphi^+ \varphi^- M^+_\frac{1}{2}, \nonumber \\
& &  M^-_\frac{1}{2} x + \varphi^- + \varphi^- M^+_\frac{1}{2} M^-_\frac{1}{2} x - \varphi^+ \varphi^- M^-_\frac{1}{2} \Bigr) \nonumber \\
T_3 (x, \varphi^+, \varphi^-) &=& (x + A_{-1}, \varphi^+, \varphi^-).
\end{eqnarray}
By Propositions \ref{Switch} and \ref{Switch2}, we have $T(x, \varphi^+, \varphi^-) = T_1 \circ T_2 \circ T_3 (x, \varphi^+, \varphi^-)$ which gives (\ref{ex-transformation}).
\end{ex}

Note that if $d \in (\bigwedge_\infty^0)^\times$, that is if $d$ is invertible, then 
\begin{equation}
\frac{f^\pm x + h^\pm }{(cx+d)^2} = \frac{(f^\pm - \frac{h^\pm c}{d})x + \frac{h^\pm }{d} (cx + d)}{(cx+d)^2} = \frac{(f^\pm - \frac{h^\pm c}{d})x }{(cx+d)^2} + \frac{ \frac{h^\pm }{d}}{cx+d} ,
\end{equation}
and the $g^\pm(x)$ term (using the notation (\ref{formal-superconformal-condition1}) -- (\ref{formal-superconformal-condition3}) for a superconformal function) i.e., the coefficient of $\varphi^\pm$ in $\tilde{\varphi}^\pm$ in the $N=2$ superprojective transformation can be written 
\begin{eqnarray}
\lefteqn{\frac{e^\pm}{cx+d} \mp \frac{e^\pm  ( \gamma^+ \gamma^-d x - \delta^+ \delta^-c  + (\gamma^+ \delta^- + \delta^+ \gamma^-)d \pm \delta^+ \delta^- \gamma^+ \gamma^- d)}{(cx+d)^2} } \\
&=& \frac{e^\pm(1 \pm (  \delta^+ \delta^-\frac{c}{d}  - (\gamma^+ \delta^- + \delta^+ \gamma^- ) \mp \delta^+ \delta^- \gamma^+ \gamma^- ))}{cx+d} \nonumber\\
& & \quad \mp e^\pm \frac{(    \gamma^+ \gamma^-d   +  \delta^+ \delta^-\frac{c^2}{d}  - (\gamma^+ \delta^- + \delta^+ \gamma^-)c  \mp \delta^+ \delta^- \gamma^+ \gamma^-c)x}{(cx+d)^2} \nonumber \\
&=& \frac{e^\pm(1 \pm (  \delta^+ \delta^-\frac{c}{d}  - (\gamma^+ \delta^- + \delta^+ \gamma^- ) \mp \delta^+ \delta^- \gamma^+ \gamma^- ))}{cx+d} \nonumber\\
& & \quad \mp e^\pm \frac{(   (\gamma^+ d   -  \delta^+ c)(\gamma^- d - \delta^- c) \mp \delta^+ \delta^- \gamma^+ \gamma^-cd)x}{d(cx+d)^2} .\nonumber 
\end{eqnarray}

In \cite{Nogueira}, P. Nogueira claims that all $N=2$ superprojective transformations are $N=2$ superconformal functions of the form (\ref{formal-superconformal-condition1}) -- (\ref{formal-superconformal-condition3}) with 
\begin{eqnarray*}
f(x) \! \! \! \! &=& \! \! \! \!  \frac{ax +b}{cx+d}, \qquad \quad \psi^\pm(x)  =   \frac{\gamma^\pm x + \delta^\pm}{cx+d}, \\
g^\pm(x) \! \! \! \! &=& \! \! \! \! e^{\mp i \alpha}\!  \left(\!  \frac{1 + \frac{1}{2} (\delta^+ \gamma^- + \delta^- \gamma^+) - \frac{1}{4} \delta^+ \delta^- \gamma^+ \gamma^- }{cx+d} \mp \frac{(\gamma^+ d - \delta^+c)(\gamma^- d - \delta^- c) x}{d(cx+d)^2} \! \right)
\end{eqnarray*}
with $a,b,c,d, \alpha \in \bigwedge_\infty^0$, and $\gamma^\pm, \delta^\pm \in \bigwedge_\infty^1$ satisfying $ad-bc = 1$.  Setting
\begin{eqnarray*}
e^\pm \!\! \! &=& \! \! \! e^{\mp i \alpha} \frac{1 + \frac{1}{2} (\delta^+ \gamma^- + \delta^- \gamma^+) - \frac{1}{4} \delta^+ \delta^- \gamma^+ \gamma^- }{1 \pm (\delta^+ \delta^-\frac{c}{d} - (\gamma^+ \delta^- + \delta^+ \gamma^-)  \mp \delta^+ \delta^- \gamma^+ \gamma^-   )} \\
&=& \! \! \! e^{\mp i \alpha} (1 + \frac{1}{2} (\delta^+ \gamma^- + \delta^- \gamma^+) - \frac{5}{4} \delta^+ \delta^- \gamma^+ \gamma^-  \mp ( \delta^+ \delta^-\frac{c}{d} - (\gamma^+ \delta^- + \delta^+ \gamma^-) ) \nonumber,
\end{eqnarray*}
i.e.,
\[e^{\mp i \alpha}  = e^\pm  (1 - \frac{1}{2} (\delta^+ \gamma^- + \delta^- \gamma^+) - \frac{1}{4} \delta^+ \delta^- \gamma^+ \gamma^-  \pm ( \delta^+ \delta^-\frac{c}{d} - (\gamma^+ \delta^- + \delta^+ \gamma^-) )), \]
it is easy to see that Nogueira's $N=2$ superprojective transformations are equivalent to ours if and only if $d$ is invertible.    Since in general though, one can not assume $d$ is invertible (even if it is nonzero), Nogueira's $N=2$ superprojective transformations do not exhaust the entire group of $N=2$ superprojective transformations and in fact do not even form a subgroup.   In Example \ref{ex1} above, letting $A_1, A_{-1} \in \bigwedge_\infty^0$ such that $(A_1A_{-1})_B = 1$ gives a concrete counter example to Nogueira's claim of providing a presentation of the group of $N=2$ superprojective transformations.  In addition, it is not possible to deduce the general form of an $N=2$ superprojective transformation {}from Nogueira's subset of $N=2$ superprojective transformations.  A similar situation occurs if one assumes $c$ is invertible. 

Presentations of the group of $N=2$ superprojective transformations such as \cite{C} and  \cite{Ki},  miss crucial terms such as the $f^\pm$ or $h^\pm$. The $N=2$ superprojective transformation $T$ given in Example \ref{ex1} gives a counterexample to the claims of Cohn in \cite{C} and Kiritsis in \cite{Ki} of giving the form of a general $N=2$ superprojective transformation.  In fact the transformation $T_2$ given in (\ref{T2-for-example}) is itself a counter example to Cohn's and Kiritsis' transformations as is the $N=2$ superprojective transformation given by (\ref{symmetric-action-transformation}) which arises {}from putting an action of the symmetric group on $n \in \mathbb{N}$ letters on the moduli space of $N=2$ super-Riemann spheres with $n$ incoming punctures.   In \cite{BL}, it is easy to see that the transformations given are not $N=2$ superconformal in our sense, i.e., they do not in general satisfy conditions (\ref{nice-superconformal-condition}).  

In \cite{Melzer} and \cite{Sc}, $N=2$ superprojective transformations are presented for the ``nonhomogeneous" coordinate system.  In Section \ref{nonhomo-section}, we discuss this coordinate system and show that the presentations of $N=2$ superprojective transformations given in \cite{Melzer} and \cite{Sc} are not the $N=2$ superconformal automorphisms of the super-Riemann sphere in the sense studied in this paper.

\section[A reformulation of the moduli space]{A reformulation of the moduli space of $N=2$ super-Riemann spheres with tubes}\label{reformulate-moduli-section}

In this section, using the characterization of local $N=2$ superconformal coordinates in terms of 
exponentials of certain infinite sums of superderivations proved in  Section \ref{infinitesimals-section}, we show that a canonical $N=2$ super-Riemann sphere with tubes can be identified with certain data concerning the punctures and the coefficients of the infinite sums of superderivations appearing in  
the expressions for the local coordinates.  Thus we can identify the moduli space of $N=2$ super-Riemann spheres with tubes with the set of this data.  
 
Let
\begin{multline*}
\mathcal{H}  =  \bigl\{ (A^+, A^-,M^+, M^-) \in \mbox{$(\bigwedge_\infty^\infty)^2$} 
\; | \; \tilde{E}(A^+, A^-,M^+, M^-) \hspace{.1cm} \mbox{is an absolutely } \bigr. \\
\bigl. \mbox{ convergent power series in some neighborhood of } 0 \bigr\},
\end{multline*} 
and for $n \in \Z$, let
\begin{multline*} 
S^2M^{n - 1} = \bigl\{ \bigl((z_1, \theta_1^+, \theta_1^-),\dots,(z_{n-1}, \theta_{n-1}^+, \theta_{n-1}^-)\bigr) \; | \; \\
(z_j, \theta_j^+, \theta_j^-) \in \mbox{$(\bigwedge_\infty^0)^\times \times (\bigwedge_\infty^1)^2$}, 
 \; (z_j)_B \neq (z_k)_B , \; \mbox{for} \; j \neq k \bigr\} . 
\end{multline*} 
Note that for $n=1$, the set $S^2M^0$ has exactly one element.

For any canonical $N=2$ supersphere with $1 + n$ tubes, we can write the power series expansion of the local coordinate at the $k$-th puncture $(z_k,\theta_k^+, \theta_k^-)$, for $k = 1,\dots,n$, as 
\begin{eqnarray*}
\lefteqn{H_k (w, \rho^+, \rho^-)}\\
& \! \! =& \! \!  \! \exp \Biggl(\! - \! \sum_{j \in \Z} \Bigl( A^{+,(k)}_j L_j(x,\varphi^+, \varphi^-) + A^{-,(k)}_j J_j(x,\varphi^+, \varphi^-) \\ 
& &  +  M^{+,(k)}_{j - \frac{1}{2}} G^+_{j -\frac{1}{2}} (x, \varphi^+, \varphi^-) +  M^{-,(k)}_{j - \frac{1}{2}} G^-_{j -\frac{1}{2}} (x, \varphi^+, \varphi^-) \Bigr) \! \Biggr)\!  \cdot  \! (a_0^{+,(k)})^{-2L_0(x,\varphi^+, \varphi^-)} \cdot  \\
& &    \left.  (a_0^{-,(k)})^{-J_0(x,\varphi^+, \varphi^-)} \cdot (x, \varphi^+, \varphi^-) \right|_{(x,\varphi^+, \varphi^-) = (w - z_k - \rho^+ \theta^-_k - \rho^- \theta^+_k, \rho^+ - \theta^+_k,  \rho^- - \theta^-_k)}\\
&\! \! =&\! \!  \! \hat{E}(a_0^{+,(k)}, a_0^{-,(k)},A^{+,(k)}, A^{-,(k)}, M^{+,(k)}, M^{-,(k)}) (w - z_k - \rho^+ \theta^-_k - \rho^- \theta^+_k, \\
& & \rho^+ - \theta^+_k,  \rho^- - \theta^-_k) 
\end{eqnarray*}
for $(a_0^{+,(k)}, a_0^{-,(k)},A^{+,(k)}, A^{-,(k)}, M^{+,(k)}, M^{-,(k)}) \in ((\bigwedge_\infty^0)^\times)^2 /\langle \pm 1 \rangle \times \mathcal{H}$, where $(z_n$, $\theta_n^+, \theta_n^-) = 0$, and we can write the power series expansion of the local coordinate at $\infty$ as  
\begin{eqnarray*}
\lefteqn{H_0 (w, \rho^+, \rho^-)} \\
&=& \! \! \! \exp \Biggl(\sum_{j \in \Z} \Bigl( A^{+,(0)}_j L_{-j}(w,\rho^+, \rho^-) + A^{-,(0)}_j J_{-j}(w,\rho^+, \rho^-) \\
& & \  +  M^{+,(0)}_{j - \frac{1}{2}} G^+_{-j +\frac{1}{2}} (w, \rho^+, \rho^-) +  M^{-,(0)}_{j - \frac{1}{2}} G^-_{-j +\frac{1}{2}} (w, \rho^+, \rho^-)
\Bigr) \! \Biggr) \cdot \Bigl(\frac{1}{w}, \frac{i \rho^+}{w}, \frac{i \rho^-}{w}\Bigr)  \\ 
&=& \! \! \!  \tilde{E}(A^{+,(0)}, -A^{-,(0)}, -iM^{+,(0)}, -iM^{-,(0)}) \Bigl(\frac{1}{w}, \frac{i \rho^+}{w}, \frac{i \rho^-}{w} \Bigr)
\end{eqnarray*}
for $(A^{+,(0)}, -A^{-,(0)}, -iM^{+,(0)}, -iM^{-,(0)}) \in \mathcal{H}$.  But $(A^{+,(0)},- A^{-,(0)}, -iM^{+,(0)},$ $-iM^{-,(0)}) \in \mathcal{H}$ if and only if $(A^{+,(0)}, A^{-,(0)}, M^{+,(0)}, M^{-,(0)}) \in \mathcal{H}$.
Thus we have the following two theorems which follow {}from Theorem \ref{bijection} and the characterization of local coordinates in terms of infinitesimals given above.

\begin{thm}\label{moduli2}
The moduli space of $N=2$ super-Riemann spheres with $1 + n$ tubes, for $n \in \Z$, can be identified with the set   
\begin{equation}
S^2K(n) = S^2M^{n-1} \times \mathcal{H} \times \bigl(((\mbox{$\bigwedge_\infty^0$})^\times)^2 /\langle \pm 1 \rangle \times \mathcal{H}\bigr)^n .
\end{equation}
\end{thm}

For $N=2$ super-Riemann spheres with one tube, {}from  Theorem \ref{bijection}  and equations (\ref{system-of-equations3}) -- (\ref{system-of-equations5}), we have:
 
\begin{thm}\label{moduli1}
The moduli space of $N=2$ super-Riemann spheres with one tube can be identified with the set
\[S^2K(0) = \bigl\{(A^+, A^-,M^+, M^-) \in \mathcal{H} \; | \; A_1^+ = M^+_{\frac{1}{2}} =M^-_{\frac{1}{2}} = 0 \bigr\} .\] 
\end{thm}

By Theorems \ref{moduli2} and \ref{moduli1}, we can identify $S^2K(n)$ with the moduli space of $N=2$ super-Riemann spheres with $1 + n$ tubes for $n \in \mathbb{N}$, and the set 
\[S^2K = \bigcup_{n \in \mathbb{N}} S^2K(n) \]
can be identified with the moduli space of $N=2$ super-Riemann spheres with tubes.  The actual elements of $S^2K$ give the data for a canonical $N=2$ supersphere representative of a given equivalence class of $N=2$ super-Riemann spheres with tubes modulo $N=2$ superconformal equivalence. {}From now on it will be convenient to refer to $S^2K$ as the moduli space of $N=2$ super-Riemann spheres with tubes.  Any element of $S^2K(n)$, for $n \in \Z$, can be written as 
\begin{multline} 
\Bigl((z_1, \theta_1^+, \theta_1^-),\dots,(z_{n-1}, \theta_{n-1}^+, \theta_{n-1}^-); (A^{+,(0)},A^{-,(0)}, M^{+,(0)}, M^{-,(0)}), \\
(a_0^{+,(1)}, a_0^{-,(1)}, A^{+,(1)}, A^{-,(1)}, M^{+,(1)}, M^{-,(1)}), \dots,\\
(a_0^{+,(n)} a_0^{-,(n)}, A^{+,(n)}, A^{-,(n)}, M^{+,(n)},M^{-,(n)}  )\Bigr)  
\end{multline}
where $((z_1, \theta_1^+, \theta_1^-),\dots,(z_{n-1}, \theta_{n-1}^+, \theta_{n-1}^-)) \in S^2M^{n-1}$,  $(a_0^{+,(k)}, a_0^{-,(k)}) \in (\mbox{$\bigwedge_\infty^0$})^\times)^2 /$ $\langle \pm 1 \rangle$ for $k = 1,\dots, n$, and $( A^{+,(k)}, A^{-,(k)}, M^{+,(k)}, M^{-,(k)}) \in \mathcal{H}$ for $k = 0,\dots, n$.  Thus for an element $Q \in S^2K$, we can think of $Q$ as consisting of the above data, or as being a canonical $N=2$ supersphere with tubes corresponding to that data.

The element of $\mathcal{H}$ with all components equal to $0$ will be denoted by $(\mathbf{0},\mathbf{0}, \mathbf{0},\mathbf{0})$ or just $\mathbf{0}$.  Note that in terms of the chart $(U_\sou, \sou)$ of $S^2\hat{\mathbb{C}}$, the local coordinate chart corresponding to $(1,1, \mathbf{0},\mathbf{0}, \mathbf{0},\mathbf{0}) \in ((\bigwedge_\infty^0)^\times)^2/\langle \pm 1 \rangle \times \mathcal{H}$ is the identity map on $\bigwedge_\infty^0 \times (\bigwedge_\infty^1)^2$ if the puncture is at $0$, is the shift $s_{(z_k,\theta_k^+, \theta_k^-)}(w,\rho^+, \rho^-) =  (w - z_k -\rho^+ \theta_k^- - \rho^- \theta_k^+, \rho^+ - \theta_k^+, \rho^- - \theta_k^-)$ if the puncture is at $(z_k,\theta_k^+, \theta_k^-)$, and is $I(w,\rho^+, \rho^-) = (1/w,i\rho^+/w, i\rho^-/w)$ if the puncture is at $\infty$.  We call such coordinates {\it standard local coordinates} in analogy to the nonsuper case \cite{H-book} and $N=1$ super case \cite{B-memoir}.

\section{Generalized $N=2$ super-Riemann spheres with tubes and group structures on subsets of the moduli space of $N=2$ super-Riemann spheres with one incoming tube}\label{moduli-group-section}

Let 
\[\overline{S^2K}(n) = S^2M^{n-1} \times \mbox{$(\bigwedge^\infty_\infty)^2/\langle \pm 1 \rangle$} \times
((\mbox{$\bigwedge_\infty^0$})^\times \times \mbox{$\bigwedge^\infty_\infty$})^{2n} , \]
for $n \in \Z$, and 
\[\overline{S^2K}(0) = \bigl\{ \left. (A^+, A^-,M^+, M^-) \in \mbox{$(\bigwedge^\infty_\infty)^2$} \; 
\right| \;
A_1^+ = M^+_{\frac{1}{2}} = M^-_\frac{1}{2} = 0 \bigr\} .\]
Elements of $\overline{S^2K} = \bigcup_{n \in \mathbb{N}}\overline{S^2K}(n)$ are called {\it generalized} (or {\it formal}) {\it $N=2$ super-Riemann spheres with tubes}.  Note that $S^2K \subsetneq \overline{S^2K}$.

Let $Q_1 \in \overline{S^2K}(n)$, for $n \in \mathbb{Z}_+$, be given by 
\begin{multline*}
Q_1 = \Bigl((z_1, \theta_1^+, \theta_1^-),...,(z_{n-1}, \theta_{n-1}^+, \theta_{n-1}^-); (A^{+,(0)},A^{-,(0)}, M^{+,(0)}, M^{-,(0)}), \\ 
(a_0^{+,(1)}, a_0^{-,(1)}, A^{+,(1)}, A^{-,(1)}, M^{+,(1)}, M^{-,(1)}),\dots,\\
(a_0^{+,(n)},a_0^{-,(n)},  A^{+, (n)}, A^{-, (n)},M^{+, (n)}, M^{-,(n)}) \Bigr), 
\end{multline*}
and let $Q_2 = (\mathbf{0},(b_0^{+,(1)}, b_0^{-,(1)}, B^{+,(1)}, B^{-,(1)}, N^{+,(1)}, N^{-,(1)}))$ be an element of $\overline{S^2K}(1)$, where $\mathbf{0}$ denotes the sequence in $\bigwedge_\infty^\infty$ consisting of all zeros.   For $1\leq k \leq n$, we define the {\it sewing together of the puncture at infinity (the $0$-th puncture) of  $Q_2$ to the $k$-th puncture of $Q_1$}, denoted $Q_1 \; _k\infty_0 \; Q_2 \in \overline{S^2K}(n)$, by 
\begin{multline}
Q_1 \; _k\infty_0 \; Q_2 =  
\Bigl((z_1, \theta_1^+, \theta_1^-),...,(z_{n-1}, \theta_{n-1}^+, \theta_{n-1}^-); (A^{+,(0)},A^{-,(0)}, M^{+,(0)}, \\ 
M^{-,(0)}), (a_0^{+,(1)}, a_0^{-,(1)}, A^{+,(1)}, A^{-,(1)}, M^{+,(1)}, M^{-,(1)}),\dots,  \\
 (a_0^{+,(k-1)}, a_0^{-,(k-1)}, A^{+,(k-1)}, A^{-,(k-1)}, M^{+,(k-1)}, M^{-,(k-1)}), (a_0^{+,(k)},  \\
a_0^{-,(k)},  A^{+,(k)}, A^{-,(k)}, M^{+,(k)}, M^{-,(k)}) \circ (b_0^{+,(1)}, b_0^{-,(1)}, B^{+,(1)}, B^{-,(1)}, \\
N^{+,(1)}, N^{-,(1)}), (a_0^{+,(k+1)}, a_0^{-,(k+1)}, A^{+,(k+1)}, A^{-,(k+1)},M^{+,(k+1)},  \\
M^{-,(k+1)}),\dots,  (a_0^{+,(n)},a_0^{-,(n)},  A^{+, (n)}, A^{-, (n)},M^{+, (n)}, M^{-,(n)}) \Bigr),
\end{multline}
where the operation $\circ$ on $((\bigwedge_\infty^0)^\times)^2/\langle \pm 1 \rangle \times (\bigwedge_\infty^\infty)^2$ is defined  by (\ref{formal-composition-for-R}) with $R = \bigwedge_\infty$.

Similarly, for $Q_3 = ((B^{+,(0)}, B^{-,(0)}, N^{+,(0)}, N^{-,(0)}),(1,1,\mathbf{0})) \in
\overline{S^2K}(1)$, we define the {\it sewing together of the puncture at infinity (the $0$-th puncture) of $Q_1$ to the $1$-st puncture of $Q_3$}, denoted $Q_3 \; _1\infty_0 \; Q_1\in \overline{S^2K}(n)$, by 
\begin{multline}
Q_3 \; _1\infty_0 \; Q_1 = 
\Bigl((z_1, \theta_1^+, \theta_1^-),...,(z_{n-1}, \theta_{n-1}^+, \theta_{n-1}^-); (A^{+,(0)},A^{-,(0)}, M^{+,(0)}, \\
M^{-,(0)}) \circ_\infty (B^{+,(0)}, B^{-,(0)}, N^{+,(0)}, N^{-,(0)}), (a_0^{+,(1)}, a_0^{-,(1)}, A^{+,(1)}, A^{-,(1)}, \\
M^{+,(1)}, M^{-,(1)}),\dots,(a_0^{+,(n)},a_0^{-,(n)},  A^{+, (n)}, A^{-, (n)},M^{+, (n)}, M^{-,(n)}) \Bigr),
\end{multline}
where the operation $\circ_\infty$ on $( \bigwedge_\infty^\infty)^2$ is defined by (\ref{sequencecompositiondef-infty}) and (\ref{sequencecompositiondef-infty2})  for $R = \bigwedge_\infty$.

\begin{prop}
The subsets of $\overline{S^2K}(1)$ given by
\[\overline{S^2K}(1)|_{0} = \{ (\mathbf{0},(a_0^+, a_0^-, A^+, A^-,  M^+, M^-)) \in \overline{S^2K}(1) \} \]
and
\[\overline{S^2K}(1) |_\infty = \{ ((A^+,A^-, M^+, M^-), (1,1,\mathbf{0}) ) \in \overline{S^2K}(1) \} \]
with the sewing operation $\; _1\infty_0 \;$ as defined above are groups.   The group $\overline{S^2K}(1)|_0$ is isomorphic to the group $((\bigwedge_\infty^0)^\times )^2/\langle \pm 1 \rangle \times (\bigwedge_\infty^\infty)^2$ discussed in Proposition \ref{sequencecompositionprop} and Remark \ref{group-remark} with $R = \bigwedge_\infty$.  The subgroup  of $\overline{S^2K}(1)|_0$ given by 
\begin{equation}\label{generalized-at-zero-subgroup}
\{(\mathbf{0},(1,1, A^+, A^-,M^+, M^-)) \in \overline{S^2K}(1)|_0 \}\
\end{equation}
is a subgroup isomorphic to the subgroup $(\bigwedge_\infty^\infty)^2$ of $((\bigwedge_\infty^0)^\times)^2/\langle \pm 1 \rangle \times (\bigwedge_\infty^\infty)^2$ under the group operation $\circ$.  The group $\overline{S^2K}(1)|_\infty$ is isomorphic to the group $(\bigwedge_\infty^\infty)^2$ with group operation $\circ_\infty$ discussed in Corollary \ref{sequence-composition-cor-infty} and Remark \ref{group-remark-infty} with $R = \bigwedge_\infty$.   In addition, under these group isomorphisms, we have the subgroups of  (\ref{generalized-at-zero-subgroup}) and $\overline{S^2K}(1)|_\infty$ corresponding to those given in Corollary \ref{subgroups-corollary}. 

The subset  $S^2K(1)|_0= ((\bigwedge_\infty^0)^\times)^2/\langle \pm 1 \rangle \times \mathcal{H}$ is a subgroup of $\overline{S^2K}(1)|_0 = ((\bigwedge_\infty^0)^\times)^2$ $/\langle \pm 1 \rangle  \times (\bigwedge_\infty^\infty)^2$.  This subgroup is isomorphic to the group of local $N=2$ superconformal transformations vanishing at 0.   Similarly, the subset $S^2K(1)|_\infty = \mathcal{H}$ is a subgroup of $\overline{S^2K}(1)|_\infty = (\bigwedge_\infty^\infty)^2$.  This subgroup is isomorphic to the group of local $N=2$ superconformal transformations vanishing at $\infty$.

If $(A^+,A^-, M^+, M^-) \in (\bigwedge^\infty_\infty)^2$, $a_0^\pm \in (\bigwedge^0_\infty)^\times$, and $t \in \bigwedge_\infty^0$, then
\begin{eqnarray*}
t &\mapsto& (\mathbf{0}, ((a_0^+)^t,1, \mathbf{0}))\\
t &\mapsto& (\mathbf{0}, (1,(a_0^-)^t, \mathbf{0}))\\
t &\mapsto& (\mathbf{0}, (1, 1,t(A^+. A^-,M^+, M^-)))
\end{eqnarray*}
give homomorphisms {}from the additive group $\bigwedge_\infty^0$ to $\overline{S^2K}(1)|_0$, and 
\[t \mapsto (t(A^+, A^-,M^+, M^-),(1,1,\mathbf{0}))\]
gives a homomorphism {}from $\bigwedge_\infty^0$ to $\overline{S^2K}(1)|_\infty$. 

Let $t_1,t_2 \in \bigwedge_\infty^0$, and $(A^+, A^-,M^+, M^-) \in ( \bigwedge_\infty^\infty)^2$.
Assume 
\[t_1(A^+, A^-,M^+, M^-), t_2(A^+. A^-,M^+, M^-) \in \mathcal{H}.\]  
Then $(t_1 + t_2)(A^+, A^-,M^+, M^-) \in \mathcal{H}$, and we have
\begin{eqnarray}
\lefteqn{\quad \bigl(\mathbf{0}, (1,1,(t_1 + t_2)(A^+, A^-,M^+, M^-))\bigr) } \label{sequencesew1} \\
&=& \bigl(\mathbf{0}, (1,1,t_1(A^+, A^-,M^+, M^-)\bigr) \; _1\infty_0 \; \bigl(\mathbf{0}, (1,1,t_2(A^+, A^-,M^+, M^-)\bigr), \nonumber \\
\lefteqn{\quad \bigl((t_1+ t_2)(A^+, A^-,M^+, M^-), (1,1, \mathbf{0})\bigr) }  \label{sequencesew2}  \\
&=& \bigl(t_1(A^+, A^-,M^+, M^-),(1,1, \mathbf{0})\bigr) \; _1\infty_0 \; \bigl(t_2(A^+, A^-,M^+, M^-),(1,1, \mathbf{0})\bigr) . \nonumber
\end{eqnarray}
In particular, for $(A^+, A^-,M^+, M^-) \in \mathcal{H}$, the sets 
\begin{equation}\label{first subgroup}
\bigl\{\bigl(\mathbf{0},(1,1,t(A^+,A^-,M^+, M^-)\bigr) \; | \; t \in \mbox{$\bigwedge_\infty^0$},
\;  t(A^+, A^-,M^+,M^-) \in \mathcal{H} \bigr\}
\end{equation} 
and
\begin{equation}\label{second subgroup}
\bigl\{\bigl(t(A^+, A^-,M^+, M^-), (1,1,\mathbf{0})\bigr) \; | \; t \in \mbox{$\bigwedge_\infty^0$}, \; 
t(A^+, A^-,M^+, M^-) \in \mathcal{H} \bigr\}
\end{equation} 
are subgroups of $S^2K(1)|_0$ and $S^2K(1)|_\infty$, respectively.   
\end{prop}

\begin{proof} {}From the definition of the sewing operation defined above, it follows {}from Proposition \ref{sequencecompositionprop}  and Corollary \ref{sequence-composition-cor-infty}  that  $(\overline{S^2K}(1)|_0, \; _1\infty_0 \;)$ and $(\overline{S^2K}(1)|_\infty, \; _1\infty_0 \;)$ are groups isomorphic to $(((R^0)^\times)^2/\langle \pm 1 \rangle \times (R^\infty)^2, \circ)$ and $((R^\infty)^2, \circ_\infty)$, respectively, with $R = \bigwedge_\infty$.

From the definition of the operation $\; _1\infty_0 \;$, Proposition \ref{sequencecompositionprop} and equation (\ref{composition-with-a}), we have
\begin{multline*}
\bigl(\mathbf{0}, ((a_0^+)^{t_1}, (a_0^-)^{t_1}, t_1(A^+, A^-,M^+, M^-))) \; _1\infty_0 \; (\mathbf{0},
((a_0^+)^{t_2}, (a_0^-)^{t_2}, t_2(A^+, A^-,\\
M^+, M^-)) \bigr) 
\end{multline*}
\begin{multline*}
= \Bigl(\mathbf{0}, \bigl((a_0^+)^{t_1}(a_0^+)^{t_2}, (a_0^-)^{t_1}(a_0^-)^{t_2},  t_1(A^+,A^-,M^+, M^-) \circ t_2 \bigl\{ (a_0^+)^{t_1 2j} A_j^+, \\
(a_0^+)^{t_12j} A_j^-, (a_0^+)^{t_1(2j-1)} (a_0^-)^{-t_1} M_{j - \frac{1}{2}}^+,  (a_0^+)^{t_1(2j-1)} (a_0^-)^{t_1} M_{j - \frac{1}{2}}^- \bigr\}_{j \in \Z} \bigr)\Bigr) . 
\end{multline*}
Therefore
\begin{eqnarray}
\left(\mathbf{0}, ((a_0^+)^{t_1}, 1, \mathbf{0})) \; _1\infty_0 \; (\mathbf{0}, ((a_0^+)^{t_2},1, \mathbf{0}) \right) &=& \left(\mathbf{0}, ((a_0^+)^{t_1 + t_2}, 1, \mathbf{0})\right) \\
\left(\mathbf{0}, (1, (a_0^-)^{t_1},  \mathbf{0})) \; _1\infty_0 \; (\mathbf{0}, (1, (a_0^-)^{t_2},\mathbf{0}) \right) &=& \left(\mathbf{0}, (1, (a_0^-)^{t_1 + t_2}, \mathbf{0})\right) 
\end{eqnarray}
and 
\begin{multline}
\left(\mathbf{0}, (1, 1, t_1(A^+, A^-,M^+, M^-))) \; _1\infty_0 \; (\mathbf{0}, (1, 1, t_2(A^+, A^-,M^+, M^-)) \right) \\
=  \left(\mathbf{0}, (1, 1, (t_1 + t_2)(A^+, A^-,M^+, M^-))\right). 
\end{multline}
Thus the maps $t \mapsto (\mathbf{0}, ((a_0^+)^t,1, \mathbf{0}))$, $t \mapsto (\mathbf{0}, (1,(a_0^-)^t, \mathbf{0}))$ and $t \mapsto (\mathbf{0}, (1, 1, t(A^+, A^-,$ $ M^+, M^-))$ are homomorphisms {}from $\bigwedge_\infty^0$ to $\overline{S^2K}(1)|_0$.  In addition, by Corollary \ref{sequence-composition-cor-infty}
\begin{eqnarray*}
\lefteqn{\left(t_1(A^+, A^-,M^+, M^-),(1,1,\mathbf{0})) \; _1\infty_0 \; (t_2(A^+, A^-,M^+, M^-),(1,1,\mathbf{0})\right) }\\
&=& \left(t_1(A^+, A^-,M^+, M^-) \circ_\infty t_2(A^+, A^-,M^+, M^-),(1,1,\mathbf{0})\right) \\
&=& \left((t_1 + t_2)(A^+, A^-,M^+, M^-), (1,1,\mathbf{0})\right) .
\end{eqnarray*}
Thus $t \mapsto (t(A^+, A^-,M^+, M^-), (1,1,\mathbf{0}))$ is a homomorphism {}from $\bigwedge_\infty^0$ to $\overline{S^2K}(1)|_\infty$. 

Let
\begin{eqnarray}
H_{t_1}(x, \varphi^+, \varphi^-) &=& \tilde{E}(t_1A^+,t_!A^-,t_1M^+, t_1 M^-), \\
H_{t_2} (x, \varphi^+, \varphi^-) &=& \tilde{E}({t_2}A^+,t_2A^-,t_2M^+, t_2M^-) . 
\end{eqnarray}
Since $t_1(A^+, A^-,M^+, M^-)$, and $t_2(A^+, A^-,M^+, M^-)$ are in $\mathcal{H}$, the superfunctions
$H_{t_1}(w, \rho^+, \rho^-)$ and $H_{t_2}(w, \rho^+, \rho^-)$ are convergent in a neighborhood of
$0 \in \bigwedge_\infty^0 \oplus (\bigwedge_\infty^1)^2$.  Hence $H_{t_1} \circ H_{t_2} (w, \rho^+, \rho^-)$ and $H_{t_2} \circ H_{t_1} (w, \rho^+, \rho^-)$ are convergent in a neighborhood of $0$.  Thus
by (\ref{sequencecompositiondef}), we see that $t_1(A^+, A^-,M^+, M^-) \circ t_2(A^+, A^-,M^+,M^-)$ and $t_2(A^+, A^-,M^+, M^-) \circ t_1(A^+, A^-,M^+, M^-)$ are in $\mathcal{H}$.  By Proposition
\ref{sequencecompositionprop}, 
\begin{eqnarray*}
(t_1+  t_2)(A^+, A^-,M^+,M^-) &=& t_1(A^+, A^-,M^+, M^-) \circ t_2(A^+, A^-,M^+, M^-) \\
&=& t_2(A^+, A^-,M^+, M^-) \circ t_1(A^+. A^-,M^+, M^-) .
\end{eqnarray*}
Thus $(t_1 + t_2)(A^+, A^-,M^+,M^-) \in \mathcal{H}$. {}From the definition of the
operation $\; _1\infty_0 \;$, we see that
\begin{multline*}
\bigl(\mathbf{0}, (1,1,t_1(A^+, A^-,M^+, M^-) \circ t_2(A^+,A^-M^+, M^-))\bigr)  \\
= \bigl(\mathbf{0}, (1,1,t_1(A^+, A^-,M^+, M^-)\bigr) \; _1\infty_0 \; \bigl(\mathbf{0}, (1,1,t_2(A^+, A^-,M^+, M^-)\bigr). 
\end{multline*}
This then gives (\ref{sequencesew1}). 
                     
Similarly, let
\begin{eqnarray}
H_{t_1}\circ I^{-1} (x, \varphi^+, \varphi^-) &=& \tilde{E}(t_1A^+,- t_1 A^-, - it_1M^+, -i t_1 M^-), \\
H_{t_2} \circ I^{-1} (x, \varphi^+, \varphi^-) &=& \tilde{E}(t_2A^+,-t_2A^-,-it_2M^+, -it_2M^-) . 
\end{eqnarray}
Since $t_1(A^+, A^-,M^+, M^-)$, and $t_2(A^+, A^-,M^+, M^-)$ are in $\mathcal{H}$, the superfunctions
$H_{t_1}(w, \rho^+, \rho^-)$ and $H_{t_2}(w, \rho^+, \rho^-)$ are convergent in a neighborhood of
$\infty \in \bigwedge_\infty^0 \oplus (\bigwedge_\infty^1)^2$.  Hence $H_{t_1} \circ I^{-1} \circ H_{t_2} (w, \rho^+, \rho^-)$ and $H_{t_2} \circ I^{-1} \circ H_{t_1} (w, \rho^+, \rho^-)$ are convergent in a neighborhood of $\infty$.  Thus by (\ref{sequencecompositiondef-infty}) and (\ref{sequencecompositiondef-infty2})  we see that $t_1(A^+, A^-,M^+, M^-) \circ_\infty t_2(A^+, A^-,M^+,M^-)$ and $t_2(A^+, A^-,M^+, M^-) \circ_\infty t_1(A^+,$ $A^-,M^+, M^-)$ are in $\mathcal{H}$.  By Corollary
\ref{sequence-composition-cor-infty}, 
\begin{eqnarray*}
(t_1+  t_2)(A^+, A^-,M^+,M^-) &=& t_1(A^+, A^-,M^+, M^-) \circ_\infty t_2(A^+, A^-,M^+, M^-) \\
&=& t_2(A^+, A^-,M^+, M^-) \circ_\infty  t_1(A^+. A^-,M^+, M^-) .
\end{eqnarray*}
Thus $(t_1 + t_2)(A^+, A^-,M^+,M^-) \in \mathcal{H}$. {}From the definition of the
operation $\; _1\infty_0 \;$, we see that
\begin{multline*}
\bigl((t_1(A^+, A^-,M^+, M^-) \circ_\infty t_2(A^+,A^-M^+, M^-)), (1,1, \mathbf{0}) \bigr)  \\
= \bigl(t_1(A^+, A^-,M^+, M^-), (1,1,\mathbf{0}) \bigr) \; _1\infty_0 \; \bigl(t_2(A^+, A^-,M^+, M^-), (1,1,\mathbf{0}) \bigr). 
\end{multline*}
This then gives (\ref{sequencesew2}).
\end{proof}

\begin{rema}  {\em By the proposition above, 
\begin{eqnarray}
\qquad & \bigl\{ (\mathbf{0}, (1,1,t(A^+, A^-,M^+,M^-)) \; | \; t \in \mbox{$\bigwedge_\infty^0$},
\; (A^+, A^-,M^+, M^-) \in \mbox{$(\bigwedge_\infty^\infty)^2$}\bigr\} \label{first set}\\
\qquad & \bigl\{ (t(A^+, A^-,M^+, M^-), (1,1,\mathbf{0})) \; | \; t \in \mbox{$\bigwedge_\infty^0$},
\; (A^+, A^-,M^+, M^-) \in \mbox{$(\bigwedge_\infty^\infty)^2$}\bigr\} \label{second set}
\end{eqnarray}
are $(1,0)$-dimensional Lie supergroups over $\bigwedge_\infty$, i.e., they are abstract groups which are also $(1,0)$-dimensional superanalytic supermanifolds with superanalytic group multiplication and inverse mappings, (cf. \cite{Var}).   But elements of (\ref{first set}) and (\ref{second set}) are in general not in $S^2K(1)$, even if $(A^+, A^-,M^+, M^-) \in \mathcal{H}$.  This is one of our motivations for introducing generalized $N=2$ super-Riemann spheres with tubes.}
\end{rema} 

\begin{rema}{\em We use the notation $\; _1\infty_0 \;$ for the group operation on $\overline{S^2K}(1)|_0$ and $\overline{S^2K}(1)|_\infty$ in analogy to the notation of \cite{V}, \cite{H-book}, and \cite{B-memoir} where this notation is used to denote the sewing together of elements in the moduli space of spheres with tubes or, in the latter case, $N=1$ superspheres with tubes.  In subsequent work, we will define a sewing operation on the moduli space of $N=2$ super-Riemann spheres with tubes, and this sewing operation will coincide with the group operation $\; _1\infty_0 \;$ on $S^2K(1)|_0$ and $S^2K(1)|_\infty$ as defined above. }
\end{rema}

\section[An action of the symmetric groups on the moduli space]{An action of the symmetric 
group $S_n$ on the moduli space $S^2K(n)$}

In this section, we introduce an action of the symmetric group on $n$ letters on the moduli space of $N=2$ super-Riemann spheres with $1 + n$ tubes. 

Let $S_n$ be the permutation group on $n$ letters, for $n \in \Z$.   There is a natural action of $S_{n - 1}$ on $S^2K(n)$ defined by permuting the ordering on the first $n - 1$ positively oriented 
punctures and their local coordinates.  More explicitly, for $\sigma \in S_{n-1}$, and $Q \in S^2K(n)$ given by
\begin{multline*}
Q = \bigl((z_1, \theta_1^+, \theta_1^-),\dots,(z_{n-1}, \theta_{n-1}^+, \theta_{n-1}^-); (A^{+,(0)},A^{-,(0)}, M^{+,(0)}, M^{-,(0)}), \\
(a_0^{+,(1)},a_0^{-,(1)}, A^{+,(1)}, A^{-,(1)}, M^{+,(1)},M^{+,(1)}),\dots, \\
(a_0^{+,(n)},a_0^{-,(n)}, A^{+,(n)}, A^{-,(n)}, M^{+,(n)},M^{+,(n)})\bigr), 
\end{multline*}
we define
\begin{multline*}
\sigma(Q) = \bigl((z_{\sigma^{-1}(1)}, \theta_{\sigma^{-1}(1)}^+,  \theta_{\sigma^{-1}(1)}^-), \dots,
(z_{\sigma^{-1}(n-1)}, \theta_{\sigma^{-1}(n-1)}^+,\theta_{\sigma^{-1}(n-1)}^-); \\
(A^{+,(0)},A^{-,(0)}, M^{+,(0)}, M^{-,(0)}),
(a_0^{+,(\sigma^{-1}(1))}, a_0^{-,(\sigma^{-1}(1))}, A^{+,(\sigma^{-1}(1))}, \\
A^{-,(\sigma^{-1}(1))}, M^{+,(\sigma^{-1}(1))}, M^{-,(\sigma^{-1}(1))}), \dots, (a_0^{+,(\sigma^{-1}(n-1))}, a_0^{-,(\sigma^{-1}(n-1))}, \\
A^{+,(\sigma^{-1}(n-1))},  A^{-,(\sigma^{-1}(n-1))}, M^{+, (\sigma^{-1}(n-1))}, M^{-, (\sigma^{-1}(n-1))}), \\
(a_0^{+,(n)},a_0^{-,(n)}, A^{+,(n)}, A^{-,(n)}, M^{+,(n)},M^{+,(n)})\bigr) .  
\end{multline*}

To extend this to an action of $S_n$ on $S^2K(n)$, we first note that $S_n$ is generated by the symmetric group on the first $n-1$ letters $S_{n - 1}$ and the transposition $(n-1 \; n)$.  We can let $(n-1 \; n)$ act on $S^2K(n)$ by permuting the $(n- 1)$-th and $n$-th punctures and their local coordinates for a canonical $N=2$ supersphere with $1 + n$ tubes but the resulting $N=2$ super-Riemann sphere with $1+n$ tubes is not canonical.  To obtain the $N=2$ superconformally equivalent canonical supersphere, we have to translate the new $n$-th puncture to 0.  This translation will change the coordinates of all of the positively oriented punctures, and the local coordinates vanishing at those punctures will be shifted to vanish at the new puncture, but the infinitesimals will remain the same.  On the other hand, the translation will not change the location of the negatively oriented puncture (at infinity) but will change the local coordinate at $\infty$.  This translation is given by
\begin{eqnarray*}
T_\sou : \mbox{$\bigwedge_\infty^0$} \oplus (\mbox{$\bigwedge_\infty^1$})^2 & \longrightarrow & \mbox{$\bigwedge_\infty^0$} \oplus (\mbox{$\bigwedge_\infty^1$})^2  \\
(w,\rho^+, \rho^-) & \mapsto & (w - z_{n-1} - \rho^+\theta_{n-1}^- - \rho^- \theta_{n-1}^+, \rho^+ - \theta_{n-1}^+, \rho^- - \theta_{n-1}^-), 
\end{eqnarray*}
that is
\begin{equation}
T_\sou(w, \rho^+, \rho^-) = e^{z_{n-1} L_{-1}(w, \rho^+, \rho^-) + \theta_{n-1}^+ G_{-\frac{1}{2}}^+(w, \rho^+, \rho^-) + \theta_{n-1}^- G_{-\frac{1}{2}}^-(w, \rho^+, \rho^-) } \cdot (w, \rho^+, \rho^-),
\end{equation}
and thus by (\ref{T1}), we have $T : S^2\hat{\mathbb{C}} \rightarrow S^2\hat{\mathbb{C}}$ is given by
\begin{eqnarray*}
T(p) = \left\{
  \begin{array}{ll} 
      \sou^{-1} \circ T_\sou \circ \sou (p) & \mbox{if $p \in
           U_\sou$}, \\  
      \nor^{-1} \circ T_\nor \circ \nor (p) & \mbox{if $p \in
           U_\nor \smallsetminus \nor^{-1}(( \{(\frac{1}{z_{n-1}})_B \} \times
             (\bigwedge_\infty^0)_S ) \oplus (\bigwedge_\infty^1)^2 )$ }  
\end{array} \right.
\end{eqnarray*}
where 
\begin{eqnarray}\label{symmetric-action-transformation}
\lefteqn{\quad  T_\nor (w, \rho^+, \rho^-) }\\
&=& \biggl(\frac{w}{1 - wz_{n-1}} + \frac{i \rho^+ \theta_{n-1}^-w}{(1 - wz_{n-1})^2} + \frac{i \rho^- \theta_{n-1}^+w}{(1 - wz_{n-1})^2} - \frac{ 2\rho^+ \rho^- \theta_{n-1}^+ \theta_{n-1}^- w}{(1 - wz_{n-1})^3}, \nonumber \\
& & \quad \frac{i\theta_{n-1}^+w}{1 - wz_{n-1}} + \frac{\rho^+}{1 - wz_{n-1}} + \frac{\rho^+ \theta_{n-1}^+ \theta_{n-1}^- w}{(1-wz_{n-1})^2} + \frac{i \rho^+ \rho^- \theta_{n-1}^+}{(1-w z_{n-1})^2}  \nonumber \\
& & \quad \frac{i\theta_{n-1}^-w}{1 - wz_{n-1}} + \frac{\rho^-}{1 - wz_{n-1}} -  \frac{\rho^- \theta_{n-1}^+ \theta_{n-1}^- w}{(1-wz_{n-1})^2} - \frac{i \rho^+ \rho^- \theta_{n-1}^-}{(1-w z_{n-1})^2}  \biggr) \nonumber \\
&=& e^{-z_{n-1} L_1(w, \rho^+, \rho^-) -  i\theta_{n-1}^+ G_{\frac{1}{2}}^+(w, \rho^+, \rho^-) - i \theta_{n-1}^- G_{\frac{1}{2}}^-(w, \rho^+, \rho^-) } \cdot (w, \rho^+, \rho^-).  \nonumber
\end{eqnarray}

The new local coordinate at infinity can be written as 
\begin{equation}
\tilde{E}(\tilde{A}^{+,(0)}, - \tilde{A}^{-,(0)}, -i \tilde{M}^{+,(0)}, -i \tilde{M}^{+,(0)})  (1/w,i\rho^+/w, i\rho^-/w),
\end{equation} 
and it is determined by the old local coordinate at infinity and the superprojective transformation $T$ via
\begin{multline}
\tilde{E}( \tilde{A}^{+,(0)}, - \tilde{A}^{-,(0)}, -i \tilde{M}^{+,(0)}, -i \tilde{M}^{+,(0)}) 
\Bigl(\frac{1}{w},\frac{i\rho^+}{w}, \frac{i\rho^-}{w} \Bigr) \\
= \tilde{E}(A^{+,(0)},-A^{-,(0)}, -i M^{+,(0)}, -i M^{-,(0)}) \circ  I \circ T_\sou^{-1} (w,\rho^+, \rho^-) .
\end{multline}
Using Proposition \ref{Switch2}, we have
\begin{eqnarray}
\lefteqn{ \qquad \quad  \ \tilde{E}( \tilde{A}^{+,(0)}, - \tilde{A}^{-,(0)}, -i \tilde{M}^{+,(0)}, -i \tilde{M}^{+,(0)}) 
\Bigl(\frac{1}{w},\frac{i\rho^+}{w}, \frac{i\rho^-}{w} \Bigr) } \label{Sn-action}\\
&=& \exp \Biggl(\sum_{j \in \Z} \biggl(\tilde{A}^{+,(0)}_j L_{-j}(w,\rho^+,\rho^-) + \tilde{A}^{+,(0)}_j J_{-j}(w,\rho^+,\rho^-) \nonumber \\
& & \ \  + \tilde{M}^{+, (0)}_{j - \frac{1}{2}} G^+_{-j + \frac{1}{2}} (w,\rho^+, \rho^-) + \tilde{M}^{-, (0)}_{j - \frac{1}{2}} G^-_{-j + \frac{1}{2}} (w,\rho^+, \rho^-)  \biggr) \!\Biggr) \! \cdot \! \Bigl(\frac{1}{w},\frac{i\rho^+}{w}, \frac{i\rho^-}{w}\Bigr)  \nonumber\\
&=&  \exp \Biggl(\sum_{j \in \Z} \biggl(A^{+,(0)}_j L_{-j}(x,\varphi^+, \varphi^-) + A^{-,(0)}_j J_{-j}(x,\varphi^+, \varphi^-) \nonumber  \\
& & \ \  + M^{+,(0)}_{j - \frac{1}{2}} G^+_{-j + \frac{1}{2}}(x,\varphi^+, \varphi^-) + M^{-,(0)}_{j - \frac{1}{2}} G^-_{-j + \frac{1}{2}}(x,\varphi^+, \varphi^-) \biggr) \! \Biggr) \! \cdot \nonumber \\
& & \ \  \biggl.  \cdot \Bigl(\frac{1}{x},\frac{i\varphi^+}{x}, \frac{i\varphi^-}{x} \Bigr) 
\biggr|_{(x,\varphi^+, \varphi^-) = (w + z_{n-1} + \rho^+ \theta_{n-1}^- + \rho^- \theta_{n-1}^+, \rho^+ + \theta_{n-1}^+,  \rho^- + \theta_{n-1}^-)} \nonumber \\
&=& e^{- z_{n-1} L_{-1}(w, \rho^+, \rho^-) -  \theta_{n-1}^+ G_{-\frac{1}{2}}^+(w, \rho^+, \rho^-) - \theta_{n-1}^- G_{-\frac{1}{2}}^-(w, \rho^+, \rho^-) } \cdot  \nonumber\\
& & \ \   \exp \Biggl(\sum_{j \in \Z} \biggl(A^{+,(0)}_j L_{-j}(w,\rho^+, \rho^-) + A^{-,(0)}_j J_{-j}(w,\rho^+, \rho^-) \nonumber  \\
& & \ \  + M^{+,(0)}_{j - \frac{1}{2}} G^+_{-j + \frac{1}{2}}(w,\rho^+, \rho^-) + M^{-,(0)}_{j - \frac{1}{2}} G^-_{-j + \frac{1}{2}}(w,\rho^+, \rho^-) \biggr) \! \Biggr) \! \cdot \Bigl(\frac{1}{w},\frac{i\rho^+}{w}, \frac{i\rho^-}{w} \Bigr)  . \nonumber
\end{eqnarray} 
Then 
\begin{eqnarray}
\lefteqn{\qquad \quad \ (n-1 \; n)(Q)} \\
&=& \! \! \Bigl( \infty,(z_1, \theta^+_1,\theta_1^-),\dots,(z_{n-2}, \theta_{n-2}^+, \theta_{n-2}^-), 0, (z_{n-1},
\theta_{n-1}^+, \theta_{n-1}^-); (A^{+,(0)},A^{-,(0)}, \nonumber \\
& & \ \   M^{+, (0)}, M^{-,(0)} ), (a_0^{+,(1)},a_0^{-,(1)}, A^{+,(1)}, A^{-,(1)}, M^{+,(1)}, M^{-,(1)}), \dots \nonumber  \\
& & \ \   (a_0^{+,(n-2)}, a_0^{-,(n-2)}, A^{+,(n-2)},A^{-,(n-2)}, M^{+,(n-2)}, M^{-,(n-2)}), (a_0^{+,(n)}, \nonumber  \\
& & \ \ a_0^{-,(n)}, A^{+,(n)}, A^{-,(n)}, M^{+,(n)}, M^{-,(n)}), (a_0^{+,(n-1)}, a_0^{-,(n-1)}, A^{+,(n-1)}, \nonumber \\
& & \ \  A^{-,(n-1)}, M^{+,(n-1)}, M^{-,(n-1)}) \Bigr) \nonumber  \\ 
&=& \! \! \Bigl(\infty, (z_1 - z_{n-1} - \theta_1^+ \theta_{n-1}^- -  \theta_1^- \theta_{n-1}^+, \theta_1^+ -
\theta_{n-1}^+, \theta_1^- - \theta_{n-1}^-),(z_2 - z_{n-1}  \nonumber \\
& & \ \  - \theta_2^+\theta_{n-1}^- - \theta_2^-\theta_{n-1}^+ , \theta_2^+ - \theta_{n-1}^+,  \theta_2^- - \theta_{n-1}^- ),\dots, (z_{n-2} - z_{n-1} - \theta_{n-2}^+ \theta_{n-1}^- \nonumber \\
& & \ \  - \theta_{n-2}^- \theta_{n-1}^+  , \theta_{n-2}^+ - \theta_{n-1}^+, \theta_{n-2}^- - \theta_{n-1}^- ), (-z_{n-1},-\theta_{n-1}^+,  - \theta_{n-1}^-),0 ; (\tilde{A}^{+,(0)}, \nonumber \\
& & \quad   \tilde{A}^{-,(0)}, \tilde{M}^{+,(0)}, \tilde{M}^{-, (0)}),  (a_0^{+,(1)}, a_0^{-,(1)}, A^{+,(1)}, A^{-,(1)}, M^{+,(1)}, M^{-,(1)}), \dots, \nonumber 
\end{eqnarray}
\begin{eqnarray}
& & \ \  (a_0^{+,(n-2)}, a_0^{-,(n-2)}, A^{+,(n-2)}, A^{-,(n-2)}, M^{+,(n-2)}, M^{-,(n-2)}), (a_0^{+,(n)}, \nonumber \\
& & \ \  a_0^{-,(n)}, A^{+,(n)}, A^{-,(n)}, M^{+,(n)}, M^{-,(n)}),  (a_0^{+,(n-1)}, a_0^{-,(n-1)}, A^{+,(n-1)}, \nonumber  \\
& & \ \  A^{-,(n-1)},  M^{+,(n-1)},  M^{-,(n-1)}) \Bigr) \in S^2K(n) .  \nonumber 
\end{eqnarray}     
Thus we have an action of $S_n$ on $S^2K(n)$.

\section{The nonhomogeneous coordinate system}\label{nonhomo-section}

So far we have restricted our attention to what we call the ``homogeneous" coordinate system, denoted by $N=2$ supercoordinates $(z, \theta^+, \theta^-)$ or $(w, \rho^+, \rho^-)$ and by formal $N=2$ variables $(x, \varphi^+, \varphi^-)$.   In this section, we transfer some of our results to the ``nonhomogeneous" $N=2$ supercoordinates which we denote by $(z, \theta, \theta^*)$ or $(w, \rho, \rho^*)$ and by formal $N=2$ variables $(x, \varphi, \varphi^*)$.   This is a standard transformation in $N=2$ superconformal field theory (cf. \cite{C}, \cite{DRS}, \cite{Nogueira}).  We call these coordinate systems ``homogeneous" and ``nonhomogeneous", respectively, due to the transformation properties of a nonhomogeneous $N=2$ superconformal function on the corresponding nonhomogeneous superconformal operators $D$ and $D^*$ defined below as described in Remark \ref{transform-D-remark}, and due to the action of the $J_j$ terms, for $j \in \mathbb{Z}$,  in the algebra of infinitesimals as described in Remark \ref{homo-infinitesimals-remark}.   The transformation {}from homogeneous to nonhomogeneous coordinates is given by
\begin{equation}\label{transform-homo-inhomo}
\theta = \frac{1}{\sqrt{2}} \left( \theta^+ + \theta^- \right) \qquad \mbox{and} \qquad
\theta^* = - \frac{i}{\sqrt{2}} \left( \theta^+ - \theta^- \right) ,
\end{equation}
or equivalently
\begin{equation}
\theta^\pm  = \frac{1}{\sqrt{2}} \left( \theta \pm  i\theta^* \right).
\end{equation}
Then we have that
\begin{equation*}
\frac{\partial}{\partial \theta^\pm} =  \frac{1}{\sqrt{2}} \Bigl(\frac{\partial}{\partial \theta} \mp i \frac{\partial}{\partial \theta^*} \Bigr)
\end{equation*}
or equivalently
\begin{equation*}
\frac{\partial}{\partial \theta} =  \frac{1}{\sqrt{2}} \Bigl(\frac{\partial}{\partial \theta^+} + \frac{\partial}{\partial \theta^-} \Bigr) \qquad \mathrm{and} \qquad 
\frac{\partial}{\partial \theta^*} = \frac{i}{\sqrt{2}} \Bigl( \frac{\partial}{\partial \theta^+} - \frac{\partial}{\partial \theta^-} \Bigr) .
\end{equation*}
Thus 
\begin{equation}
D^\pm = \frac{1}{\sqrt{2}} \Bigl(\frac{\partial}{\partial \theta} \mp  i \frac{\partial}{\partial
\theta^*} \Bigr)  +  \frac{1}{\sqrt{2}} \left( \theta \mp  i\theta^* \right)  \frac{\partial}{\partial z} .\\
\end{equation}
Then we define the superderivations $D$ and $D^*$ by 
\begin{eqnarray}
D \ =&  \frac{1}{\sqrt{2}} \left( D^+ + D^- \right) & = \  \frac{\partial}{\partial \theta} +  \theta   \frac{\partial}{\partial z}  \\
D^*  =&  - \frac{i}{\sqrt{2}} \left( D^+ - D^- \right) & = \  -  \biggl( \frac{\partial}{\partial \theta^*}  + \theta^*  \frac{\partial}{\partial z} \biggr) .
\end{eqnarray}
Note that then 
\begin{eqnarray}
D^\pm &=& \frac{1}{\sqrt{2}} ( D \pm i  D^* ) \\
\left[D, D^*\right] &=&0\\
\left[D,D\right] \ = \  [D^*, D^*] &=& 2\frac{\partial}{\partial z} .
\end{eqnarray}

A homogeneous $N=2$ superconformal function $H(z, \theta^+,\theta^-)$ transforms $D^\pm$ homogeneously of degree one.    Under the transformation {}from homogeneous to nonhomogeneous 
coordinates we see that the conditions for $H(z, \theta^+, \theta^-) = (\tilde{z}, \tilde{\theta}^+, \tilde{\theta}^-)$ to transform $D^\pm$ homogeneously of degree one, i.e., that $H$ satisfies conditions (\ref{basic-superconformal-condition1}) -- (\ref{basic-superconformal-condition4}), is equivalent to $H(z, \theta, \theta^*) = (\tilde{z}, \tilde{\theta}, \tilde{\theta}^*)$ satisfying
\begin{eqnarray}
D \tilde{\theta}  + D^* \tilde{\theta}^* &=& 0 \label{inhomo-condition1} \\
D \tilde{\theta}^* -  D^* \tilde{\theta} &=& 0\\
D \tilde{z} -  \tilde{\theta} D \tilde{\theta} -  \tilde{\theta}^*D \tilde{\theta}^*  &=& 0 \\
D^* \tilde{z} -    \tilde{\theta} D^* \tilde{\theta} -  \tilde{\theta}^* D^* \tilde{\theta}^*  &=& 0 . \label{inhomo-condition4}
\end{eqnarray}
In general, an $N=2$ superanalytic superfunction $H(z, \theta, \theta^*)$ {}from $\bigwedge_\infty^0 \oplus (\bigwedge_\infty^1)^2$ to $\bigwedge_\infty^0 \oplus (\bigwedge_\infty^1)^2$ transforms $D$ and $D^*$ by
\begin{eqnarray*}
D  &=& (D \tilde{\theta})\tilde{D} - (D \tilde{\theta}^* ) \tilde{D}^* + \left( D \tilde{z} - \tilde{\theta}D\tilde{\theta} -  \tilde{\theta}^*D \tilde{\theta}^*\right) D^2 \\
D^* &=& (D^* \tilde{\theta} )\tilde{D}  - (D^* \tilde{\theta}^* ) \tilde{D}^* + \left( D^* \tilde{z} - \tilde{\theta} D^* \tilde{\theta} - \tilde{\theta}^*  D^* \tilde{\theta}^*  \right) ( \tilde{D}^*)^2.
\end{eqnarray*}
Thus a homogeneous $N=2$ superconformal function is equivalent to a nonhomogeneous $N=2$ superanalytic superfunction transforming $D$, respectively $D^*$, as $D = (D \tilde{\theta}) \tilde{D} - (D\tilde{\theta}^* ) \tilde{D}^* =  (D \tilde{\theta}) \tilde{D} - (D^*\tilde{\theta} ) \tilde{D}^*$, respectively 
$D^* = (D^* \tilde{\theta} )\tilde{D}  - (D^* \tilde{\theta}^* ) \tilde{D}^* = (D^* \tilde{\theta} )\tilde{D}  + (D \tilde{\theta} ) \tilde{D}^*$.   We will call a nonhomogeneous $N=2$ superanalytic function satisfying  conditions (\ref{inhomo-condition1}) -- (\ref{inhomo-condition4}), a {\it nonhomogeneous $N=2$ superconformal function}.   Note that with this definition, a homogeneous $N=2$ superfunction $H(z, \theta^+, \theta^-)$ is superconformal if and only if the corresponding nonhomogeneous  $N=2$ superfunction $H(z, \theta, \theta^*)$ is superconformal in the above sense.

\begin{rema}\label{transform-D-remark}
{\em {}From the properties derived above for a nonhomogeneous $N=2$ superconformal function, we see one of the reasons for our terminology.  Namely,  that  a nonhomogeneous $N=2$ superconformal function does not transform the superderivations $D$ and $D^*$, respectively, homogeneously of degree one.  Instead it transforms them as  $D = (D \tilde{\theta}) \tilde{D} - (D\tilde{\theta}^* ) \tilde{D}^*$ and $D^* = (D^* \tilde{\theta} )\tilde{D}  - (D^* \tilde{\theta}^* ) \tilde{D}^*$, respectively -- unlike the homogeneous nature of the transformation of $D^\pm$ under a homogeneous $N=2$ superconformal function.  In the latter case the superderivations transform homogeneously as $D^\pm = (D^\pm \tilde{\theta}^\pm) \tilde{D}^\pm$.}
\end{rema}

The conditions (\ref{inhomo-condition1}) -- (\ref{inhomo-condition4}) imply that a nonhomogeneous $N=2$ superconformal function $H(z, \theta, \theta^*) = (\tilde{z}, \tilde{\theta}, \tilde{\theta}^*)$ is of the form
\begin{eqnarray}
\tilde{z} &=& f(z) + \theta (g(z) \psi(z) +   g^*(z)\psi^*(z))  + \theta^* (g(z)\psi^*(z) -g^*(z) \psi(z) )\label{nonhomo-superconformal-condition1}  \\
& &  \quad - \ \theta \theta^* ( \psi(z)  \psi^*(z))'   \nonumber \\
\tilde{\theta} &=& \psi (z) + \theta  g(z) - \theta^* g^*(z) + \theta \theta^*  (\psi^*)'(z)  \\
\qquad \tilde{\theta}^*&=& \psi^*(z) + \theta g^*(z) + \theta^* g(z) - \theta \theta^* \psi'(z) .
\end{eqnarray}
satisfying
\begin{equation}
f'(z) =   g^2(z) +   (g^*)^2(z)  -  \psi(z) \psi'(z)  -    \psi^*(z) (\psi^*)'(z)  , \label{nonhomo-superconformal-condition4} 
\end{equation}
for even superanalytic $(1,0)$-superfunctions $f, g$ and $g^*$ and odd superanalytic $(1,0)$-superfunctions $\psi, \psi^*$.  Note that, in the correspondence between homogeneous and nonhomogeneous $N=2$ superconformal functions given by (\ref{superconformal-condition1}) -- (\ref{superconformal-condition4}), and (\ref{nonhomo-superconformal-condition1}) -- (\ref{nonhomo-superconformal-condition4}), respectively, we have that $\psi^\pm (x) = \frac{1}{\sqrt{2}} (\psi(x) \pm i \psi^*(x))$ and $g^\pm(x) = g(x) \pm i g^*(x)$.

Under the transformation {}from homogeneous to nonhomogeneous coordinates, we see that the formal infinitesimal homogeneous $N=2$ superconformal transformations generated by the superderivations (\ref{L-notation}) -- (\ref{G-notation}) in $\mbox{Der} (\mathbb{C}  [x, x^{-1}, \varphi^+, \varphi^-])$ give the formal infinitesimal nonhomogeneous $N=2$ superconformal transformations as being generated by the following superderivations in $\mbox{Der} (\mathbb{C}  [x, x^{-1}, \varphi, \varphi^*])$
\begin{eqnarray}
L_j(x,\varphi,\varphi^*) &=& - \biggl( x^{j + 1} \frac{\partial}{\partial x} + (\frac{j + 1}{2})x^j \Bigl(  \varphi \frac{\partial }{\partial \varphi}  + \varphi^* \frac{\partial }{\partial \varphi^* }  \Bigr) \biggr) \label{L-nonhomo} \\
J_j(x,\varphi,\varphi^*) &=& i x^j \biggl( \varphi  \frac{\partial }{\partial \varphi^*} -   \varphi^*  \frac{\partial }{\partial \varphi} \biggr)    \\
G_{j -\frac{1}{2}} (x,\varphi,\varphi^*) &=& \frac{1}{\sqrt{2}} \left( G^+_{j - \frac{1}{2}} (x, \varphi^+, \varphi^-) + G^-_{j - \frac{1}{2}} (x, \varphi^+, \varphi^-) \right) \\
&=& - \biggl( x^j \Bigl(   \frac{\partial }{\partial \varphi}  -  \varphi  \frac{\partial}{\partial x} \Bigr)   - j x^{j-1} \varphi \varphi^*   \frac{\partial }{\partial \varphi^*}  \biggr)  \nonumber  \\
G^*_{j -\frac{1}{2}} (x,\varphi,\varphi^*) &=& - \frac{i}{\sqrt{2}} \left( G^+_{j - \frac{1}{2}} (x, \varphi^+, \varphi^-) - G^-_{j - \frac{1}{2}} (x, \varphi^+, \varphi^-) \right) \label{G-nonhom} \\
&=&  \biggl( x^j \Bigl( \frac{\partial}{\partial \varphi^*} - \varphi^* 
\frac{\partial}{\partial x} \Bigr)  + j x^{j-1}  \varphi \varphi^* \frac{\partial}{\partial \varphi} \biggr) \nonumber 
\end{eqnarray}
for formal nonhomogeneous variables $(x, \varphi, \varphi^*)$ and for $j \in \mathbb{Z}$.

These superderivations satisfy the relations for the ``nonhomogeneous $N=2$ Neveu-Schwarz algebra with central charge zero".  The nonhomogeneous $N=2$ Neveu-Schwarz algebra with central charge $d$ is the Lie superalgebra with basis consisting of the central charge $d$, even elements $L_n$ and $J_n$ and odd elements $G_{n + 1/2}$ and $G^*_{n + 1/2}$, for $n \in \mathbb{Z}$, and commutation relations
\begin{eqnarray}
\left[L_m ,L_n \right] &=& (m - n)L_{m + n} +  \frac{1}{12} (m^3 - m) \delta_{m + n, 0} \; d  , \label{Virasoro-nonhomo}\\
\left[J_m , J_n \right] &=& \frac{1}{3} m \delta_{m + n , 0} \;d, \\ 
\left[L_m , J_n \right] &=&  -nJ_{m+n} ,\\
\bigl[L_m,G_{n + \frac{1}{2}}\bigr] &=& \Bigl(\frac{m}{2} - n - \frac{1}{2} \Bigr) G _{m + n + \frac{1}{2}} ,  \\
\bigl[L_m,G^*_{n + \frac{1}{2}}\bigr] &=& \Bigl(\frac{m}{2} - n - \frac{1}{2} \Bigr) G^*_{m + n + \frac{1}{2}} ,  \\
\bigl[ J_m , G _{n + \frac{1}{2}}\bigr] &=& i G^*_{m + n + \frac{1}{2}},\\
\bigl[ J_m , G^*_{n + \frac{1}{2}}\bigr] &=& - i G_{m + n + \frac{1}{2}},\\
\bigl[ G_{m + \frac{1}{2}} , G_{n - \frac{1}{2}}\bigr] &=& 2L_{m + n} + \frac{1}{3}(m^2 + m) \delta_{m + n , 0} \;d , \\
\bigl[ G^*_{m + \frac{1}{2}} , G^*_{n - \frac{1}{2}} \bigr] &=& 2L_{m + n}  + \frac{1}{3}(m^2 + m) \delta_{m + n , 0} \;d , \\
\bigl[ G_{m + \frac{1}{2}} , G^*_{n - \frac{1}{2}} \bigr] &=& i(m-n+1) J_{m+n} \label{NS-nonhomo}. 
\end{eqnarray}
Note that in general, if $L(n), J(n), G^\pm(n + 1/2)$ give a representation of the (homogeneous) $N=2$ Neveu-Schwarz algebra with central charge $c$ (i.e., satisfy the relations (\ref{Virasoro-relation}) -- (\ref{NS-relation6})), then letting
\begin{eqnarray}
G(n+ \frac{1}{2}) &=& \frac{1}{\sqrt{2}} \left( G^+(n + \frac{1}{2}) + G^-(n + \frac{1}{2}) \right) \\
G^*(n+ \frac{1}{2}) &=& \frac{-i}{\sqrt{2}} \left( G^+(n + \frac{1}{2}) - G^-(n + \frac{1}{2}) \right),
\end{eqnarray}
(or equivalently $G^\pm (n + 1/2) = \frac{1}{\sqrt{2}} ( G(n + 1/2) \pm i G^*(n + 1/2))$) we have that $L(n), J(n), G(n + 1/2), G^*(n + 1/2)$ give a representation of the nonhomogeneous $N=2$ Neveu-Schwarz algebra with central charge $c$. 

\begin{rema}\label{homo-infinitesimals-remark}
{\em  We see in the commutation relations (\ref{Virasoro-relation}) -- (\ref{NS-relation6}) for the (homogeneous) $N=2$ Neveu-Schwarz algebra in comparison to the commutation relations (\ref{Virasoro-nonhomo}) -- (\ref{NS-nonhomo}) for the nonhomogeneous $N=2$ Neveu-Schwarz algebra another justification for our terminology ``homogeneous" and ``nonhomogeneous".   In the former case, the commutator of $J_m$ with $G^\pm_{n + 1/2}$ is homogeneous in $G^+$ and $G^-$ terms, respectively, whereas in the later case, the commutator of $J_m$ with $G_{n + 1/2}$ and with $G^*_{n + 1/2}$, respectively, is nonhomogeneous. }
\end{rema} 

Carrying through our transformation {}from homogeneous coordinates to nonhomogeneous coordinates we see that the $N=2$ super-Riemann sphere in nonhomogeneous coordinates is given by the usual coordinate charts $(U_\sou, \sou)$ and $(U_\nor, \nor)$, and coordinate transition function $\sou \circ \nor^{-1} (w, \rho, \rho^*) = I(w, \rho, \rho^*) = (1/w, i\rho/w, i\rho^*/w)$.   Local coordinates vanishing at a point $p$ on the nonhomogeneous $N=2$ super-Riemann sphere with $\sou(p) = (z, \theta, \theta^*) \in \bigwedge_\infty^0 \oplus (\bigwedge_\infty^1)^2$ are given by 
\begin{multline}
H (w, \rho, \rho^*) =  \exp \Biggl(\! - \! \sum_{j \in \Z} \Bigl( A_j L_j(x,\varphi, \varphi^*) + A^*_j J_j(x,\varphi, \varphi^*)  \\ 
+  M_{j - \frac{1}{2}} G_{j -\frac{1}{2}} (x, \varphi, \varphi^*) +  M^*_{j - \frac{1}{2}} G^*_{j -\frac{1}{2}} (x, \varphi, \varphi^*) \Bigr) \! \Biggr)\!  \cdot  \! a_0^{-2L_0(x,\varphi, \varphi^*)} \cdot   \\
\left.  (a_0^*)^{-J_0(x,\varphi, \varphi^*)} \cdot (x, \varphi, \varphi^*) \right|_{(x,\varphi, \varphi^*) = (w - z - \rho \theta - \rho^* \theta^*, \rho - \theta,  \rho^* - \theta^*)}
\end{multline}
for $(a_0 , a_0^*) \in ((\bigwedge_\infty^0)^\times)^2/\langle \pm 1 \rangle$, $A_j, A^*_j \in \bigwedge_\infty^0$, and $M_{j- 1/2}, M^*_{j-1/2} \in \bigwedge_\infty^1$, for $j \in \Z$.  And we can write the power series expansion of the local coordinate vanishing at $\infty$ with leading coefficient of $\rho$ in $\tilde{\rho}$ and $\rho^*$ in $\tilde{\rho}^*$, respectively, equal to $i$ as  
\begin{multline}
H (w, \rho, \rho^*) = \exp \Biggl(\sum_{j \in \Z} \Bigl( B_j L_{-j}(w,\rho, \rho^*) + B^*_j J_{-j}(w,\rho, \rho^*)  \\ 
+  N_{j - \frac{1}{2}} G_{-j +\frac{1}{2}} (w, \rho, \rho^*) +  N^*_{j - \frac{1}{2}} G^*_{-j +\frac{1}{2}} (w, \rho, \rho^*) \Bigr) \! \Biggr) \cdot \Bigl(\frac{1}{w}, \frac{i \rho}{w}, \frac{i \rho^*}{w}\Bigr) 
\end{multline}
for $B_j, B^*_j \in \bigwedge_\infty^0$, and $N_{j- 1/2}, N^*_{j-1/2} \in \bigwedge_\infty^1$, for $j \in \Z$.
This characterization allows us to formulate the moduli space of nonhomogeneous $N=2$ super-Riemann spheres with tubes in a way completely analogous to that for homogeneous $N=2$ super-Riemann spheres with tubes as is done in Section \ref{reformulate-moduli-section}.

Now, we specifically formulate the Lie supergroup of nonhomogeneous $N=2$ superprojective transformations {}from the homogeneous $N=2$ superprojective transformations given by  (\ref{transform1}) -- (\ref{transform-condition4}) (or equivalently (\ref{short-transform1}) -- (\ref{short-transform-condition2})) and discussed in detail in Section \ref{NS-algebra-section}.  Under the transform {}from homogeneous coordinates to nonhomogeneous coordinates, we see that the automorphism group for the nonhomogeneous $N=2$ super-Riemann sphere is given by transformations $(w, \rho, \rho^*) \mapsto (\tilde{w}, \tilde{\rho}, \tilde{\rho}^*)$ of the form
\begin{eqnarray}
\tilde{w} &=& \frac{aw + b}{cw + d} +   \rho \frac{  e (\gamma w +  \delta)  + e^* (\gamma^*  w +  \delta^*) }{(cw+d)^2}  \label{nonhomo-transform1}\\
& & +    \rho \frac{  (f w + h) (\gamma w    + \delta)   +  (f^*w + h^*) ( \gamma^* w  +   \delta^*)   }{(cw+d)^3} \nonumber\\
& & +  \rho^* \frac{ e (\gamma^*w  + \delta^*) - e^*( \gamma w +  \delta)  }{(cw+d)^2}  \nonumber  \\
& & + \rho^* \frac{ (f w + h)( \gamma^*w   +  \delta^*)  - ( f^*w + h^*)  (\gamma w  +  \delta)    }{(cw+d)^3} \nonumber \\
& & -\rho\rho^* \frac{ 2\gamma \gamma^*  dw  - (  \gamma \delta^* - \gamma^* \delta  )(cw - d) - 2  \delta \delta^*c}{(cw + d)^3} \nonumber \\
\tilde{\rho} &=& \frac{\gamma w + \delta}{cw + d}  + \rho  \frac{e}{cw+d}  +  \rho \frac{fw+ h}{(cw+d)^2}   -  \rho^*   \frac{e^*}{cw+d}\\
& & - \rho^* \frac{f^* w+ h^*}{(cw+d)^2}   + \rho  \rho^*  \frac{\gamma^* d - \delta^*  c}{(cw+d)^2}  \nonumber  \\
\tilde{\rho}^*&=& \frac{\gamma^* w + \delta^* }{cw + d}  + \rho \frac{e^* }{cw+d}  + \rho \frac{f^* w+ h^* }{(cw+d)^2}  +  \rho^*   \frac{e}{cw+d} \\
& & + \rho^* \frac{f w+ h}{(cw+d)^2}     -  \rho  \rho^*  \frac{\gamma d -  \delta  c}{(cw+d)^2} \nonumber
\end{eqnarray}
with $a,b,c,d,e,e^*, f, f^*, h, h^* \in \bigwedge_\infty^0$ and $\gamma, \gamma^*, \delta, \delta^* \in \bigwedge_\infty^1$ satisfying
\begin{eqnarray}
ad - bc &= & 1 \\
e^2 + (e^*)^2  &=& 1-  \gamma  \delta   + \delta^* \gamma^*  \\
f  &=&  - e^* \gamma \gamma^*  d  \\
f^*  &=& e  \gamma  \gamma^* d  \\
h &=& e^* (\delta  \delta^*c  - ( \gamma \delta^* -  \gamma^* \delta   )d)   + e\delta \delta^* \gamma \gamma^*  d \\
h^* &=& - e (  \delta \delta^* c -  (\gamma \delta^* - \gamma^* \delta  )d)  + e^* \delta \delta^*\gamma \gamma^* d . \label{nonhomo-condition6}
\end{eqnarray} 

Note that, in the correspondence between homogeneous and nonhomogeneous $N=2$ superprojective transformations given by (\ref{transform1}) -- (\ref{transform-condition4}) and (\ref{nonhomo-transform1}) -- (\ref{nonhomo-condition6}), respectively, we have that 
\begin{eqnarray*}
\gamma &= \frac{1}{\sqrt{2}} (\gamma^+ + \gamma^-) \qquad &\gamma^*  =  - \frac{i}{\sqrt{2}} (\gamma^+ - \gamma^-) \\
\delta &=  \frac{1}{\sqrt{2}} (\delta^+ + \delta^-)  \qquad & \delta^* =  - \frac{i}{\sqrt{2}} (\delta^+ - \delta^-) \\
e &=  \frac{1}{2} (e^+ + e^-)  \qquad   &  e^* = - \frac{i}{2} (e^+ - e^-) \\
f &= \frac{1}{2} (f^+ + f^-) \qquad & f^* = - \frac{i}{2} (f^+ - f^-) \\
h &= \frac{1}{2} (h^+ + h^-) \qquad & h^* = - \frac{i}{2} (h^+ - h^-).
\end{eqnarray*}

Equivalently, the nonhomogeneous $N=2$ superprojective transformations can be written as
\begin{eqnarray}
\tilde{w} &=& \frac{aw + b}{cw + d} +   \rho \frac{  e (\gamma w +  \delta)  + e^* (\gamma^*  w +  \delta^*) + e^* \delta  \delta^* \gamma     - e   \delta \delta^* \gamma^*}{(cw+d)^2} \label{short-nonhomo-trans1}  \\
& & +  \rho^* \frac{ e (\gamma^*w  + \delta^*) - e^*( \gamma w +  \delta) + e^* \delta  \delta^*  \gamma^*  +  e   \delta \delta^* \gamma     }{(cw+d)^2 }   \nonumber \\
& & -\rho\rho^* \frac{ 2\gamma \gamma^*  dw  - (  \gamma \delta^* - \gamma^* \delta  )(cw - d) - 2  \delta \delta^*c}{(cw + d)^3} \nonumber  \\
\tilde{\rho} &=& \frac{\gamma w + \delta}{cw + d}  + \rho  \frac{e}{cw+d}  \\
& & +  \rho \frac{ - e^* \gamma \gamma^*  d w+ e^* (\delta  \delta^*c  - ( \gamma \delta^* -  \gamma^* \delta   )d)   + e\delta \delta^* \gamma \gamma^*  d }{(cw+d)^2}  \nonumber\\
& &    -  \rho^*  \frac{e^*}{cw+d} \nonumber\\
& & - \rho^*  \frac{e  \gamma  \gamma^* d  w - e (  \delta \delta^* c -  (\gamma \delta^* - \gamma^* \delta  )d)  + e^* \delta \delta^*\gamma \gamma^* d }{(cw+d)^2}     \nonumber\\
&  & + \rho  \rho^*  \frac{\gamma^* d - \delta^*  c}{(cw+d)^2}  \nonumber  \\
\tilde{\rho}^*&=& \frac{\gamma^* w + \delta^* }{cw + d}  + \rho  \frac{e^* }{cw+d}  \\
& & + \rho  \frac{e  \gamma  \gamma^* d  w - e (  \delta \delta^* c -  (\gamma \delta^* - \gamma^* \delta  )d)  + e^* \delta \delta^*\gamma \gamma^* d  }{(cw+d)^2} \nonumber   \\
&  &+  \rho^*  \frac{e}{cw+d} \nonumber\\
& & + \rho^*  \frac{ - e^* \gamma \gamma^*  d w+ e^* (\delta  \delta^*c  - ( \gamma \delta^* -  \gamma^* \delta   )d)   + e\delta \delta^* \gamma \gamma^*  d }{(cw+d)^2}   \nonumber\\
& &  -  \rho  \rho^*  \frac{\gamma d -  \delta  c}{(cw+d)^2} \nonumber
\end{eqnarray}
with $a,b,c,d,e,e^* \in \bigwedge_\infty^0$ and $\gamma, \gamma^*, \delta, \delta^* \in \bigwedge_\infty^1$ satisfying
\begin{eqnarray}
ad-bc&=& 1 \\
e^2 + (e^*)^2  &=& 1-  \gamma  \delta   + \delta^* \gamma^*  . \label{short-nonhomo-trans-condition2}
\end{eqnarray} 

\begin{rema} {\em 
The nonhomogeneous $N=2$ superprojective transformations above do not coincide with those given in \cite{Melzer} or \cite{Sc}.  Those given in \cite{Melzer}  do not give all possible transformations.  An example of a nonhomogeneous $N=2$ superprojective transformation which is of the form (\ref{short-nonhomo-trans1}) -- (\ref{short-nonhomo-trans-condition2}) but cannot in general be put in the form proposed in \cite{Melzer} can be obtained by transforming the homogeneous $N=2$ superprojective transformation given in Example \ref{ex1} via the coordinate transformation {}from homogeneous to nonhomogeneous coordinates (\ref{transform-homo-inhomo}).   On the other hand, in \cite{Sc}, Schoutens' presentation of the nonhomogeneous $N=2$ superprojective transformations implies that a transformation $(w, \rho, \rho^*) \mapsto (\tilde{w} , \tilde{\rho}, \tilde{\rho}^*)$ of the form
\begin{eqnarray}
\qquad (\tilde{w} , \tilde{\rho}, \tilde{\rho}^*) &=& \left( \frac{aw+b}{cw+d}, \ \rho \frac{t^{1 1}}{cw +d} + \rho^* \frac{ t^{12}}{cw+d},  \ \rho \frac{t^{2 1}}{cw +d} + \rho^* \frac{ t^{22}}{cw+d} \right)
\end{eqnarray}
with $a,b,c,d, t^{jk} \in \bigwedge_\infty^0$ for $j,k = 1,2$ satisfying $ad-bc = 1$ and $t^{ij} t^{ik} = \delta_{i,k}$ for $i,j,k = 1,2$, is $N=2$ superprojective.   (Here for simplicity we have set Schoutens' odd parameters $\varepsilon^i_j = 0$ for $i, j = 1,2$, as of course we are allowed to do.)  But clearly, in general, such a transformation is not of the form (\ref{short-nonhomo-trans1}) -- (\ref{short-nonhomo-trans-condition2}), and in fact is not even a nonhomogeneous $N=2$ superconformal function since it does not satisfy (\ref{nonhomo-superconformal-condition1}) -- (\ref{nonhomo-superconformal-condition4}). 
}
\end{rema}

\end{document}